\documentclass[12pt,russian]{article}
\usepackage[T2A]{fontenc}
\usepackage[russian,english]{babel}
\usepackage{amssymb,amsfonts,amsmath,amsthm}
\usepackage{citehack}

\theoremstyle{plain}

\newtheorem{savkoz}{Савчук -- чудак!!!}[section]
\newtheorem{Lemm}[savkoz]{Лемма}
\newtheorem{Proposition}[savkoz]{Предложение}
\newtheorem{Theorem}[savkoz]{Теорема}

\newtheorem{Note}[savkoz]{Замечание}

\newcommand{\intl}{\mathop{\int}\limits}

\newcommand{\suml}{\mathop{\sum}\limits}
\renewcommand{\Re}{\operatorname{Re}}
\renewcommand{\Im}{\operatorname{Im}}

\renewcommand{\le}{\leqslant}
\renewcommand{\ge}{\geqslant}

\def\Dom{{\frak D}}

\def\al{\alpha}
\def\la{\lambda}

\def\Lal{\Lambda_\alpha}

\def\tg{\operatorname{tg}}
\def\eps{\varepsilon}

\def\sign{\operatorname{sign}}

\def\lar{\mathop{\longrightarrow}\limits}

\def\R{{\mathbb R}}
\def\N{{\mathbb N}}

\def\C{{\mathbb C}}
\def\z{{\mathbf z}}
\def\y{{\mathbf y}}

\def\D{{\mathcal D}}
\def\A{{\mathbf A}}
\def\f{{\mathbf f}}
\def\xxi{{\mathbf \xi}}
\def\b{B}
\def\a{A}

\def\w{{\rm w}}
\def\wf{\widetilde{f}}
\newlength{\lenun}
\newlength{\lendu}
\def\Wo{\makebox{\parbox[c][\lendu][b]{\lenun}{\(\stackrel{o}{W}\)}}}

\makeatletter
\@addtoreset{equation}{section}
\makeatother

\begin{document}

{\Large
{\bf \centerline{Операторы Штурма--Лиувилля}}
{\bf \centerline{ с потенциалами --- распределениями}}}\footnote{Работа
выполнена при поддержке Российского Фонда Фундаментальных Исследований.
А.М.Савчук поддержан грантом No. 01-01-00691, А.А.Шкаликов поддержан грантом
No. 01-01-00958. }

\bigskip

{\centerline{Савчук А.~М., Шкаликов А.~А.}}

\medskip
{\footnotesize В этой статье предлагаются четыре различных подхода к
определению оператора Штурма--Лиувилля \(Ly=-y''+q(x)y\) на интервале \((a,b)\)
в случае, когда потенциал \(q(x)\) является распределением из соболевского
пространства с негативным индексом гладкости, а именно, \(q\in W_2^{-\theta}\)
при \(\theta\le 1\). Получены главные и вторые члены асимптотик для
собственных значений и собственных функций определяемых операторов и проведены
оценки остатков в этих асимптотических формулах в зависимости от класса
потенциала.  Соответствующие операторы определяются и изучаются также для
некоторых аналитических семейств потенциалов высокой сингулярности, не
принадлежащих пространству \(W_2^{-1}\).}

\medskip

{\bf\centerline{Оглавление}}
\noindent
0. Введение\hfill\break
1. Четыре подхода к определению операторов с
потенциалами-распределениями\hfill\break
2. Асимптотика собственных значений и собственных функций для регулярных
краевых условий\hfill\break
3. Вторые члены в асимптотиках для собственных значений и собственных
функций\hfill\break
4. Операторы Штурма--Лиувилля с потенциалами высокой сингулярности\hfill\break
\vskip 1cm

\section*{{\bf \S 0.\ \  Введение}}\refstepcounter{section}

В классической теории операторов Штурма--Лиувилля, порождаемых на интервале
\((a,b)\subset\mathbb R\) дифференциальным выражением
\begin{equation}\label{de1}
    l(y)=-y''+q(x)y,
\end{equation}
стандартным условием на функцию \(q(x)\) является условие \(q(x)\in
L_{1,loc}(a,b)\), т.~е. функция предполагается суммируемой на любом отрезке,
компактно вложенном в \((a,b)\). Сингулярные операторы Штурма--Лиувилля в
классической теории характеризуются тем, что либо функция \(q(x)\) не
суммируема на отрезке \([a,b]\) (имеется неинтегрируемая особенность по
крайней мере на одном из концов отрезка), либо интервал \((a,b)\) бесконечен.
В этой работе мы изучаем операторы с потенциалами, имеющими неинтегрируемые
особенности внутри интервала. При этом естественным оказывается язык теории
распределений, в частности, естественно рассматривать потенциалы из
пространств Соболева с негативными индексами гладкости.

Задачи об изучении оператора Штурма--Лиувилля и его многомерных аналогов
\(-\Delta+q(x)\) с потенциалом короткого взаимодействия (типа
\(\delta\)-функции) возникли в физической литературе. Математические
исследования соответствующих физических моделей были инициированы в начале
60-х годов в работах Березина, Фаддеева и Минлоса~\cite{B}, \cite{BF},
\cite{MF}. Эта тематика интенсивно развивалась в последние два десятилетия.
Имеются монографии Альбеверио, Гештези, Хоэг-Крона и Хольдена~\cite{AGHH},
Кошманенко~\cite{Ko}, Альбеверио и Курасова~\cite{AK}, где можно познакомиться
с подробностями теории Березина--Минлоса--Фаддеева в ее современном
состоянии и другими новыми направлениями, возникшими на основе этой
теории.  Там же можно познакомиться с обширной библиографией.

Другой подход к изучению операторов Штурма--Лиувилля с неклассическими
потенциалами \(q(x)\), являющимися производными от функций
ограниченной вариации (зарядами), был предпринят Крейном~\cite{Kr1},
Кацем~\cite{Kac}, Аткинсоном~\cite{At} и Жиковым~\cite{Zhi}.  На этом
пути в недавней работе Винокурова и Садовничего \cite{ВС} получены
асимптотические формулы для собственных значений и собственных функций
такого класса операторов.  Из потенциалов, не принадлежащих последнему
классу, изучался кулоновский потенциал \(q(x)=1/x\) на отрезке \([-1,1]\)
или на прямой \(\mathbb R\), например, в работах Гунсона~\cite{Gu},
Курасова~\cite{Ku}, Аткинсона, Эверитта и Зеттла~\cite{AEZ}.

В работе авторов~\cite{SSh1} (см. также работу Нейман-заде и
Шкаликова~\cite{NS}) было показано, что оператор Штурма--Лиувилля можно
корректно определить для существенно более общего класса потенциалов \(q(x)\),
являющихся сингулярными распределениями первого порядка. Авторы предприняли
дальнейшее изучение операторов с такими потенциалами в работах~\cite{Chud},
\cite{diss} и \cite{SSh2}. Вскоре появились работы Гринива и
Микитюка~\cite{HM1}--\cite{HM5}, где этот подход получил существенное
развитие, в особенности при решении обратной задачи Штурма--Лиувилля с
неклассическими потенциалами. Отметим также недавнюю работу Каппелера и
Моора~\cite{KM}, в которой предложен другой подход к изучению операторов с
потенциалами-распределениями из пространств \(W_2^{-\theta}(a,b)\) при
\(\theta<1\).

Настоящая статья является продолжением и существенным развитием работ
авторов~\cite{SSh1}, \cite{Chud}, \cite{diss}. Мы ставим три цели: 1)~дать
несколько подходов к определению оператора Штурма--Лиувилля с
потенциалом-распределением и обсудить их взаимосвязь; 2)~получить
асимптотические формулы для собственных значений и собственных функций
определяемых операторов, выписывая по возможности вторые члены асимптотик и
оценивая остатки в зависимости от классов потенциалов; 3)~предложить
подходы к определению операторов с потенциалами высокой сингулярности
(для потенциалов \(q(x)\not\in W_2^{-1}\)), когда однозначного
определения оператора с помощью четырех ранее предложенных подходов не
существует. В последнем случае мы только намечаем подходы, которые, как
мы надеемся, должны получить дальнейшее развитие.

\medskip

\section*{{\bfseries \S 1.\ \ Четыре подхода к определению операторов с
потенциалами-распределениями}}\refstepcounter{section}
\textbf{1.1. Метод регуляризации.} Далее через \(\mathcal D\)
обозначается пространство тест-функций на интервале \((0,\pi)\) (т.~е.
бесконечно дифференцируемых функций с компактным носителем на
\((0,\pi)\)), а через \(\mathcal D'\) --- пространство распределений на
\(\mathcal D\). Через \(W^{-1}_2[0,\pi]\) (сокращенно \(W^{-1}_2\))
обозначаем пространство, состоящее из функций \(q(x)\in\mathcal D'\),
для которых первообразная \(u(x)=\int q(\xi)\,d\xi\) (в смысле
распределений) принадлежит \(L_2[0,\pi]\). Норму в \(W^{-1}_2\) определяем
равенством \(\|q\|_{-1}=\inf\|u(x)+c\|_{L_2}\), где \(\inf\) берется по всем
константам \(c\). Нетрудно показать, что пространство \(W_2^{-1}\)
совпадает с дуальным к пространству \(\Wo_2^1[0,\pi]\) по отношению к
скалярному произведению в \(L_2[0,\pi]\). Здесь
\[
    \Wo_2^1[0,\pi]=\left\{y\;\vline\;y\in W_2^1[0,\pi],\;y(0)=y(\pi)=0
    \right\},
\]
где через \(W_p^k\) обозначается соболевское пространство с нормой
\(\|y\|_{k,p}=\|y\|_{L_p}+\|y^{(k)}\|_{L_p}\). Далее,
если норма \(\|\cdot\|\) пишется без индексов, предполагается, что она берется
в пространстве \(L_2\).

Пусть в дифференциальном выражении~\eqref{de1} \(q(x)\in W_2^{-1}\), а
\(u(x)=\int q(\xi)\,d\xi\) --- первообразная из пространства \(L_2\). Введем
квазипроизводную
\[
    y^{[1]}(x)=y'(x)-u(x)y(x)
\]
и перепишем выражение~\eqref{de1} в виде
\begin{equation}\label{eq:1.1}
    l(y)=-\left(y^{[1]}\right)'-u(x)y^{[1]}-u^2(x)y.
\end{equation}
Несложно видеть, что для гладкой функции \(u(x)\) дифференциальное
выражение~\eqref{de1} и квазидифференциальное выражение~\eqref{eq:1.1}
совпадают. Однако выражение~\eqref{eq:1.1} обладает тем преимуществом, что не
содержит распределений, а потому с ним можно оперировать, по существу,
так же, как в классической теории. Методы построения операторов на основе
квазидифференциальных выражений можно найти в монографии Наймарка~\cite{Naj},
статьях Эверитта, Маркуса и Зеттла~\cite{E}, \cite{EM}, \cite{EZ}. Мы
будем пользоваться конструкцией, приведенной в работе
авторов~\cite{SSh1}. Для удобства читателей мы напомним эту конструкцию,
опуская доказательства, которые можно найти в~\cite{SSh1}.

С выражением~\eqref{eq:1.1} свяжем максимальный оператор \(L_M\), определенный
равенствами
\begin{equation}\label{eq:Lm}
    \left\{
        \begin{aligned}
        L_My&=l(y),\\
        \mathfrak D(L_M)&=\left\{y|\ y,y^{[1]}\in W_1^1[0,\pi],\,
        l(y)\in L_2[0,\pi]\right\},
    \end{aligned}\right.
\end{equation}
а также минимальный оператор \(L_m\), являющийся сужением максимального
оператора на область
\[
    \mathfrak D(L_m)=\left\{y|\ y\in\mathfrak D(L_M),
    y(0)=y(\pi)=y^{[1]}(0)=
    y^{[1]}(\pi)=0\right\}.
\]
Заметим, что в случае $u'(x)\in L_1$, определения максимального и
минимального операторов совпадают с классическим их определениями, т.~к.
в этом случае условие $y^{[1]}\in W_1^1$ влечет $y'\in W_1^1$, и
наоборот.

Мы не предполагаем, что функция \(u(x)\) вещественна. Через \(\overline{L_M}\)
и \(\overline{L_m}\) будем обозначать максимальный и минимальный операторы,
порожденные сопряженным дифференциальным выражением \(\overline{l}(y)\) (в
котором функция \(u(x)\) заменена на \(\overline{u(x)}\)). Прямым вычислением
получаем следующее предложение.

\begin{Lemm}[Формула Лагранжа]\label{lem:1.1} Для функций \(f\in\mathfrak
D(L_M)\), \(g\in\mathfrak D(\overline{L_M})\) справедливо тождество
\[
    (L_Mf,g)=(f,\overline{L_M}g)+[f,g]^{\pi}_0,
\]
где \([f,g]_0^{\pi}=-\left.f^{[1]}(x)\overline{g(x)}\right|_0^{\pi}-
\left.\overline{g^{[1]}(x)}f(x)\right|_0^{\pi}\).
\end{Lemm}

Из этой формулы, в частности, получаем
\begin{align*}
    (L_Mf,g)&=(f,\overline{L_m}g),&f&\in\mathfrak D(L_M),\;g\in
    \mathfrak D(\overline{L_m}),
\end{align*}
т.~е. операторы \(L_M\) и \(\overline{L_m}\) взаимно сопряженны.

Заметим, что уравнение
\begin{align*}
    L_My&:=-y''+u'(x)y=\lambda y+f,&\lambda\in\mathbb C,\;f\in L_2,
\end{align*}
можно записать в виде системы
\begin{equation}\label{eq:S1}
    \begin{pmatrix}y_1\\y_2\end{pmatrix}=\begin{pmatrix}u&1\\-\lambda-u^2
    &-u\end{pmatrix}\begin{pmatrix}y_1\\y_2\end{pmatrix}+\begin{pmatrix}
    0\\f\end{pmatrix},
\end{equation}
где \(y_1=y\), \(y_2=y^{[1]}\). При этом элементы матрицы
\[
    \mathbf A(x)=\begin{pmatrix}u&1\\-\lambda-u^2&-u\end{pmatrix}
\]
являются функциями из \(L_1[0,\pi]\). Это обстоятельство позволяет использовать
следующее важное утверждение, доказательство первой части которого можно найти
в~\cite{Naj}, а второй части --- в~\cite{SSh1}.
\begin{Theorem}\label{tm:0} Пусть $\A(x)$ --- матрица размера $n\times n$,
элементы которой являются функциями пространства $L_1[0,\pi]$, а
$\f\in \big [ L_1(0,\pi)\big ]^n$ --- вектор-функция. Тогда при любом
$c\in[0,\pi]$ уравнение
\begin{align*}
\y'&=\A(x)\y+\f,& \y(c)&=\xxi\in\C^n,
\end{align*}
имеет единственное решение $\y(x)$, причем $\y(x)$ --- абсолютно непрерывная
на $[0,\pi]$ вектор--функция. Если последовательность матриц
$\A_\varepsilon(x)$ с элементами из $L_1[0,\pi]$ такова, что
$\|\A_\varepsilon(x)-\A(x)\|_{L_1}\to 0$ при $\varepsilon\to 0$, то
решения уравнений
\begin{align*}
    \y_\varepsilon'&=\A_\varepsilon(x)\y_\varepsilon+\f,&
    \y_\varepsilon(c)=\xxi,
\end{align*}
сходятся к $\y(x)$ равномерно на $[0,\pi]$ (и даже в метрике пространства
$W_1^1[0,\pi]$).  Кроме того, справедлива оценка
\begin{equation*}
    \|\y(x)-\y_\varepsilon(x)\|_{1,1}\le
    C\|\f\|_{L_1}\|\A(x)-\A_\varepsilon(x)\|_{L_1}
\end{equation*}
с постоянной $C$, не зависящей от $\f$ и $\varepsilon$.
\end{Theorem}
Напомним, что оператор \(F\), действующий в гильбертовом (или банаховом)
пространстве \(\mathfrak H\), называется фредгольмовым, если его область
определения \(\mathfrak D(F)\) плотна в \(\mathfrak H\), образ замкнут, а
дефектные числа \(\{\alpha,\beta\}\), равные размерностям ядра и коядра,
конечны. Теорема~\ref{tm:0} позволяет легко доказать следующее утверждение
(см.~\cite{SSh1}).

\begin{Theorem}\label{tm:1.2} При любом \(\lambda\in\mathbb C\) операторы
\(L_M-\lambda\) и \(\overline{L_m}-\overline{\lambda}\) фредгольмовы, являются
сопряженными друг к другу, а их дефектные числа равны \(\{0,2\}\) и
\(\{2,0\}\), соответственно.
\end{Theorem}

Далее нам понадобится следующее утверждение, которое мы докажем ниже, в \S\,2.
\begin{Lemm}\label{lem:1.3}
Пусть \(\varphi(x,\lambda)\) и \(\psi(x,\lambda)\) --- решения уравнения
\(l(y)=\lambda y\) с начальными условиями
\begin{equation}\label{eq:1.4+}
    \varphi(0,\lambda)=1,\quad\varphi^{[1]}(0,\lambda)=0,\quad
    \psi(0,\lambda)=0,\quad\psi^{[1]}(0,\lambda)=1.
\end{equation}
Тогда при \(\lambda\to+\infty\) внутри любой параболы
\[
    P_a=\left\{\lambda\in\mathbb C\;\vline\;|\Im\sqrt{\lambda}|\le a\right\}
\]
справедливы асимптотические представления
\begin{equation}\label{zz}
\begin{array}{cc}
\Phi(x,\lambda)=\cos(\lambda^{1/2}x)+\varphi(x,\la),\quad
&\Phi^{[1]}(x,\lambda)=-\lambda^{1/2}\sin(\lambda^{1/2}x)
+\lambda^{1/2}\varphi_1(x,\la),\\
\Psi(x,\lambda)=\lambda^{-1/2}\sin(\lambda^{1/2}x)+
\lambda^{-1/2}\psi(x,\la), \quad
&\Psi^{[1]}(x,\lambda)=\cos(\lambda^{1/2}x)+\psi_1(x,\la),
\end{array}
\end{equation}
где каждая из функций $\varphi$, $\varphi_1$, $\psi$, $\psi_1$
мажорируется величиной $M\Upsilon(\la)\to0$ при $\la\to\infty$,
где $\Upsilon(\la)=\Upsilon(\pi/2,\la)$.
\end{Lemm}
Заметим, что всякий оператор \(L\), подчиненный условию \(L_m\subset L\subset
L_M\), имеет область
\[
    \mathfrak D(L)=\left\{y\;\vline\;y\in\mathfrak D(L_M),\;U_j(y)=0,
    1\le j\le\nu\right\},
\]
где \(U_j(y)\) --- линейные формы от переменных \(y(0)\), \(y(\pi)\),
\(y^{[1]}(0)\), \(y^{[1]}(\pi)\). Эти формы можно считать линейно
независимыми, и тогда их число \(\nu\) заключено между \(0\) и \(4\). Из
теоремы~\ref{tm:1.2} следует, что нужно брать две линейные формы, если хотим,
чтобы оператор \(L\) имел непустое резольвентное множество (для этого
необходимо, чтобы индексы \(L\) были равны \(\{0,0\}\)).

Теперь мы готовы доказать следующий результат.
\begin{Theorem}\label{tm:1.4}
Пусть оператор \(L\) является сужением оператора \(L_M\) на область
\[
    \mathfrak D(L)=\left\{y\;\vline\;y\in\mathfrak D(L_M),\;U_1(y)=
    U_2(y)=0\right\},
\]
где
\begin{align}\label{eq:T.1.4}
    U_j(y)&=a_{j1}y(0)+a_{j2}y^{[1]}(0)+b_{j1}y(\pi)+b_{j2}y^{[1]}(\pi),&
    j&=1,2.
\end{align}
Обозначим через \(J_{\alpha\beta}\) определитель, составленный из \(\alpha\)-го
и \(\beta\)-го столбца матрицы
\[
    \begin{pmatrix}a_{11}&a_{12}&b_{11}&b_{12}\\
    a_{21}&a_{22}&b_{21}&b_{22}\end{pmatrix}.
\]
Тогда оператор \(L\) имеет непустое резольвентное множество и спектр его
дискретен, если выполнено одно из следующих условий:
\begin{enumerate}
\item \(J_{42}\neq 0\),
\item \(J_{42}=0\), \(J_{14}+J_{32}\neq 0\),
\item \(J_{42}=J_{14}=J_{32}=0\), \(J_{12}+J_{34}=0\),
\(J_{13}\neq 0\).
\end{enumerate}
\end{Theorem}
\begin{proof}
Пусть \(\varphi(x,\lambda)\), \(\psi(x,\lambda)\) --- решения уравнения
\(l(y)=\lambda y\), определенные в лемме~\ref{lem:1.3}. Согласно
теореме~\ref{tm:0}, для любой функции \(f\in L_2\) существует решение
\(z(x,\lambda)\) уравнения \(l(y)-\lambda y=f\), подчиненное условию
\(z(0,\lambda)=z^{[1]}(0,\lambda)=0\). Общее решение этого уравнения
имеет вид \(y=c_1\varphi+c_2\psi+z\), где \(c_1\), \(c_2\) ---
постоянные.

Условие \(y\in\mathfrak D(L)\) влечет \(U_j(y)=0\), \(j=1,2\). Поэтому решение
\(y\) находится однозначно, если не равен нулю характеристический определитель
\[
    \Delta(\lambda)=\begin{vmatrix}U_1(\varphi)&U_1(\psi)\\
    U_2(\varphi)&U_2(\psi)\end{vmatrix}=\begin{vmatrix}a_{11}+b_{11}
    \varphi(\pi,\lambda)+b_{12}\varphi^{[1]}(\pi,\lambda)&a_{12}+b_{11}
    \psi(\pi,\lambda)+b_{12}\psi^{[1]}(\pi,\lambda)\\
    a_{21}+b_{21}\varphi(\pi,\lambda)+b_{22}\varphi^{[1]}(\pi,\lambda)&
    a_{22}+b_{21}\psi(\pi,\lambda)+b_{22}\psi^{[1]}(\pi,\lambda)
    \end{vmatrix}.
\]
Прямые вычисления приводят к равенству
\begin{equation}\label{eq:1.5A}
    \Delta(\lambda)=J_{12}+J_{34}+J_{13}\psi(\pi,\lambda)+
    J_{14}\psi^{[1]}(\pi,\lambda)+J_{32}\varphi(\pi,\lambda)+
    J_{42}\varphi^{[1]}(\pi,\lambda).
\end{equation}
Из этого представления и асимптотических равенств, сформулированных в
лемме~\ref{lem:1.3}, непосредственно следует, что \(\Delta(\lambda)\not\equiv
0\) внутри любой параболы \(\Pi_a\), если выполнено одно из условий
1)\,--\,3).

Так как \(\Delta(\lambda)\) есть голоморфная функция от \(\lambda\) во всей
комплексной плоскости \(\mathbb C\), то нули \(\Delta(\lambda)\) образуют
последовательность, не имеющую конечных предельных точек. Если
\(\Delta(\lambda_0)\neq 0\), то в силу теоремы~\ref{tm:0} оператор
\((L-\lambda_0)^{-1}\) отображает единичный шар пространства \(L_2\) в
ограниченное множество пространства \(W_1^1\). Следовательно,
\((L-\lambda_0)^{-1}\) компактен. Но тогда \(L\) имеет дискретный спектр,
который совпадает с нулями определителя \(\Delta(\lambda)\). Теорема доказана.
\end{proof}

\begin{Note}
При изучении операторов с регулярными потенциалами,
краевые условия, для которых выполняется одно из условий 1)\,--\,3),
называют регулярными по Биркгофу. Если
\(u(x)\) --- гладкая функция, то замена в краевых условиях переменных
\(y'(0)\), \(y'(\pi)\) на квазипроизводные \(y^{[1]}(0)\), \(y^{[1]}(\pi)\)
сохраняет свойство регулярности краевых условий.
Можно отметить также, что утверждение
предыдущей теоремы сохраняется для так называемых невырожденных
краевых условий --- эти условия отличаются от регулярных тем, что пункт
3) заменяется следующим\hfill\break
3') \(J_{42}=0\), \(J_{14}+J_{32}=0\), \(J_{13}\neq 0\).\hfill\break
(см. \cite{Mar}). Для доказательства такого уточнения теоремы~\ref{tm:1.4}
недостаточно асимптотики~\eqref{zz} --- необходимо раскрыть символы
$o(1)$ так, как это сделано в лемме~\ref{lem:2.4}. Однако в дальнейшем в
данной работе нас будут интересовать регулярные краевые условия, так как
для рассматриваемой задачи при \(u(x)\in L_2\) операторы с регулярными
краевыми условиями сохраняют классические асимптотики для собственных
значений и собственных функций, кроме того, система их собственных и
присоединенных функций образует базис Рисса.
\end{Note}

В случае вещественности функции \(u(x)\) минимальный оператор \(L_m\)
симметричен с индексами дефекта \((2,2)\). Несложно описать все
самосопряженные расширения \(L_m\).
\begin{Theorem}\label{tm:1.5}
Если функция \(u(x)\) вещественна, то произвольное самосопряженное расширение
\(L\) симметрического оператора \(L_m\) является сужением оператора \(L_M\) на
область
\[
    \mathfrak D(L)=\left\{y\;\vline\;y\in\mathfrak D(L_M),\;
    U_1(y)=U_2(y)=0\right\},
\]
где линейные формы \(U_1\), \(U_2\) имеют представление~\eqref{eq:T.1.4}, для
коэффициентов которых выполнены равенства
\begin{align}\label{eq:T.1.5}
    AJA^*-BJB^*&=0,&A&=\begin{pmatrix}a_{11}&a_{12}\\a_{21}&a_{22}
    \end{pmatrix},\;B=\begin{pmatrix}b_{11}&b_{12}\\b_{21}&b_{22}
    \end{pmatrix},\;J=\begin{pmatrix}0&1\\-1&0\end{pmatrix}.
\end{align}
Наоборот, любые краевые условия вида~\eqref{eq:T.1.4}, \eqref{eq:T.1.5}
определяют самосопряженный оператор \(L\).
\end{Theorem}
\textit{Доказательство} этой теоремы получается почти дословным повторением
рассуждений из работы Крейна~\cite[\S\,3]{Kr} и здесь не приводится.
\hfill\(\square\)
\begin{Note}\label{not:1.2}
Полезно отметить, что краевые условия, определяющие самосопряженные расширения,
можно записать также в форме
\begin{equation*}
    (U-1)\begin{pmatrix}y^{[1]}(0)\\-y^{[1]}(\pi)\end{pmatrix}+
    i(U+1)\begin{pmatrix}y(0)\\y(\pi)\end{pmatrix}=0,
\end{equation*}
где \(U\) --- произвольная унитарная матрица второго порядка. Доказательство
эквивалентности такой записи предыдущей проводится так же, как в монографии
Рофе--Бекетова и Холькина~\cite{RH}. Полезно также отметить, что краевые
условия, определяющие самосопряженный оператор, обязательно удовлетворяют
одному из условий~1\,--\,3 теоремы~\eqref{tm:1.4}, т.~е. являются регулярными
по Биркгофу. Доказательство этого утверждения получается элементарной
проверкой.
\end{Note}

\textbf{1.2. Аппроксимация гладкими потенциалами.} Другой подход к определению
оператора \(L\) с потенциалом-распределением основан на идее аппроксимации
операторами с гладкими потенциалами.

Пусть \(q(x)\in W_2^{-1}[0,\pi]\), \(u(x)=\int q(t)\,dt\). Пусть
\(q_{\varepsilon}(x)\) --- семейство гладких функций, таких, что
\(\|q_{\varepsilon}(x)-q(x)\|_{W_2^{-1}}\to 0\) при \(\varepsilon\to 0\). Это
условие эквивалентно тому, что \(u_{\varepsilon}(x)=\int
q_{\varepsilon}(t)\,dt\to u(x)\) при \(\varepsilon\to 0\) в пространстве
\(L_2\).

Обозначим через \(L_{\varepsilon}\) оператор, порожденный дифференциальным
выражением \(l_{\varepsilon}(y)=-y''+q_{\varepsilon}(x)y\) и регулярными
краевыми условиями~\eqref{eq:T.1.4}, в которых переменные \(y^{[1]}(0)\),
\(y^{[1]}(\pi)\) определяются равенством
\begin{equation}\label{eq:Q}
    y^{[1]}(x)=y'(x)-u_{\varepsilon}(x)y(x).
\end{equation}
В случае гладких функций \(u_{\varepsilon}(x)\) подстановка в краевые условия
переменных \(y^{[1]}(0)\), \(y^{[1]}(\pi)\) вместо обычных производных
сохраняет регулярность краевых условий. Поэтому (см.~\cite[Гл.~1]{Naj})
операторы \(L_{\varepsilon}\) корректно определены и имеют дискретный спектр.
Оказывается справедливым следующий результат.
\begin{Theorem}\label{tm:1.6}
Существуют значения \(\lambda\in\mathbb C\) такие, что при всех достаточно
малых \(\varepsilon>0\) значение \(\lambda\) принадлежит резольвентным
множествам операторов \(L_{\varepsilon}\), а последовательность
\((L_{\varepsilon}-\lambda)^{-1}\) фундаментальна при \(\varepsilon\to 0\) в
равномерной операторной топологии, т.~е.
\begin{align*}
    \|(L_{\varepsilon}-\lambda)^{-1}-(L_{\delta}-\lambda)^{-1}\|&\to 0
    &\text{при }\varepsilon,\delta&\to 0.
\end{align*}
Оператор \(T\), являющийся пределом последовательности
\((L_{\varepsilon}-\lambda)^{-1}\) не имеет ядра, а потому на области значений
\(T\) определен оператор \(T^{-1}\). При этом оператор \(T^{-1}+\lambda\)
совпадает с оператором \(L\), определенным в теореме~\ref{tm:1.5}.
\end{Theorem}
\begin{proof}
Здесь, как и в первом подходе, вновь основную роль играет теорема~\ref{tm:0}.
Пусть
\begin{align*}
    u_{\varepsilon}(x)&\stackrel{L_2}{\to}u_0(x)&\text{при }
    \varepsilon&\to 0,
\end{align*}
а \(L_{\varepsilon}\) --- определенные выше операторы с регулярными краевыми
условиями. Обозначим через \(\varphi_{\varepsilon}(x)\),
\(\psi_{\varepsilon}(x)\) пару решений уравнения
\begin{equation}\label{eq:T.1.7}
    -y''+u_{\varepsilon}'(x)y=\lambda y,
\end{equation}
удовлетворяющих условиям
\[
    \varphi_{\varepsilon}(0)=1,\quad\varphi_{\varepsilon}^{[1]}(0)=0;\quad
    \psi_{\varepsilon}(0)=0,\quad\psi_{\varepsilon}^{[1]}(0)=1.
\]
При \(\varepsilon=0\) левую часть~\eqref{eq:T.1.7} понимаем так же, как
в~\eqref{eq:1.1}.

Уравнение~\eqref{eq:T.1.7} перепишем в виде системы~\eqref{eq:S1}, при этом
матрица коэффициентов \(\mathbf A\) соответствующей системы будет содержать
функции \(u_{\varepsilon}(x)\) и \(u_{\varepsilon}^2(x)\). Согласно
теореме~\ref{tm:0}, имеем (мы полагаем \(z^{[0]}=z\))
\begin{multline}\label{eq:bas}
    \|\varphi_{\varepsilon}^{[j]}(x)-\varphi_{\delta}^{[j]}(x)\|_{1,1}+
    \|\psi_{\varepsilon}^{[j]}(x)-\psi_{\delta}^{[j]}(x)\|_{1,1}\le
    C\left(\|u_{\varepsilon}(x)-u_{\delta}(x)\|_{L_1}+
    \|u_{\varepsilon}^2(x)-u_{\delta}^2(x)\|_{L_1}\right)\le\\
    \le C_1\left(\|u_{\varepsilon}(x)-u_{\delta}(x)\|_{L_2}\right),
    \quad j=0,1,
\end{multline}
где \(\|\cdot\|_{1,1}\) --- норма в \(W_1^1\).

Так как соответствующая система имеет решение при \(\delta=0\), то
неравенство~\eqref{eq:bas} остается справедливым и при \(\delta=0\).

Учитывая, что вронскиан пары \(\varphi_{\varepsilon}\), \(\psi_{\varepsilon}\)
тождественно равен \(1\), непосредственной проверкой убеждаемся, что функция
\[
    z_{\varepsilon}(x)=\int\limits_0^x\left[\varphi_{\varepsilon}(x)
    \psi_{\varepsilon}(\xi)-\psi_{\varepsilon}(x)\varphi_{\varepsilon}
    (\xi)\right]f(\xi)\,d\xi
\]
удовлетворяет резольвентному равенству
\begin{equation}\label{eq:1.6A}
    -y''+u_{\varepsilon}'(x)y-\lambda y=f(x)\in L_2[0,\pi].
\end{equation}
Общее решение~\eqref{eq:1.6A} имеет вид
\begin{equation}\label{eq:1.6B}
    y_{\varepsilon}(x)=z_{\varepsilon}(x)+a_1\varphi_{\varepsilon}(x)+
    a_2\psi_{\varepsilon}(x),
\end{equation}
где \(a_1\), \(a_2\) --- постоянные. Подставив это решение в краевые условия,
получим
\begin{gather*}
    \begin{aligned}
        a_1&=\Delta^{-1}\begin{vmatrix}U_1(z_{\varepsilon})&
        U_1(\psi_{\varepsilon})\\ U_2(z_{\varepsilon})&
        U_2(\psi_{\varepsilon})\end{vmatrix},\quad &
            a_2&=\Delta^{-1}\begin{vmatrix}U_1(\varphi_{\varepsilon})&
        U_1(z_{\varepsilon})\\ U_2(\varphi_{\varepsilon})&
        U_2(z_{\varepsilon})\end{vmatrix},\quad
    \end{aligned}
    \Delta=\begin{vmatrix}U_1(\varphi_{\varepsilon})&
    U_1(\psi_{\varepsilon})\\ U_2(\varphi_{\varepsilon})&
    U_2(\psi_{\varepsilon})\end{vmatrix}.
\end{gather*}
Поскольку краевые условия регулярны, то \(\Delta=\Delta_{\varepsilon}(\lambda)
\not\equiv 0\) при всех \(\varepsilon\ge 0\) (при \(\varepsilon=0\) это
вытекает из леммы~\ref{lem:1.3} и представления~\eqref{eq:1.5A}). Из
оценки~\eqref{eq:bas} имеем
\begin{align*}
    |U_s(\varphi_{\varepsilon})-U_s(\varphi_{\delta})|+
    |U_s(\psi_{\varepsilon})-U_s(\psi_{\delta})|&\le C\|u_{\varepsilon}-
    u_{\delta}\|_{L_2},&s&=1,2,
\end{align*}
а потому \(|\Delta_{\varepsilon}-\Delta_{\delta}|\le C\|u_{\varepsilon}-
u_{\delta}\|_{L_2}\). В силу теоремы~\ref{tm:0} для функций
\(z_{\varepsilon}(x)\) также имеем оценку
\begin{align*}
    \|z_{\varepsilon}^{[j]}(x)-z_{\delta}^{[j]}(x)\|_{1,1}\le
    C\|u_{\varepsilon}^2-u_{\delta}^2\|_{L_2}\|f\|_{L_1}&\le
    C\|u_{\varepsilon}-u_{\delta}\|_{L_2}\|f\|_{L_1},& j&=0,1.
\end{align*}
Выберем теперь число \(\lambda\) таким, чтобы \(\Delta_0(\lambda)\neq 0\).
Тогда \(|\Delta_{\varepsilon}(\lambda)|>c=c(\lambda)\) при всех достаточно
малых \(\varepsilon>0\). Следовательно, постоянные \(a_1=a_1(\varepsilon)\),
\(a_2=a_2(\varepsilon)\) в~\eqref{eq:1.6B} таковы, что
\[
    |a_1(\varepsilon)-a_1(\delta)|+|a_2(\varepsilon)-a_2(\delta)|\le
    C\|u_{\varepsilon}-u_{\delta}\|_{L_2}\|f\|_{L_2},
\]
где \(C\) зависит только от выбранного числа \(\lambda\). Полученные оценки
показывают, что для решений
\[
    y_{\varepsilon}=(L_{\varepsilon}-\lambda)^{-1}f
\]
справедливы неравенства
\[
    \|y_{\varepsilon}-y_{\delta}\|_{L_2}\le C\|y_{\varepsilon}-
    y_{\delta}\|_{1,1}\le C\|u_{\varepsilon}-u_{\delta}\|_{L_2}
    \|f\|_{L_2}.
\]
Тем самым доказано соотношение~\eqref{eq:T.1.7}, или равномерная резольвентная
сходимость операторов \(L_{\varepsilon}\) при \(\varepsilon\to 0\). Заметим,
что если \(T\) есть предел операторов \((L_{\varepsilon}-\lambda)^{-1}\), то
\(y=Tf\) есть решение уравнения~\eqref{eq:1.6A} при \(\varepsilon=0\) и
удовлетворяет соответствующим краевым условиям. Следовательно, \(y\neq 0\),
если \(f\neq 0\), т.~е. \(\ker T=\{0\}\). При этом очевидно,
\(T=(L-\lambda)^{-1}\), где оператор \(L\) был построен в
теореме~\ref{tm:1.6}. Этим завершается доказательство теоремы.
\end{proof}

\textbf{1.3. Метод квадратичных форм.} Сначала мы предположим, что
\(u(x)\) --- вещественная функция. Выпишем квадратичную форму,
отвечающую дифференциальному выражению~\eqref{eq:1.1}. Имеем
\begin{equation}\label{eq:QF}
    (l(y),y)=-((y^{[1]})',y)-(u(x)y^{[1]},y)-(u^2(x)y,y)=
    (y^{[1]},y^{[1]})-(u^2(x)y,y)+(y^{\vee},y^{\wedge}),
\end{equation}
где используются обозначения
\[
    y^{\wedge}=\begin{pmatrix}y(0)\\y(\pi)\end{pmatrix},\qquad
    y^{\vee}=\begin{pmatrix}y^{[1]}(\pi)\\-y^{[1]}(\pi)
    \end{pmatrix}
\]
(см.~замечание~\ref{not:1.2}). Пусть \(A\) --- произвольная самосопряженная
матрица размера \(2\times 2\). Положим
\begin{equation}\label{eq:1.17app}
    W_{2,U}^1=\left\{y\in W_2^1[0,\pi]\;\;\vline\;\;Uy^{\wedge}=0
    \right\},
\end{equation}
где \(U\) --- произвольная матрица размера \(2\times 2\). Очевидно,
\(W_{2,U}^1\) есть подпространство в \(W_2^1\) коразмерности \(\le 2\),
в зависимости от ранга матрицы \(U\). При \(U=0\) имеем \(W_{2,U}^1=
W_2^1\). На пространстве \(W_{2,U}^1\) определим квадратичную форму
\begin{equation}\label{eq:1.7}
    \mathcal L(y,y)=(y^{[1]},y^{[1]})-(u^2(x)y,y)+(Ay^{\wedge},
    y^{\wedge}).
\end{equation}
\begin{Lemm}\label{lem:1.7}
Квадратичная форма~\eqref{eq:1.7} определена при \(y\in W_{2,U}^1\) и
замкнута.
\end{Lemm}
\begin{proof}
Если \(y\in W_{2,U}^1\), а \(u(x)\in L_2\), то \(u(x)y(x)\) и \(y'(x)\)
принадлежат \(L_2\), а потому форма~\eqref{eq:1.7} корректно определена.
Далее,
\[
    (u^2(x)y,y)\le\|u\|_{L_2}^2\|y\|_C^2\le\varepsilon\|y\|_{1,2}^2
    +M\|y\|_{L_2}^2,
\]
где число \(\varepsilon\) можно выбрать произвольно малым, а постоянная
\(M\) зависит от \(\varepsilon\) и \(\|u\|_{L_2}\). Последняя оценка
вытекает из компактности вложения \(W_2^1\subset C[0,\pi]\). Имеем также
\[
    |(y',u(x)y)|\le\|y'\|\|u(x)y\|\le\|y'\|\|u(x)\|\|y\|_C\le
    \varepsilon\|u(x)\|\|y'\|\|y\|_{1,2}+M\|y'\|\|u(x)\|\|y\|\le
\]
\[
       \le(\varepsilon\|u(x)\|+M\eps)\|y\|_{1,2}+
    \varepsilon^{-1}M\|u(x)\|^2\|y\|^2\le\eps_1\|y\|^2_{1,2}+M_1\|y\|^2.
\]
Наконец,
\[
    |(Ay^{\wedge},y^{\wedge})|\le C\|y\|_C^2\le\varepsilon\|y'\|^2+
    M\|y\|^2.
\]
Полученные оценки позволяют получить представление
\[
    \mathcal L(y,y)=(y',y')+(y,y)+\mathcal
    Q(y,y)=\|y\|_{1,2}^2+\mathcal Q(y,y),
\]
где
\[
    |\mathcal Q(y,y)|\le\varepsilon\|y\|_{1,2}^2+M\|y\|_{L_2}^2.
\]
Выбрав \(\varepsilon<1\), из известной теоремы
(см.~\cite[Гл.~6.1]{Kato}) получим замыкаемость формы \(\mathcal
L(y,y)\) в пространстве \(W_{2,U}^1\). Лемма доказана.
\end{proof}

Замкнутая квадратичная форма~\eqref{eq:1.7} (она зависит от выбора
матриц \(A\) и \(U\)), согласно первой теореме о представлении (см.,
например, \cite[Гл.~\mbox{VI}.2]{Kato}), определяет самосопряженный
полуограниченный оператор \(L\), причем область определения
квадратного корня \((L+\alpha)^{1/2}\) (здесь \(\alpha>0\) --- достаточно
большое число) совпадает с \(W_{2,U}^1\). Конечно, таким образом могут
быть получены все описанные ранее самосопряженные расширения минимального
оператора \(L_m\). Например, если самосопряженное расширение \(L\)
описывается краевыми условиями \((U-1)y^{\vee}+i(U+1)y^{\wedge}=0\),
где матрица \(U-1\) обратима (это соответствует фигурирующему в
теореме~\ref{tm:1.4} условию \(J_{42}\neq 0\)), то этому расширению
соответствует квадратичная форма~\eqref{eq:1.7}, определенная на всем
пространстве \(W_2^1\), причем матрица \(A\) в~\eqref{eq:1.7} находится
из условия \(A=-i(U-1)^{-1}(U+1)\), т.~е. является преобразованием Кэли
от \(U\).

Подход с помощью метода квадратичных форм позволяет нам получить дополнительную
информацию об области оператора \(L\). А именно, справедливо следующее
утверждение.
\begin{Theorem}\label{tm:1.8}
Пусть \(L\) --- самосопряженное расширение минимального оператора
\(L_m\), а \(W_{2,U}^1\) --- подпространство в \(W_2^1\), состоящее
из функций, которые удовлетворяют краевым условиям нулевого порядка
(имеются в виду краевые условия вида~\eqref{eq:T.1.4}, порождающие расширение
\(L\)). Тогда
\begin{equation}\label{eq:1.8}
    \mathfrak D(L)=\left\{y\in W_{2,U}^1\;\vline\;
    l(y)\in L_2\right\},
\end{equation}
где равенство \(-y''+q(x)y=f(x)\in L_2\) понимается в смысле теории
распределений.
\end{Theorem}
\begin{proof}
Заметим, что все самосопряженные расширения \(L_m\) полуограничены,
так как \(L_m\) полуограничен и его индексы дефекта конечны. Не
ограничивая общности, считаем, что \(L>0\). Тогда \(\mathfrak D(L)
\subset\mathfrak D(L^{1/2})\), а потому \(\mathfrak D(L)\subset
W_{2,U}^1\). Условие \(l(y)\in L_2\) следует из определения оператора
\(L\). Теорема доказана.
\end{proof}
\begin{Note}
Определение~\eqref{eq:1.8} области \(\mathfrak D(L)\) вполне согласуется
с определением~\eqref{eq:Lm}. Действительно, если
\[
    -y''+q(x)y=-(y^{[1]})'-u(x)y^{[1]}-u^2(x)y=f(x)\in L_2,
\]
то \((y^{[1]})'\in L_1[0,\pi]\), а потому \(y^{[1]}\in W_1^1[0,\pi]\).
В частности, при \(y\in \mathfrak D(L)\) интегрирование по частям
в~\eqref{eq:QF} является корректным.
\end{Note}

Предложенный метод квадратичных форм может быть также использован и для
определения операторов с комплексным потенциалом-распределением \(q(x)\).
Однако этот метод не дает определения операторов с произвольными регулярными
(или более общими) краевыми условиями, а только краевыми условиями, которые
являются подчиненными возмущениями самосопряженных. В случае комплексной
функции \(u(x)\) равенство~\eqref{eq:QF} следует записать в виде
\[
    (l(y),y)=(y',y')-(u(x)y,y')-(u(x)y',y)+(y^{\vee},y^{\wedge}).
\]
Это равенство позволяет ассоциировать с дифференциальным выражением \(l(y)\)
квадратичную форму
\[
    \mathcal L(y,y)=(y',y')-(u(x)y,y')-(u(x)y',y)+
    (Ay^{\wedge},y^{\wedge}),
\]
где \(A\) --- произвольная комплексная матрица размера \(2\times 2\). Здесь
предполагается, что форма \(\mathcal L(y,y)\) определена на пространстве
\(W_{2,U}^1\). Эта квадратичная форма не является вещественной, но она
секториальна и является \(\varepsilon\)-подчиненной форме \((y',y')\) при
любом \(\varepsilon>0\). Следовательно (см.~\cite[Гл.~6.2]{Kato}),
существует максимальный секториальный оператор \(L\), порождающий эту форму.
Произвольные расширения \(L\), полученные в предыдущем пункте, не могут
получаться на этом пути. Более точно, этим методом получаются те операторы
\(L\) из теоремы~\ref{tm:1.4}, для которых равенства \(U_j(y)=0\), \(j=1,2\),
влекут возможность представления
\[
    (y^{\vee},y^{\wedge})=(Ay^{\wedge},y^{\wedge}).
\]

\textbf{1.4. Метод мультипликаторов.} Пусть \(q(x)\in\mathcal D'\), а
\(y\in\mathcal D\). Тогда
\[
    \mathcal L(y,y)=(-y''+q(x)y,y)=(y',y')+(q(x)y,y).
\]
Если справедлива оценка
\begin{equation}\label{eq:Es}
    |(q(x)y,y)|\le\varepsilon(y',y')+M(y,y),\qquad M=M(\varepsilon),
\end{equation}
то квадратичная форма секториальна и замыкаема, причем область ее замыкания
совпадает с пространством \(\Wo_2^1\). В этом случае с формой \(\mathcal L\)
можно ассоциировать оператор. Естествен вопрос: для каких функций
\(q(x)\in\mathcal D'\) оценка~\eqref{eq:Es} справедлива? Для ответа на него
полезно ввести следующее понятие. Функцию \(q(x)\in\mathcal D'\) назовем
мультипликатором из пространства \(\Wo_2^1\) в дуальное пространство
\(W_2^{-1}\), если
\begin{equation}\label{eq:Es1}
    |(q(x)y,y)|\le
    C\|y\|_{1,2}^2,\qquad\forall y\in\mathcal D,
\end{equation}
где постоянная \(C\) не зависит от \(y\), а \(\|\cdot\|_{1,2}\) --- норма в
\(\Wo_2^1\). Очевидно, мультипликаторы образуют линейное пространство с нормой
\(\|q\|=\inf C\), где инфинум берется среди постоянных \(C\)
в~\eqref{eq:Es1}.
Это пространство обозначим через \(M[1]\). В статье~\cite{NS} показано, что
\(W_2^{-1}\subset M[1]\), а в статье~\cite{BS} доказано равенство
\(M[1]=W_2^{-1}\) и эквивалентность норм в этих пространствах. Хотя в
работах~\cite{NS} и~\cite{BS} рассматриваются операторы на всей прямой
\(\mathbb R\) (и их многомерные обобщения в \(\mathbb R^n\)),
доказательства не
меняются при переходе на конечный интервал. В нашем случае справедливость
включения \(W_2^{-1}\subset M[1]\) очевидна в силу оценки
\[
    |(q(x)y,y)|\le |(q(x),y\overline{y})|\le\|q\|_{-1,2}
    \|y\overline{y}\|_{1,2}\le\|q\|_{-1,2}\|y\|_{1,2}^2,\quad
    y\in\mathcal D.
\]
Далее, заметим, что гладкие функции \(\varphi\in\mathcal D\) плотны в
пространстве \(W_2^{-1}=M[1]\), следовательно, для любого мультипликатора
\(q\in W_2^{-1}\) выполнена оценка~\eqref{eq:Es}. Тем самым, для любой \(q\in
W_2^{-1}\) определен оператор \(L\), ассоциированный с формой \(\mathcal
L(y,y)\). Очевидно, этот оператор \(L\) совпадает с прежними определениями
оператора \(L\), отвечающего краевым условиям Дирихле
\[
    y(0)=y(\pi)=0.
\]
Однако этот метод можно распространить для определения операторов с более
общими краевыми условиями, хотя краевые условия либо не будут фигурировать
вовсе, либо будет одно условие нулевого порядка.

Рассмотрим подпространство \(W_{2,U}^1\subseteq W_2^1\) коразмерности \(1\)
или \(0\) (см.~\eqref{eq:1.17app}. На этом подпространстве определим
квадратичную форму
\begin{align}\label{eq:NQ} \mathcal
    L(y,y)&=(y',y')+(q(x),y\overline{y}),&y&\in W_{2,U}^1.
\end{align}
Если \(W_{2,U}^1=W_2^1\), то условие \(y\in W_2^1\) влечет \(y\overline{y}\in
W_2^1\). Это свойство сохраняется для \(W_{2,U}^1\), если краевое условие,
порождающее это пространство, имеет вид \(y(0)=0\) или \(y(\pi)=0\), либо
\(y(0)-\alpha y(\pi)=0\) и \(\alpha=\pm 1\). При \(\alpha\neq\pm 1\) функция
\(y\overline{y}\in W_{2,V}^1\), где индекс \(V\) означает новое краевое условие
\(V(\varphi)=\varphi(0)-|\alpha|^2\varphi(\pi)=0\). Из определения
формы~\eqref{eq:NQ} и включения \(y\overline{y}\in W_{2,V}^1\) следует, что
форма \(\mathcal L\) корректно определена на \(W_{2,U}^1\), если
\(q\in(W_{2,V}^1)'\) --- дуальному пространству к \(W_{2,V}^1\) по отношению к
скалярному произведению в \(L_2\).

Пространство \(\Wo_2^1\) имеет коразмерность \(1\) или \(2\) в \(W_{2,V}^1\),
поэтому \(W_2^{-1}\) (дуальное к \(\Wo_2^1\)) имеет коразмерность \(1\) или
\(2\) в \((W_{2,V}^1)'\). Например, в случае \(W_{2,U}^1=W_2^1\) имеем
\((W_2^1)'=W_2^{-1}\oplus\mathfrak N\), где \(\mathfrak N\) --- двумерное
подпространство, содержащее функционалы \(F_0(y)=y(0)\) и \(F_1(y)=y(\pi)\).
Дуальным к пространству \(W_{2,V}^1\) с краевым условием
\(y(0)-\alpha^2y(\pi)=0\) будет пространство \(W_2^{-1}\oplus\mathfrak N\),
где одномерное пространство \(\mathfrak N\) содержит функционал
\(F(y)=\alpha^2y(0)+y(\pi)\). Из плотности гладких функций в \(W_2^{-1}\)
можно вывести оценку
\[
    |(q,y\overline{y})|\le\varepsilon(y',y')+M(y,y),\qquad y\in W_{2,U}^1,
\]
а потому квадратичная форма~\eqref{eq:NQ} определяет некоторый секториальный
оператор \(L\).

Недостаток такого определения оператора \(L\) состоит в том, что мы не
указываем явную формулу для квадратичной формы~\eqref{eq:NQ} через регулярную
функцию \(u(x)=\int q(\xi)\,d\xi\). Однако этот недостаток можно исправить и
написать явную формулу, используя представление
\[
    y\overline{y}=(y-\psi)\overline{y}+\psi(\overline{y}-\overline{\psi})+
    \psi\overline{\psi},
\]
где \(\psi=y(0)+\pi^{-1}(y(\pi)-y(0))x\). Функция \(y-\psi\) аннулируется
на концах отрезка, а потому справедливо равенство
\[
    (q(x),y\overline{y})=-(u(x),[y\overline{y}-\psi\overline{\psi}]')+
    (q(x),\psi\overline{\psi})=-(u(x),(y\overline{y})')+(u(x),(\psi
    \overline{\psi})')+(q,\psi\overline{\psi}),
\]
при этом выбор значения \((q,\psi\overline{\psi})\) находится в нашей власти.
Несложно видеть, что описанный метод позволяет определить такой же класс
операторов, как метод квадратичных форм.

\textbf{1.5 Обсуждение условия на потенциал \(q(x)\).}
Во всех предложенных методах определения оператора \(L\) фигурировало условие
\(q(x)\in W_2^{-1}\).Мы не видим возможности получить похожие результаты,
например, для функций $u(x)\in L_p$, $p<2$. Поставим
следующий вопрос: можно ли определить оператор $L$ для $u(x)\in L_p$,
$p<2$, с помощью предельного перехода? Точнее, если последовательность
гладких функций $u_\varepsilon(x)$ такова, что
$\|u_\varepsilon-u\|_{L_p}\to 0$ при $\varepsilon\to 0$, то будет ли
соответствующая
последовательность операторов $L_\varepsilon$ иметь сильный или равномерный
резольвентный предел, не зависящий от выбора последовательности
$u_\varepsilon$? Покажем, что ответ на этот вопрос отрицателен. Далее нам
будет удобно вместо отрезка \([0,\pi]\) рассматривать отрезок \([-1,1]\).
\begin{Proposition}
Рассмотрим последовательность потенциалов
\begin{equation}
q_\eps(x)=u'_\varepsilon(x)=\left\{\begin{array}{lll}
0,\quad&\mbox{если}\ x\in[-1,-\eps]\cup[\eps,1];\\
\eps^{-3/2},\quad&\mbox{если}\ x\in(-\eps,0);\\
-\eps^{-3/2},\quad&\mbox{если}\ x\in(0,\eps).
\end{array}
\right.
\label{aa}
\end{equation}
на отрезке $[-1,1]$. Последовательность операторов Штурма--Лиувилля
\(L_{\varepsilon}\) с потенциалами \(q_{\varepsilon}(x)\) и граничными
условиями Дирихле $y(-1)=y(1)=0$ имеет предел в смысле сильной резольвентной
сходимости. Предельный оператор \(L_0\) задается дифференциальным выражением
\[
    -y''-\dfrac23\delta(x)y
\]
и краевыми условиями Дирихле.
\end{Proposition}
\begin{proof}
Отметим, что потенциалы вида~\eqref{aa} появились в связи с задачей
определения оператора Штурма--Лиувилля с потенциалом $\delta'$
(см. \cite{S})). Нам нужно доказать сильную резольвентную сходимость
$$
L_\eps:=-\frac{d^2}{dx^2}+q_\eps(x)\lar_{\eps\to0}L_0:=-\frac{d^2}{dx^2}-
\frac23\delta(x).
$$
В силу вещественности потенциала, операторы $L_\eps$, $\eps\ge0$
самосопряженны. Таким образом, достаточно доказать сильную сходимость
резольвент $R_\eps(\lambda)\lar_{\eps\to0}R_0(\lambda)$ для
$\lambda\in\C\backslash\R$ (на самом деле, достаточно доказать это
соотношение для $\lambda=i$ и $\lambda=-i$ (см. \cite{Kato})).

Пусть $\varphi_\varepsilon(x)$ и
$\psi_\varepsilon(x)$-- пара решений однородного уравнения
\begin{displaymath}
-y''+u'_\varepsilon(x)y=\lambda y,
\end{displaymath}
таких, что
$$
\varphi_\varepsilon(-1)=1,\quad\varphi_\varepsilon^{[1]}(-1)=0;\quad
\psi_\varepsilon(-1)=0,\quad\psi_\varepsilon^{[1]}(-1)=1.
$$
Здесь квазипроизводные понимаются в смысле равенства \eqref{eq:Q}. При
\(\varepsilon=0\) полагаем \(u_0'(x)=-2/3\delta(x)\).

Функция
\begin{equation}
z_\varepsilon(x)=\int\limits_0^x(\varphi_\varepsilon(x)\psi_\varepsilon(\xi)
-\psi_\varepsilon(x)\varphi_\varepsilon(\xi))f(\xi)d\xi
\label{a3}
\end{equation}
удовлетворяет резольвентному равенству
\begin{equation}
-y''+u'_\varepsilon(x)y-\lambda y=f(x)\in L_2[0,\pi]
\label{a1}
\end{equation}
(см. доказательство теоремы~\ref{tm:1.6}).
Общее решение уравнения (\ref{a1}) имеет вид
$$
y_\varepsilon(x)=z_\varepsilon(x)+a_1\varphi_\varepsilon(x)+
a_2\varphi_\varepsilon(x).
$$
Подставив его в краевые условия, получим
$$
y_\eps(-1)=a_1=0,\quad y_\eps(1)=z_\eps(1)+a_2\psi_\eps(1)=0,
$$
т.е. $a_2=-z_\eps(1)/\psi_\eps(1)$. Покажем, что
\begin{equation}
\|\varphi_\eps(x)-\varphi_0(x)\|_{L_2[-1,1]}+
\|\psi_\eps(x)-\psi_0(x)\|_{L_2[-1,1]}\lar_{\eps\to0}0.
\label{a2}
\end{equation}
Тогда из формулы \eqref{a3} следует
$\|z_\eps(x)-z_0(x)\|_{L_2}\lar_{\eps\to0}0$, а поэтому
$|a_2(\eps)-a_2(0)|\lar_{\eps\to0}0$ и
$\|y_\eps(x)-y_0(x)\|_{L_2}\lar_{\eps\to0}0$. Таким образом,
доказательство сводится к проверке соотношения (\ref{a2}). Учитывая
простой вид потенциалов $q_\eps(x)$, это соотношение можно доказать
прямыми вычислениями, но для того, чтобы избежать утомительных подсчетов,
мы докажем его другим способом.

Пусть
$$
u_\eps(x)=\int\limits_{-1}^x q_\eps(t)dt=\left\{\begin{array}{lll}
0,\quad&\mbox{при}\ x\in[-1,-\eps]\cup[\eps,1];\\
\eps^{-3/2}(x+\eps),\quad&\mbox{при}\ x\in[-\eps,0];\\
\eps^{-3/2}(\eps-x),\quad&\mbox{при}\ x\in[0,\eps].
\end{array}\right.
$$
Обозначим также
$$
v_\eps(x):=\int\limits_{-1}^x u^2_\eps(t)dt=\left\{\begin{array}{llll}
0,\quad&\mbox{при}\ x\in[-1,-\eps];\\
1/3\eps^{-3}(x+\eps)^3,\quad&\mbox{при}\ x\in[-\eps,0];\\
2/3+1/3\eps^{-3}(x-\eps)^3,\quad&\mbox{при}\
x\in[0,\eps];\\
2/3,\quad&\mbox{при}\ x\in[\eps,1].
\end{array}\right.
$$
При $\eps>0$ от уравнения \eqref{a1} можно перейти к системе двух
уравнений первого порядка с помощью замены
\begin{gather*}
y_{1,\eps}(x)=y_\eps(x),\\
y_{2,\eps}(x)=y'_\eps(x)+(v_\eps(x)-u_\eps(x))y_\eps(x)
\end{gather*}
(такая замена является на самом деле повторной регуляризацией
дифференциального выражения (\ref{de1}) (см. \S\,4 настоящей работы)).
Система будет иметь вид
\begin{equation}
\left(\begin{array}{c} y_1 \\ y_2 \end{array}\right)'=
\left(\begin{array}{cc} u_\eps-v_\eps&1\\
-\lambda+2v_\eps u_\eps-v^2_\eps\hskip40pt &v_\eps-u_\eps
\end{array}\right)
\left(\begin{array}{c} y_1\\y_2
\end{array}\right).
\label{a5}
\end{equation}

При $\eps\to0$ видим, что в пространстве $L_1[-1,1]$
$$
u_\eps(x)\to0,\quad v_\eps(x)\to\frac23\chi(x),\quad
v^2_\eps(x)\to\frac49\chi(x),\quad
u_\eps(x)v_\eps(x)\to0,
$$
где $\chi(x)$ --- функция Хевисайда. Тогда матрицы системы (\ref{a5})
сходятся к матрице
\begin{equation}
\left(\begin{array}{ccc} -\frac23\chi(x)&1\\\phantom{.}\\
-\lambda-(\frac23\chi(x))^2 &\frac23\chi(x)
\end{array}\right)
\label{a6}
\end{equation}
В силу теоремы~\ref{tm:0}, решение системы (\ref{a5}) сходится в пространстве
$W_1^1[-1,1]$ к решению системы с матрицей (\ref{a6}). Остается заметить,
что первая компонента вектора решения такой системы есть решение уравнения
$$
-y''-\frac23\delta(x)y=\lambda y
$$
с соответствующими начальными условиями.
\end{proof}

Итак, мы показали, что сильным резольвентным пределом последовательности
операторов Штурма--Лиувилля с потенциалами (\ref{aa}) и граничными
условиями Дирихле будет оператор с потенциалом $-\frac23\delta(x)$.
С другой стороны, очевидно, что $u_\varepsilon(x)\to 0$ в
$L_p$ при $p<2$. Таким образом, резольвентный предел операторов
$L_\varepsilon$ (если он существует) зависит от выбора последовательности
$u_\varepsilon(x)$ если сходимость понимается в $L_p$ при $p<2$. Это
конечно не означает, что определить оператор Штурма--Лувилля с
потенциалом $q(x)=u'(x)$, $u(x)\notin L_2$ невозможно.  Для "разумного"
определения таких операторов необходимо привлекать дополнительные
соображения (именно на таком пути в \S\,4 настоящей работы
определяются операторы с потенциалами вида $x^\alpha$). Если же говорить
об определении таких операторов посредством предельного перехода, то для
корректности этого способа необходимо "разумным образом" сузить класс
последовательностей гладких функций приближающих потенциал, чтобы
избежать неоднозначности в пределе.

\textbf{1.6. Примеры.} Классическими примерами потенциалов-распределений,
которые были ранее изучены, являются \(q(x)=\delta(x-x_0)\) и
\(q(x)=1/(x-x_0)\), где \(x_0\in(0,\pi)\). Приведем другие примеры.

{\bf Пример 1.} $l(y)=-y''+\frac{\sign x}{|x|^\alpha}y$, $x\in[-1,1]$.
С помощью канонической регуляризации (см. \cite{GSh}), функциям
$\frac{\sign x}{|x|^\alpha}$ можно сопоставить обобщенные функции, при
$\alpha\ne 2,4,6,\dots$.  При $\alpha<3/2$, полученные обобщенные функции
являются обобщенными производными функций из $L_2$ и значит можно
определить оператор, задаваемый выражением $l(y)$. Этот оператор при
$\alpha\ne 1$ определяется следующим образом.
$$
\begin{array}{l}
Ly=-y''+|x|^{-\alpha}\sign x\, y,\\
\mathfrak D(L)\hskip-2pt=\left\{y\in W_1^1\left|
\begin{array}{l}
 l(y)\in L_2,\ \  y(\pm1)=0\\
y(x)-y(0)\left(1+\frac{x}{(1-\alpha)(2-\alpha)}|x|^{1-\alpha}\right)\in
W^2_1[-1,1]
\end{array}\right\}\right. ,
\end{array}
$$
здесь предполагается, что $l(y)\in L_2$ на каждом из интервалов
$(-1,0)$ и $(0,1)$.
При $\alpha=1$ последнее соотношение нужно заменить следующим
$$
y(x)-y(0)(1+x\ln|x|)\in W^2_1[-1,1].
$$
Такое же определение оператора $L$ при $\alpha=1$ дано в [14], [15].
\hfill\break
{\bf Пример 2.}  Предложенный метод позволяет определить оператор
Штурма--Лиувилля для потенциалов $q(x)$, имеющих сколь угодно высокую
сингулярность во внутренней точке, при условии, что сильный рост потенциала
компенсируется сильной осцилляцией. В качестве примера положим
$q(x)=x^{-3}(\exp x^{-4})\sin(\exp x^{-4})$, $x\in[-1,1]$. Имеем
$$
\int q(\xi)d\xi=4^{-1}\cos(\exp x^{-4})\in L_2,
$$
поэтому оператор Штурма--Лиувилля для такой функции $q(x)$ вполне определен.

\section*{{\bf \S 2.\ \ Асимптотика собственных значений и собственных
функций для регулярных краевых условий}}\refstepcounter{section}
\addcontentsline{toc}{cont}{{\bf часть II.1}}

Здесь мы получим асимптотические формулы для решений уравнения
\begin{equation}
\label{main2}
-y''+q(x)y=\la y
\end{equation}
при $q(x)\in W_2^{-1}$, пользуясь которыми найдем
главные члены в асимптотических формулах для
собственных значений и собственных функций построенного в \S\,1
дифференциального оператора \(L=-d^2/dx^2+q\)
с регулярными краевыми условиями.
В действительности, в полученных в этом параграфе асимптотиках для
решений выделены не только главные, но и вторые члены. Именно в этом
состоит техническая сложность соответствующих результатов. Однако, запись
вторых членов в асимптотиках собственных значений и собственных функций
приводит к громоздким формулировкам. Формулировки можно упростить, если
ограничиться краевыми условиями частного вида. Эту работу мы проведем в
\S 3, а здесь, для общих условий, ограничимся главными членами.
В конце этого параграфа мы докажем
теорему о базисности Рисса собственных и присоединенных функций таких
операторов.

\textbf{2.1. Обозначения.}
Доказательство результатов этого параграфа потребует серьезной
технической
работы. Для удобства записи будем использовать следующие
обозначения:
$$
b(c,x,\lambda):=\intl_0^x u(t)\sin(2c+2\lambda^{1/2}t)dt, \quad
a(c,x,\lambda):=\intl_0^x u(t)\cos(2c+2\lambda^{1/2}t)dt,
$$
$$
B(c,x,\lambda):=\intl_0^x u^2(t)\sin(2c+2\lambda^{1/2}t)dt,\quad
A(c,x,\lambda):=\intl_0^x u^2(t)\cos(2c+2\lambda^{1/2}t)dt,
$$
$$
U(x):=\intl_0^x u^2(t)dt,\quad
\w(c,x,\lambda):=\intl_0^x\intl_0^tu(t)u(s)\cos(2c+2\lambda^{1/2}t)
\sin(2c+2\lambda^{1/2}s)dsdt,
$$
$$
\upsilon(c,x,\lambda):=b(c,x,\lambda)+\frac12\lambda^{-1/2}U(x)+
2\w(c,x,\lambda)-\frac12\lambda^{-1/2}\a(c,x,\lambda),
$$
$$
\Upsilon(c,\lambda):=\sup\limits_{0\le
x\le\pi}\left(|b(c,x,\lambda)|+
|a(c,x,\la|+2|\w(c,x,\la)|+\frac12|\la^{-1/2}A(c,x,\la)|\right)+
|\la|^{-1/2}\|u\|^2_{L_2}.
$$
Далее мы будем проводить оценки в комплексной $\la$--плоскости
внутри
областей, ограниченных параболами
$$
P_\al=\{\la\in\C\,\vert\ \Re\la>1,\ |\Im\sqrt\la|<\al\},
$$
а в $z$--плоскости ($z=\sqrt\la$) внутри полуполос
$$
\Pi_\al=\{z\in\C\,\vert\ \Re z>1,\ |\Im z|<\al\}.
$$
Всюду далее, рассматривая функцию $z=\sqrt{\la}$, подразумеваем выбор
ее главной ветви, принимающей положительные значения при $\la>0$.

Введенные выше функции будут участвовать в асимптотиках решений,
причем
$\Upsilon^2(\la)$ будет служить для оценки остатков. Функцию
$\Upsilon(\la)$ можно заменить более простым выражением
$$
\Upsilon_1(\la)=\sup\limits_{0\le x\le\pi}\left(|b(c,x,\la)|+
|a(c,x,\la)|\right)+|\la|^{-1/2}\|u\|^2_{L_2},
$$
так как внутри парабол $P_\al$ легко получить оценку
$$
\Upsilon_1(\la)\le\Upsilon(\la)\le M\Upsilon_1(\la)
$$
с постоянной $M$, зависящей только от $\al$ и $u(x)$. Нам, однако,
будет
удобнее оперировать с функцией $\Upsilon(\la)$.

Положим
$$
\Wo^1_2[0,\pi]=\{f(x)\vert\,f(x)\in W^1_2[0,\pi],\ f(0)=f(\pi)=0\},
$$
где $W^1_2[0,\pi]=W^1_2$ --- соболевское пространство. Через
$$
\begin{array}{ccc}
W^\theta_2=\left[W^1_2,L_2\right]_\theta,\quad
0\le\theta\le1 ,\\\phantom{.}\\
\Wo^\theta_2=\left[\Wo^1_2,L_2\right]_\theta,\quad
0\le\theta\le1 ,
\end{array}
$$
обозначаем интерполяционные пространства. Отметим (см.,
например~\cite[гл. 4.3.3]{LM}), что $W_2^\theta=\Wo_2^\theta$ при
$0\le\theta<1/2$.  Через $l_p^\theta$ обозначим пространства
последовательностей
$$
l_p^\theta=\left\{x=(x_1,x_2,x_3,\dots)\left\vert\
\|x\|_{\theta,p}^p=\suml_{n=1}^\infty |n^\theta
x_n|^p<\infty\right\}\right.,\quad 0\le\theta\le1,\ p\ge1.
$$
Последовательность точек $\{z_n\}$ назовем {\it
несгущающейся}, если найдется число $N$ такое, что внутри любого
круга
$K(z,r)$ с центром в точке $z$ фиксированного радиуса $r$ лежит не
более
$N$ точек последовательности ($N=N(r)$ не зависит от $z\in\C$).

\textbf{2.2. Асимптотика функций Прюфера.}
Мы уже показали в \S\ 1, что уравнение~\eqref{main2} можно записать
в
виде системы
\begin{equation}
\left(\begin{array}{cc} y_1 \\ y_2\end{array}\right)'
=\left(\begin{array}{cc} u&1\\-\lambda-u^2&-u \end{array}\right)
\left(\begin{array}{cc}
y_1\\y_2 \end{array}\right) ,\quad y_1=y,\  y_2=y'-u(x)y.
\label{sys}
\end{equation}
Сделаем замену
\begin{equation}
y_1(x,\lambda)=r(x,\la)\sin\theta(x,\la),\
y_2(x,\lambda)=\lambda^{\frac1{2}}r(x,\la)\cos\theta(x,\la)
\label{pr}
\end{equation}
(ее можно трактовать как переход к обобщенным полярным координатам),
которая является модификацией замены Прюфера
(см. \cite{Ha}).
Тогда систему~\eqref{sys} можно записать в виде
\begin{equation}
\begin{array}{ll}
r'\sin\theta+r\theta'\cos\theta&=
ur\sin\theta+\la^{\frac12}r\cos\theta,\\\hskip-40pt
\la^{\frac12}r'\cos\theta-\la^{\frac12}r\theta'\sin\theta&=
-\la r\sin\theta-u^2r\sin\theta-\la^{\frac12}ur\cos\theta,
\end{array}
\label{polar}
\end{equation}
где $r=r(x,\la)$, $\theta=\theta(x,\la)$, $u=u(x)$, а производные
функций
$r$ и $\theta$ берутся по переменной $x$.
Умножим первое уравнение в~\eqref{polar} на
$\la^{\frac12}\cos\theta$
и вычтем второе уравнение, умноженное на $\sin\theta$. В результате
получим уравнение для функции $\theta(x,\lambda)$
\begin{equation}
\theta'(x,\lambda)=\lambda^{\frac1{2}}+
\lambda^{-\frac1{2}}u^2(x)\sin^2\theta(x,\lambda)+
u(x)\sin2\theta(x,\lambda).
\label{th}
\end{equation}
Если мы сложим первое уравнение в~\eqref{polar}, умноженное на
$\la^{\frac12}\sin\theta$ со вторым уравнением, умноженным на
$\cos\theta$, то получим уравнение на функцию $r(x,\lambda)$
\begin{equation}
r'(x,\la)=-r(x,\la)\left[u(x)\cos 2\theta(x,\la)+\frac12\la^{-
1/2}u^2(x)
\sin 2\theta(x,\la)\right].
\label{r}
\end{equation}
Таким образом, от системы~\eqref{sys} мы перешли к системе двух
уравнений~\eqref{th}, \eqref{r}. Эти уравнения не являются линейными. Но
достоинство новой системы состоит в том, что уравнение~\eqref{th} не
содержит
неизвестной функции $r(x,\la)$ и является независимым
дифференциальным
уравнением на функцию $\theta(x,\la)$. Наша ближайшая цель ---
найти асимптотические представления для функций
$\theta(x,\lambda)$ и $r(x,\lambda)$. Эти формулы будут получены в
нижеследующих леммах. Первая из этих лемм является ключевой. На ней
базируются результаты как настоящего, так и следующего параграфов. Далее
без напоминаний используем введенные в п.2.1 обозначения.
\begin{Lemm}\label{lem:2.1}
Пусть $\al>0$ --- произвольное фиксированное число, $P_\al$ ---
область, ограниченная параболой $|\Im\sqrt{\lambda}|<\al$. Тогда
существует число $\mu$ (зависящее только от $u(x)$ и $\al$) такое,
что при
любых $c\in\R$ и $\la\in P_\al$, $\Re\la>\mu$,
уравнение \eqref{th} имеет единственное решение $\theta(x,\la)$,
определенное при всех $0\le x\le\pi$ и удовлетворяющее начальному
условию $\theta(0,\lambda)=c$. Это решение допускает представление
\begin{equation}
\theta(x,\lambda)=c+\lambda^{1/2}x+b(c,x,\lambda)+\frac12\lambda^{-
1/2}U(x)
+2\w(c,x,\lambda)-\frac12\lambda^{-
1/2}\a(c,x,\lambda)+\rho(c,x,\lambda),
\label{asth}
\end{equation}
где 
$$
|\rho(c,x,\la)|\le M\Upsilon^2(\lambda),\quad \la\in P_\la,\quad
\Re\la>\mu,
$$
причем постоянная $M$ зависит от $u(x)$ и $\al$, но
не зависит от $c$, $x$, $\la$.
\end{Lemm}

\begin{proof}
Разобьем его на несколько этапов.\hfill\break{\it Этап 1.}
Перепишем уравнение \eqref{th} в интегральном виде
\begin{equation}
\label{integr}
\theta(x,\lambda)=c+\lambda^{1/2}x+\intl_0^xu(t)\sin2\theta(t,\lambda)dt+
\frac12\lambda^{-1/2}\intl_0^xu^2(t)dt-
\frac12\lambda^{-1/2}\intl_0^xu^2(t)\cos2\theta(t,\lambda)dt.
\end{equation}
Обозначим правую часть этого уравнения через $F(\theta)$ и будем
решать
это уравнение методом последовательных приближений.
$$
\begin{array}{ccc}
\theta_0(x,\la)=c+\la^{1/2}x,\\
\theta_1(x,\la)=F(\theta_0)
=\theta_0(x,\la)+\intl_0^x u(t)\sin2\theta_0(t,\la)dt+\\+
\dfrac12\la^{-1/2}\intl_0^xu^2(t)dt-\dfrac12\la^{-
1/2}\intl_0^xu^2(t)
\cos2\theta_0(t,\la)dt,
\end{array}
$$
или, используя введенные обозначения,
$$
\theta_1(x,\la)=\theta_0(x,\la)+b(c,x,\la)+\frac12\la^{-1/2}U(x)-
\frac12\la^{-1/2}\a(c,x,\la).
$$
Далее мы выпишем приближение $\theta_2$, выделим в нем
слагаемое $2\w(c,x,\la)$, которое входит в представление~\eqref{asth}
и оценим остаток. Заметим, что
$$
|\theta_1(x,\la)-\theta_0(x,\la)|\le\Upsilon(\la),
$$
причем в силу леммы Римана--Лебега
$\Upsilon(\la)\to0$ при $\la\to\infty$, $\la\in P_\al$. Поэтому
найдется число $\mu>1$ такое, что при $\la\in P_\al$, $\Re\la>\mu$
выполняются оценки
\begin{equation}
\label{fv1}
\begin{array}{cc}
\left|\sin2\left(\theta_1-\theta_0\right)-
2\left(\theta_1-
\theta_0\right)\right|\le2\Upsilon^3(\la)<\Upsilon^2(\la),\\
\left|\cos2\left(\theta_1-\theta_0\right)-1\right|
\le3\Upsilon^2(\la).
\end{array}
\end{equation}
Воспользовавшись тригонометрическими формулами
$$
\begin{array}{cc}
\sin2\theta_1=\sin2\theta_0\cos2(\theta_1-\theta_0)+
\cos2\theta_0\sin2(\theta_1-\theta_0),\\
\cos2\theta_1=\cos2\theta_0\cos2(\theta_1-\theta_0)-
\sin2\theta_0\sin2(\theta_1-\theta_0)
\end{array}
$$
и неравенствами~\eqref{fv1}, получим
\begin{equation}
\label{fv1+}
\begin{array}{cc}
\left|\sin2\theta_1-[\sin2\theta_0+
2(\theta_1-\theta_0)\cos2\theta_0]\right|\le 4M_1\Upsilon^2(\la),\\
\left|\cos2\theta_1-[\cos2\theta_0-
2(\theta_1-\theta_0)\sin2\theta_0]\right|\le 4M_1\Upsilon^2(\la),
\end{array}
\end{equation}
где $M_1=\ch2\pi\al$ (эта константа возникает из оценок
$|\sin2\theta_0|<M_1$, $|\cos2\theta_0|<M_1$ при $\la\in P_\al$).
Теперь получаем
\begin{equation}
\label{fv2}
\begin{array}{cc}
\theta_2(x,\la)=F(\theta_1)=F(\theta_0)
+2\intl_0^xu(t)\left(\theta_1(t,\la)-\theta_0(t,\la)\right)
\cos2\theta_0(t,\la)dt-\\
-\la^{-1/2}\intl_0^xu^2(t)\left(\theta_1(t,\la)-
\theta_0(t,\la)\right)
\sin2\theta_0(t,\la)dt+\rho(x,\la),
\end{array}
\end{equation}
где остаток $\rho(x,\la)$ в силу неравенств~\eqref{fv1+} подчинен
оценке
\begin{equation}\label{2.11'}
|\rho(x,\la)|\le2\left(2\|u\|_{L_1}+\|u\|^2_{L_2}\right)M_1\Upsilon^
2(\la).
\end{equation}
Второй интеграл в правой части равенства~\eqref{fv2} допускает
аналогичную оценку
\begin{equation}\label{2.11''}
\left|\la^{-1/2}\intl_0^xu^2(t)\left(\theta_1(t,\la)-
\theta_0(t,\la)\right)
\sin2\theta_0(t,\la)dt\right|\le
\left(|\la|^{-
1/2}\|u\|_{L_2}^2\right)M_1\Upsilon(\la)<M_1\Upsilon^2(\la).
\end{equation}
Вычислим первый интеграл. Имеем
\begin{equation}\label{2-1}
\begin{array}{cc}
2\intl_0^xu(t)\left(\theta_1(t,\la)-\theta_0(t,\la)\right)
\cos2\theta_0(t,\la)dt=2\w(c,x,\la)+\\+
\la^{-1/2}\intl_0^xu(t)\cos2\theta_0(t,\la)\intl_0^tu^2(s)dsdt
-\la^{-1/2}\intl_0^xu(t)\cos2\theta_0(t,\la)
\intl_0^tu^2(s)\cos2\theta_0(s,\la)dsdt.
\end{array}
\end{equation}
Для суммы второго и третьего слагаемых в правой части имеем оценку
\begin{equation}\label{es2}
\begin{array}{ccc}
\left|\la^{-1/2}\intl_0^xu(t)\cos2\theta_0(t,\la)\intl_0^tu^2(s)
(1-\cos2\theta_0(s,\la))dsdt\right|=\\
=|\la|^{-1/2}\left|\intl_0^xu^2(s)(1-\cos2\theta_0(s,\la))
\intl_s^xu(t)\cos2\theta_0(t,\la)dtds\right|\le\\
\le(1+M_1)\left(|\la|^{-1/2}\intl_0^x|u|^2(s)ds\right)
\sup\limits_{s,x\in[0,\pi]}
\left|\intl_s^xu(t)\cos2\theta_0(t,\la)dt\right|
\le2(1+M_1)\Upsilon^2(\la).
\end{array}
\end{equation}

Таким образом для приближения $\theta_2(x,\la)$ мы получили
представление
\begin{equation}
\label{fv5}
\theta_2(x,\la)=\theta_0(x,\la)+\upsilon(c,x,\la)+\rho_2(x,\la),
\end{equation}
где $|\rho_2(x,\la)|\le M\Upsilon^2(\la)$. Значение постоянной $M$ здесь
можно указать:
$$
M=\left(3+4\|u\|_{L_1}+2\|u\|^2_{L_2}\right)M_1+2,
$$
но далее через $M$ обозначаются различные постоянные, зависящие только от
$u(x)$ и $\al$.

{\it Этап 2.}
Здесь мы покажем, что для следующих приближений $\theta_3$ и $\theta_4$
справедливо такое же представление, как и для $\theta_2$.
Из представления~\eqref{fv5} имеем
\begin{equation}\label{fv7}
|\theta_2(x,\la)-\theta_0(x,\la)|<\Upsilon(\la)(1+M\Upsilon(\la))<
\sqrt{2}\Upsilon(\la),
\end{equation}
если число $\mu$ выбрано так, что $M\Upsilon(\la)<\sqrt{2}-1$
при $\la\in P_\al$, $\Re\la>\mu$. Поэтому справедливы
оценки~\eqref{fv1+}, в которых $\theta_1$ заменяется на $\theta_2$, а
число $4$ в правой части --- на $8$. С помощью этих оценок, также как
при выводе~\eqref{fv2}, получаем
\begin{equation}
\label{fv9}
\begin{array}{cc}
\theta_3(x,\la)=F(\theta_2)=F(\theta_0)+
2\intl_0^xu(t)\left(\theta_2(t,\la)-\theta_0(t,\la)\right)
\cos2\theta_0(t,\la)dt+\\
+\la^{-1/2}\intl_0^xu^2(t)\left(\theta_2(t,\la)-
\theta_0(t,\la)\right)
\sin2\theta_0(t,\la)dt+\rho(x,\la),
\end{array}
\end{equation}
где остаток $\rho(x,\la)$ подчинен оценке~\eqref{2.11'} с заменой $2$
на $4$.
Второй интеграл в правой части~\eqref{fv9} оценивается также, как
в~\eqref{2.11''}, т.е. величиной $2M_1\Upsilon^2(\la)$.
Перепишем представление~\eqref{fv5} в виде
$$
\theta_2(t,\la)-\theta_0(t,\la)=\theta_1(t,\la)-\theta_0(t,\la)+
2\w(t,\la)+\rho_2(t,\la).
$$
Тогда первый интеграл в правой части~\eqref{fv9} с учетом
равенства~\eqref{2-1} и оценок~\eqref{es2} предстанет в виде
\begin{equation}\label{FF}
\begin{array}{cc}
2\intl_0^xu(t)\left(\theta_2(t,\la)-\theta_0(t,\la)\right)
\cos2\theta_0(t,\la)dt=
2\w(c,x,\la)+\\
+4\intl_0^x\w(c,t,\la)\cos2\theta_0(t,\la)dt+\rho(x,\la),
\end{array}
\end{equation}
где $|\rho(x,\la)|<M\Upsilon^2(\la)$, $M=const$. Для оценки интеграла в
правой части~\eqref{FF} поменяем порядок интегрирования
$$
\begin{array}{cc}
4\intl_0^x\w(t,\la)\cos2\theta_0(t,\la)dt=
4\intl_0^x\intl_0^t\intl_0^su(t)u(s)u(\tau)
\cos2\theta_0(t,\la)\cos2\theta_0(s,\la)
\sin2\theta_0(\tau,\la)d\tau dsdt=\\
=4\intl_0^xu(\tau)\sin2\theta_0(\tau,\la)
\iint\limits_\D
u(t)u(s)\cos2\theta_0(t,\la)\cos2\theta_0(s,\la)dsdtd\tau,
\end{array}
$$
где внутренний интеграл берется по треугольнику
$$
\D=\{(s,t)\in\R\vert\ \tau\le s\le x,\ \tau\le t\le x,
\ s\le t\},
$$
который составляет половину квадрата $[\tau,x]\times[\tau,x]$.
Поэтому выписанный тройной интеграл оценивается величиной
\begin{equation}\label{fv15}
\begin{array}{cc}
2\left|\intl_0^xu(\tau)\sin2\theta_0(\tau,\la)
\intl_\tau^x\intl_\tau^x
u(t)u(s)\cos2\theta_0(t,\la)\cos2\theta_0(s,\la)dsdtd\tau\right|=\\
=\left|\intl_0^xu(\tau)\sin2\theta_0(\tau,\la)\left(
\intl_\tau^xu(t)\cos2\theta_0(t,\la)dt\right)^2d\tau\right|\le
M\Upsilon^2(\la).
\end{array}
\end{equation}
Тем самым, мы получили представление
$$
\theta_3(x,\la)=\theta_0(x,\la)+\upsilon(c,x,\la)+\rho_3(x,\la),
$$
где $|\rho_3(x,\la)|\le M\Upsilon^2(\la)$. Заметим теперь, что
неравенство~\eqref{fv7} сохраняется с заменой $\theta_2$ на $\theta_3$.
Но тогда можно полностью повторить проведенные на этом этапе рассуждения
и получить для функции $\theta_4(x,\la)$ такое же представление, как и
для $\theta_3(x,\la)$. Отсюда также вытекает оценка
\begin{equation}\label{Mres}
|\theta_4(x,\la)-\theta_3(x,\la)|+|\theta_3(x,\la)-\theta_2(x,\la)|\le
M\Upsilon^2(\la),
\end{equation}
которая понадобится в дальнейшем.\hfill\break
{\it Этап 3.} Покажем, что приближения $\theta_n$ сходятся к решению
$\theta$ уравнения~\eqref{integr}. Обозначив через $F(\theta)$ правую
часть уравнения~\eqref{integr}, перепишем его в виде
\begin{equation}\label{Phi}
f=F(\theta_0+f)-\theta_0=:\Phi(f),
\end{equation}
где $f=\theta-\theta_0$. Воспользовавшись тригонометрическими формулами
для вычисления выражения $F(\theta_0+f)-\theta_0$, получим
$$
\Phi(f)=\Phi_0+\Phi_1(f)+\Phi_2(f)+\Phi_3(f),
$$
где
$$
\begin{array}{ccc}
\Phi_0=\intl_0^xu(t)\sin2\theta_0(t,\la)dt-
\dfrac12\la^{-1/2}\intl_0^xu^2(t)\cos2\theta_0(t,\la)dt,\quad
\Phi_1(f)=2\intl_0^xf(t)u(t)\cos2\theta_0(t,\la)dt,\\
\Phi_2(f)=\intl_0^xu(t)\cos2\theta_0(t,\la)\left(\sin2f(t)-
2f(t)\right)dt+
\intl_0^xu(t)\sin2\theta_0(t,\la)\left(\cos2f(t)-1\right)dt,\\
\Phi_3(f)=\dfrac12\la^{-1/2}\left[\intl_0^x(1-\cos2f(t))u^2(t)
\cos2\theta_0(t,\la)dt+
\intl_0^x\sin2f(t)u^2(t)\sin2\theta_0(t,\la)dt\right].
\end{array}
$$
Покажем, что отображение $\Phi^2$ является сжимающим в шаре достаточно
малого радиуса с центром в нуле (шар рассматривается в пространстве
$C[0,\pi]$). Заметим, что $\Phi_0$ не зависит от $f$,
а $\Phi_1$ --- линейный оператор, причем
$$
\begin{array}{cc}
\Phi_1^2f=4\intl_0^x\intl_0^tu(t)\cos2\theta_0(t,\la)u(s)\cos2\theta_0(s,
\la)f(s)dsdt=\\
=4\intl_0^xf(s)u(s)\cos2\theta_0(s,\la)
\intl_s^xu(t)\cos2\theta_0(t,\la)dtds.
\end{array}
$$
Следовательно,
$$
\|\Phi^2_1\|_C\le8\|u\|_{L_1}M_1\Upsilon(\la).
$$
Далее, нетрудно видеть, что отображение $\Phi_2$ обладает свойством
$$
\|\Phi_2(f)-\Phi_2(g)\|_C\le M(\|f\|_C+\|g\|_C)\|f-g\|_C,\quad
\text{если}\ \|f\|_C+\|g\|_C\le1.
$$
Для отображения $\Phi_3$ имеем
$$
\|\Phi_3(f)-\Phi_3(g)\|_C\le \mu^{-1/2}M\|f-g\|_C,
$$
если $\la\in P_\al$ и $\Re\la>\mu$. Выберем теперь число $\mu$ столь
большим, а радиус шара $r$ столь малым, что все
попарные произведения $\Phi_j\Phi_k$, $1\le j,\,k\le3$,
являются сжимающими отображениями в выбранном шаре с коэффициентами
сжатия не превосходящими $1/18$. Следовательно,
\begin{equation}\label{1/2}
\|\Phi^2(f)-\Phi^2(g)\|_C\le\dfrac12\|f-g\|_C,\quad\text{если}\
\|f\|_C\le r,\ \|g\|_C\le r.
\end{equation}
В частности,
$$
\|\Phi^2(f)\|_C\le\dfrac12\|f\|_C+\|\Phi^2(0)\|_C\le\dfrac12\|f\|_C+
M\Upsilon(\la).
$$
Последнее неравенство вытекает из того, что
$\|\Phi(0)\|_C=\|\Phi_0\|_C\le\Upsilon(\la)$, и того, что отображение
$\Phi$ ограничено некоторой константой $M$. Увеличив, если нужно, число
$\mu$, можно добиться выполнения оценки $\Upsilon(\la)<r/(2M)$.  Поэтому
отображение $\Phi^2$ переводит шар радиуса $r$ в себя.

Положим $f_0=0$. Квадрат отображения $\Phi$ сжимает, а потому
последовательность $f_n=\Phi(f_{n-1})$, $n=1,2,\dots$, сходится к решению
$f$ уравнения $f=\Phi(f)$.  Согласно~\eqref{Phi}, последовательность
$\theta_n=\theta_0+f_n$ сходится к решению $\theta$ уравнения
$\theta=F(\theta)$. Следовательно,
$$
\theta-\theta_2=\suml_{n=2}^\infty(\theta_{n+1}-\theta_n)=
\suml_{n=2}^\infty(f_{n+1}-f_n)=\suml_{k=0}^\infty
(\Phi^{2k}(f_3)-\Phi^{2k}(f_2))+\suml_{k=0}^\infty
(\Phi^{2k}(f_4)-\Phi^{2k}(f_3)).
$$
Поэтому из оценок~\eqref{1/2} и~\eqref{Mres} имеем
$$
|\theta-\theta_2|\le 2|f_3-f_2|+2|f_4-f_3|=2|\theta_3-\theta_2|+
2|\theta_4-\theta_3|\le M\Upsilon^2(\la).
$$
Эта оценка, с учетом~\eqref{fv5} влечет справедливость~\eqref{asth}.
Лемма доказана.
\end{proof}

\begin{Lemm}\label{lem:2.2}
Пусть $\al>0$ --- произвольное фиксированное число, а $P_\al$ ---
область, ограниченная параболой $|\Im\sqrt{\la}|<\al$. Пусть
$\theta(x,\la)$ --- решение уравнения~\eqref{th} с начальным условием
$\theta(0,\la)=c\in\R$.  Тогда решение $r(x,\la)$ уравнения~\eqref{r} с
начальным условием $r(0,\la)=c_1$ допускает представление
\begin{equation}
r(x,\lambda)=c_1\left[1-a(c,x,\lambda)-
\frac12\lambda^{-1/2}\b(c,x,\lambda)\right]+
\rho(c,x,\lambda)\quad\la\in P_\al,\ \Re\la>\mu,
\label{asr}
\end{equation}
где остаток $\rho(c,x,\lambda)$ подчинен оценке
\begin{equation}
|\rho(c,x,\lambda)|\le M\Upsilon^2(\la).
\label{est}
\end{equation}
Здесь $\mu$ и $M$ --- числа, зависящее только от $\al$ и $u(x)$.
\end{Lemm}

\begin{proof}
Уравнение \eqref{r} решается явно
\begin{equation}
r(x,\lambda)=c_1\exp\left(-\intl_0^xu(t)\cos(2\theta(t,\lambda))dt-
\frac12\lambda^{-1/2}\intl_0^xu^2(t)\sin(2\theta(t,\lambda))dt\right).
\label{2.2a}
\end{equation}
Докажем оценку
\begin{equation}
\left|\intl_0^x u(t)\cos2\theta(t,\la)dt-a(c,x,\la)\right|\le
M\Upsilon^2(\la),\quad\la\in P_\al,\ \Re\la>\mu.
\label{2.2b}
\end{equation}
Вновь положим $\theta_0(x,\la)=c+\la^{1/2}x$. В силу леммы~\ref{lem:2.1}
имеем
$$
|\theta(x,\la)-\theta_0(x,\la)|\le\Upsilon(\la)+M\Upsilon^2(\la)<
2\Upsilon(\la),\quad
\la\in P_\al,\ \Re\la>\mu,
$$
если число $\mu$ выбрано достаточно большим.
Тогда с учетом~\eqref{fv1+} (где вместо $\theta_1$ берем $\theta$)
получаем
$$
\intl_0^x u(t)\cos2\theta(t,\la)dt=
a(c,x,\la)-2\intl_0^xu(t)\left(\theta(t,\la)-\theta_0(t,\la)\right)
\sin2\theta_0(t,\la)dt+\rho(c,x,\la),
$$
где
$|\rho(c,x,\la)|\le M\Upsilon^2(\la)$.
Покажем, что интеграл в правой части оценивается величиной
$M\Upsilon^2(\la)$. Имеем
$$
\theta(t,\la)-\theta_0(t,\la)=\dfrac12\la^{-1/2}U(x)-
\dfrac12\la^{-1/2}\a(c,x,\la)+2\w(c,x,\la)+b(c,x,\la)+\rho(c,x,\la),
$$
где $|\rho|\le M\Upsilon^2(\la)$. Теперь нам необходимо оценить четыре
интеграла, в которых вместо $\theta-\theta_0$ участвует любая из
функций правой части. Оценка интеграла с функцией $\la^{-1/2}U(x)$
проводится также как в~\eqref{es2}. Оценка интеграла с функцией
$\la^{-1/2}\a(c,x,\la)$ проводится аналогично. Оценка с функцией
$\w(c,x,\la)$ проводится также как в~\eqref{fv15}. Интеграл с функцией
$b(c,x,\la)$ преобразуем к виду
$$
\intl_0^x\intl_0^tu(t)u(s)\sin2\theta_0(t,\la)
\sin2\theta_0(s,\la)dsdt=\dfrac12
\left(\intl_0^xu(t)\sin2\theta_0(t,\la)dt\right)^2\le
\dfrac12\Upsilon^2(c,\la).
$$
Здесь мы воспользовались свойством симметрии подынтегральной функции
относительно диагонали $s=t$.

Тем самым, оценка~\eqref{2.2b} доказана. На таком же пути (но только
проще) получается оценка
$$
\left|\la^{-1/2}\intl_0^xu^2(t)\sin2\theta(t,\la)dt-
\la^{-1/2}\b(c,x,\la\right|<M\Upsilon^2(\la).
$$
Следовательно, выражение под знаком экспоненты в~\eqref{2.2a} можно
записать в виде
$$
f(c,x,\la):=a(c,x,\la)-\dfrac12\la^{-1/2}\b(c,x,\la)+\rho(c,x,\la).
$$
Воспользовавшись оценками
$$
|f|<M\Upsilon(\la),\quad|\rho|<M\Upsilon^2(\la),\quad
\exp{f}=1+f+O(f^2),
$$
получаем представление~\eqref{asr} с оценкой остатка~\eqref{est}.
\end{proof}

\textbf{2.3. Асимптотики собственных значений и собственных функций.}
Как было сказано в начале этого парагарафа, для регулярных краевых
условий общего вида мы ограничимся главными членами асимптотик и оценками
остатков.

Напомним, что определения пространств $\Wo_2^\theta$, $W_2^\theta$ и
$l_2^\theta$ были даны в п.2.1.
\begin{Lemm}\label{lem:2.3}
Пусть несгущающаяся последовательность комплексных чисел
$\{z_n\}$ лежит в полосе $|\Im z|<\al$, где $\al>0$ --- произвольное
число. Пусть оператор $T_x:\,W_2^\theta\to l_2^\theta$,
$0\le\theta<1/2$ определен равенством
$$
T_xf=\{c_n\}_{n=1}^\infty,\quad\text{где}\ c_n=\intl_0^xf(t)e^{iz_nt}dt.
$$
Тогда при любом $0\le\theta<1/2$ оператор $T_x$ ограничен и его норма
зависит от последовательности $\{z_n\}$ и $\theta$ (если
$\theta\to1/2$), но не зависит от $x\in[0,\pi]$.
\end{Lemm}
\begin{proof}
При $\theta=0$ имеем $\Wo^\theta_2=L_2$. В этом случае ограниченность
оператора $T$ известна (см., например~\cite{Kaz}). Если
$f(t)\in\Wo_2^1[0,\pi]$, то интеграл $\intl_0^\pi f(x)e^{iz_nx}dx$ можно
проинтегрировать по частям и воспользоваться неравенством
$|z_n|>\delta_n$ ($\delta>0$ не зависит от $n$), вытекающим из
определения несгущающейся последовательности. Тогда получим, что оператор
$$
T_\pi:\, \Wo_2^1\to\,l_2^1
$$
ограничен. Из теоремы об интерполяции (см.,
например~\cite[гл. 1]{LM}) следует, что оператор
$$
T_\pi:\, \Wo_2^\theta\to\,l_2^\theta,\quad0\le\theta\le1
$$
ограничен. Но пространства $\Wo_2^\theta$ и $W_2^\theta$ при
$0\le\theta<1/2$ совпадают. Кроме того, оператор умножения на
$\chi_{[0,x]}$ --- характеристическую функцию отрезка $[0,x]$ в
пространстве $W_2^\theta$ при $0\le\theta<1/2$ ограничен и его норма не
зависит от $x$ (см., например,~\cite[теорема 17.11]{BIN}). Зависимость от 
$\theta$ нормы этого оператора существенна только при $\theta\to1/2$.  
Это влечет утверждение леммы.  
\end{proof} 
\begin{Lemm}\label{lem:2.3+} 
Если $\{z_n\}_{n=1}^\infty$ --- несгущающаяся последовательность чисел в 
полуполосе $\Pi_\al$, а функция $u(x)\in W_2^\theta$ при некотором
$0\le\theta<1/2$, то при любом $c\in\R$ последовательность
$\{\Upsilon(c,z_n^2)\}_{n=1}^\infty\in l_2^\theta$.
\end{Lemm}
\begin{proof}
Из леммы~\ref{lem:2.3} следует
$$
\{|a(c,x,z_n^2)|+|b(c,x,z_n^2)|\}_{n=1}^\infty\in l_2^\theta,
$$
причем норма этой последовательности в $l_2^\theta$ не зависит от $x$.
Применяя лемму~\ref{lem:2.3} к внутреннему интегралу в представлении
функции $\w(c,x,\la)$, получаем также
$\{\w(c,x,z_n^2\}_{n=1}^\infty\in l_2^\theta$. Наконец, заметим, что
$|z_n|>\delta n$ при некотором $\delta>0$ для несгущающейся
последовательности $\{z_n\}_{n=1}^\infty$ и
$\{n^{-1}\}_{n=1}^\infty\in l_2^\theta$ при $0\le\theta<1/2$. Отсюда
получаем утверждение леммы.
\end{proof}
\begin{Lemm}\label{lem:2.4}
Пусть $u(x)=\intl q(\xi)d\xi\in W_2^\theta$ при некотором
$0\le\theta<1/2$. Обозначим
через $\Phi(x,\lambda)$ и $\Psi(x,\lambda)$ пару  решений
уравнения~\eqref{main2}, подчиненных начальным условиям
$\Phi(0,\lambda)=1$, $\Phi^{[1]}(0,\lambda)=0$, $\Psi(0,\lambda)=0$,
$\Psi^{[1]}(0,\lambda)=1$. Тогда внутри любой параболы $P_\al$
справедливы представления
\begin{equation}\label{z1}
\begin{array}{cc}
\Phi(x,\lambda)=\cos(\lambda^{1/2}x)+\varphi(x,\la),\quad
&\Phi^{[1]}(x,\lambda)=-\lambda^{1/2}\sin(\lambda^{1/2}x)
+\lambda^{1/2}\varphi_1(x,\la),\\
\Psi(x,\lambda)=\lambda^{-1/2}\sin(\lambda^{1/2}x)+
\lambda^{-1/2}\psi(x,\la), \quad
&\Psi^{[1]}(x,\lambda)=\cos(\lambda^{1/2}x)+\psi_1(x,\la),
\end{array}
\end{equation}
где каждая из функций $\varphi$, $\varphi_1$, $\psi$, $\psi_1$
мажорируется величиной $M\Upsilon(\la)\to0$ при $\la\to\infty$,
где $\Upsilon(\la)=\Upsilon(\pi/2,\la)$. В частности,
если $\{z_n\}$ --- несгущающаяся последовательность, лежащая в полуполосе
$\Pi_\al$, то последовательность $\la_n=z_n^2$ лежит в $P_\al$,
при этом $\{\Upsilon(\la_n)\}\in l_2^\theta$
а потому последовательности $\{\varphi(x,\la_n\}$,
$\{\varphi_1(x,\la_n\}$, $\{\psi(x,\la_n\}$, $\{\psi_1(x,\la_n\}$
принадлежат пространству $l_2^\theta$ и их нормы в этом пространстве
не зависят от $x$.
\end{Lemm}
\begin{proof}
Согласно определению~\eqref{pr} функций Прюфера имеем
$$
\Phi(x,\la)=r(c_1,x,\la)\sin\theta(c,x,\la),\quad
\Phi^{[1]}(x,\la)=r(c_1,x,\la)\cos\theta(c,x,\la).
$$
Постоянные $c$, $c_1$ определяются из начальных условий
$$
1=c_1\sin c,\quad 0=c_1\cos c,
$$
откуда находим $c=\pi/2$, $c_1=1$.

В силу лемм~\ref{lem:2.1} и~\ref{lem:2.2} имеем
$$
\begin{array}{cc}
\Phi(x,\la)=(1+\omega_1(x,\la))\sin\left(\frac\pi2+\la^{1/2}x+
\omega_2(x,\la)\right)=\\=\cos\left(\la^{1/2}x+\omega_2(x,\la)\right)+
\omega_1(x,\la)\cos\left(\la^{1/2}x+\omega_2(x,\la)\right),
\end{array}
$$
где функции $|\omega_1|$, $|\omega_2|$ мажорируются величиной
$M\Upsilon(\la)$. Тогда функция
$\cos\left(\la^{1/2}x+\omega_2(x,\la)\right)$ ограничена в области
$P_\al$, о потому функция $\Phi(x,\la)-\cos\left(\la^{1/2}x\right)$
мажорируется величиной $M\Upsilon(\la)$. Остается воспользоваться
леммой~\ref{lem:2.3+} Доказательства других асимптотических
формул~\eqref{z1} аналогичны.
\end{proof}
\begin{Theorem}\label{tm:2.5}
Пусть $u(x)\in W_2^\theta$ при некотором $0\le\theta<1/2$, а $q(x)=u'(x)$
в смысле теории распределений. Пусть $L$ --- построенный в
теореме~\ref{tm:1.4} оператор, порожденный дифференциальным
выражением $-y''+q(x)y$ и регулярными краевыми условиями (т.е.
выполнено одно из условий 1)\,--\,3) теоремы~\ref{tm:1.4}).
Обозначим через $\{\la_n\}_{n=1}^\infty$ собственные значения
оператора $L$, а через $\{\la^0_n\}_{n=1}^\infty$ ---
собственные значения оператора $L_0=-\dfrac{d^2}{dx^2}$ с такими же
краевыми условиями, в которых квазипроизводные заменены на обычные
производные. Для обоих операторов нумерацию собственных значений
проводим в порядке возрастания модулей и с учетом алгебраической
кратности. Тогда при дополнительном предположении усиленной
регулярности краевых условий справедливы равенства
\begin{equation}\label{evas}
\sqrt{\lambda_n}=\sqrt{\la^0_n}+s_n,\ n=1,2,\dots\,,\quad\text{где}\
\{s_n\}_{n=1}^\infty\in l^\theta_2.
\end{equation}
При невыполнении условия усиленной регулярности равенства остаются,
но условие $\{s_n\}_{n=1}^\infty\in l^\theta_2$ заменяется условием
$\{s_n^2\}_{n=1}^\infty\in l^\theta_2$.
\end{Theorem}

\begin{proof}
Сейчас мы проведем доказательство только части этой теоремы. А именно,
предположив, что $\{\la_n^0\}_{n=1}^\infty$ --- все собственные значения
оператора $L_0$, мы покажем, что существует серия собственных значений
$\{\la_n\}_{n=N}^\infty$, для которой справедливы равенства
$$
\sqrt{\la_n}=\sqrt{\la_n^0}+s_n,\ n=N,\,N+1,\dots\,,\quad
\{s_n\}_{n=N}^\infty\in l_2^\theta\ \left(\text{ или }
\{s_n^2\}_{n=N}^\infty\in l_2^\theta\ \right).
$$
В действительности, число оставшихся собственных значений оператора $L$
равно $N-1$ и их можно поставить во взаимнооднозначное соответствие с
$\{\la_n^0\}_{n=1}^{N-1}$. Однако, доказательство этого факта удобнее
здесь опустить и провести его в конце параграфа.

Для определенности рассмотрим только случай, когда выполнено условие 2)
теоремы~\ref{tm:1.4} (другие случаи легче). В силу теоремы~\ref{tm:1.4} и
леммы~\ref{lem:2.4} собственные значения оператора $L$ в области $P_\al$
определяются из уравнения
$$
\Delta(\la)=\cos\pi\la^{1/2}+J_0+\rho(\la)=0,\quad
J_0=(J_{12}+J_{34})(J_{14}-J_{23})^{-1},
$$
где $|\rho(\la)|\le M\Upsilon(\la)$ при $\Re\la>\mu$. Перепишем это
уравнение в виде
$$
\Delta(z^2)=\cos\pi z+J_0+\rho(z^2)=0.
$$
Пусть $\{z_n^0\}_{n=1}^\infty$ --- корни уравнения
$$
\Delta_0(z^2)=\cos\pi z+J_0=0
$$
из правой полуплоскости. Очевидно, $\la_n^0=(z_n^0)^2$. Выберем число
$\al>1$ столь большим, что все $\{z_n^0\}$ лежат в полуполосе
$\Pi_{\al-1}$. Если $J_0\ne\pm1$ (усиленная регулярность), то все нули
$\{z_n^0\}$ простые и найдутся числа $\delta>0$, $\eps>0$ такие, что
круги $K_n=\{x\vert\,|z-z_n^0|<\delta\}$ не пересекаются, причем
$|\Delta_0(z^2)|>\eps$ при $z\in\partial K_n$ (это утверждение легко
следует из периодичности рассматриваемой функции). Так как
$|\rho(z^2)|\to0$ при $z\to\infty$, $z\in\Pi_\al$, то в силу теоремы Руше
при всех достаточно больших $n\ge N$ функця $\Delta(z^2)$ имеет ровно
один нуль $z_n$ в круге $K_n$. Конечно, последовательность
$\{z_n\}_{n=N}^\infty$ является несгущающейся. Поэтому из
леммы~\ref{lem:2.3} имеем $\{\rho(z_n^2)\}\in l_2^\theta$, а тогда
$$
\{\cos\pi z_n-\cos\pi z_n^0\}_N^\infty\in l_2^\theta.
$$
Число $\delta$ можно считать выбранным столь малым, что
$\left|\left(\cos(\pi z)\right)'\right|>\eps_1>0$ при $z\in K_n$. Но
тогда
\begin{equation}\label{z_n}
|z_n-z_n^0|\le M|\cos\pi z_n-\cos\pi z_n^0|,
\end{equation}
где постоянная $M$ не зависит от $n$. Тем самым, асимптотические
равенства~\eqref{evas} доказаны.

В случае $J_0=\pm1$ корни $z_n^0$ имеют кратность $2$. Тогда можно
повторить проведенные рассуждения, но вместо неравенства~\eqref{z_n}
воспользоваться неравенством
$$
|z_n-z_n^0|^2\le M|\cos\pi z_n-\cos\pi z_n^0|,
$$
из которого получаем $\{s_n^2\}_N^\infty\in l_2^\theta$. Теорема
доказана.
\end{proof}

\begin{Theorem}\label{tm:2.6}
Пусть выполнены условия и сохранены обозначения теоремы~\ref{tm:2.5}.
Обозначим через $\{y_n(x)\}_{n=1}^\infty$ и $\{y^0_n(x)\}_{n=1}^\infty$
нормированные в пространстве $L_2[0,\pi]$ собственные функции операторов
$L$ и $L_0$, отвечающие собственным значениям $\{\la_n\}_{n=1}^\infty$
и $\{\la^0_n\}_{n=1}^\infty$ соответственно. Тогда при дополнительном
условии усиленной регулярности справедливы равенства
$$
y_n(x)=y_n^0(x)+\psi_n(x),\qquad y^{[1]}_n(x)=(y_n^0(x))'+n\psi^1_n(x),
$$
где
$$
\suml_{n=1}^\infty n^\theta\left(|\psi_n(x)|^2+|\psi^1_n(x)|^2\right)
\le C
$$
и постоянная $C$ не зависит от $x$. Если условие усиленной регулярности
не выполнено, то утверждение сохраняется, но квадраты функций $\psi_n$ и
$\psi^1_n$ в последней формуле надо заменить на четвертые
степени, а число $\theta$ заменить на $2\theta$.
\end{Theorem}
\begin{proof}
Снова рассмотрим наиболее трудный случай, когда выполнено условие 2)
теоремы~\ref{tm:1.4}. В этом случае после нормировки, одно из краевых
условий (будем считать, что второе) имеет нулевой порядок, т.е.
$$
U_2(y)=\beta y(0)+\gamma y(\pi),\quad |\beta|+|\gamma|>0.
$$
Тогда собственные функции определяются равенством $y_n(x)=y(x,\la_n)$, где
\begin{equation}\label{(y)}
y(x,\la)=\begin{vmatrix}\Phi(x,\la)&\Psi(x,\la)\\
U_2(\Phi)&U_2(\Psi)\end{vmatrix}=\begin{vmatrix}\Phi(x,\la)&\Psi(x,\la)\\
\beta+\gamma\Phi(\pi,\la)&\gamma\Psi(\pi,\la)\end{vmatrix}.
\end{equation}
Это так, поскольку $y(x,\la)$ --- решение уравнения~\eqref{main2} и
$U_j(y(x,\la_n))=0$, $j=1$, $2$. Используя обозначения
$\sqrt{\la_n}=z_n$, $\sqrt{\la^0_n}=z^0_n$, получаем
$$
y_n^0(x)=\gamma(\sin z_n^0\pi)\cos z_n^0x+
(\beta+\gamma\cos z_n^0\pi)\sin z_n^0x.
$$
Нетрудно видеть, что $y_n^0(x)$ почти нормированы, т.е. существуют такие
постоянные $M_1$ и $M_2$, что $M_1\le\|y_n^0\|_{L_2}\le M_2$.

В силу теоремы~\ref{tm:2.5} имеем $|z_n-z_n^0|<M\Upsilon(z_n^2)$.
Воспользовавшись леммой~\ref{lem:2.4}, из представления~\eqref{(y)} легко
получить оценки $|y_n(x)-y_n^0(x)|\le M\Upsilon(z_n^2)$. Очевидно, такое
же неравенство сохраняется после нормировки функций $\{y_n\}$ и
$\{y_n^0\}$. Утверждение теоремы теперь следует из леммы~\ref{lem:2.3}.
При невыполнения условия усиленной регулярности рассуждения сохраняются,
но вместо величин $\Upsilon(z_n^2)$ в оценках будут участвовать
величины $\sqrt{\Upsilon(z_n^2)}$. Теорема доказана.
\end{proof}
\begin{Note}
Нужно оговорить, что при невыполнении условия усиленной регулярности
бесконечно много собственных значений могут оказаться двукратными и им
может отвечать пара собственных функций или одна собственная и одна
присоединенная функции. В этом случае выбор собственных и присоединенных
функций нужно производить специальным образом, чтобы утверждение
теоремы~\ref{tm:2.6} сохранялось. Такой выбор возможен, но
подробное доказательство из-за его громоздкости здесь опущено. В
действительности, при наличии кратных собственных значений у
невозмущенного оператора более важной является не задача об
оценках разности собственных или присоединенных функций
возмущенного и невозмущенного операторов, а задачи об оценке
норм разности проекторов Рисса $P_n$ и $P_n^0$ на
соответствующие двумерные подпространства. Такие оценки
проводились в работе~\cite{Sh}. В рассматриваемой	 задаче
метод~\cite{Sh} позволяет без труда получить следующий результат:
последовательность $\left\{\|P_n-P_n^0\|\right\}$ принадлежит
пространству $l_2^\theta$.
\end{Note}

\textbf{2.4. Базисность Рисса.}
В заключение этого параграфа докажем теорему о базисности Рисса собственных и
присоединенных функций. Напомним, что система \(\{y_n\}\) в гильбертовом
пространстве называется базисом Рисса, если существует ограниченный и
ограниченно обратимый оператор \(A\) такой, что система \(\{Ay_n\}\) является
ортонормированной и полной. За определением понятия базиса Рисса со скобками
(или из подпространств) мы отсылаем читателя к~\cite[Гл.~6]{GK}.
\begin{Theorem}\label{tm:2.9}
Пусть \(q(x)\in W_2^{-1}\), а оператор \(L\) определен согласно
теореме~\ref{tm:1.4}. Если соответствующие краевые условия усиленно
регулярны, то система собственных и присоединенных функций оператора
\(L\) образует базис Рисса. В случае регулярных краевых условий эта
система является базисом Рисса из подпространств, причем в
подпространства нужно объединять лишь собственные функции, отвечающие
собственным значениям \(\lambda_n\) и \(\lambda_{n+1}\), для которых
\(|\sqrt{\lambda_n}-\sqrt{\lambda_{n+1}}|\to 0\).
\end{Theorem}
\begin{proof}
Сначала предположим, что краевые условия являются усиленно регулярными. Тогда
известно (см.~\cite{DS}), что система собственных функций
\(\{y_n^0(x)\}_{n=1}^\infty\) оператора \(L_0\) образует базис Рисса в
\(L_2\) (для простоты предполагаем, что присоединенные функции
отсутствуют).  Согласно теореме~\ref{tm:2.6} собственные функции
\(\{y_n(x)\}_{n=1}^{\infty}\) оператора \(L\) образуют асимптотически
квадратично близкую систему, т.~е. для любого \(\varepsilon>0\) найдется
число \(N\) такое, что
\[
    \sum\limits_{n=N}^{\infty}\|y_n(x)-y_n^0(x)\|^2<\varepsilon.
\]
В силу теоремы Бари (см.~\cite[Гл.~6]{GK}) система
\(E=\{y_n^0\}_{n=1}^{N-1}\cup\{y_n\}_{n=N}^{\infty}\) также образует базис
Рисса в пространстве \(L_2\). Занумеруем все собственные функции оператора
\(L\), не вошедшие в систему \(\{y_n\}_{n=N}^{\infty}\), следующим образом:
\(\{y_n\}_{n=\omega}^{N-1}\), где \(-\infty\le\omega<N\). Система всех
собственных функций оператора \(L\) минимальна (никакая функция не принадлежит
замкнутой линейной оболочке остальных), а дефект системы
\(\{y_n\}_{n=N}^{\infty}\) равен \(N-1\) (поскольку система \(E\) есть базис
Рисса). Следовательно, \(\omega\ge 1\). Предположим, что \(\omega>1\). Тогда
замыкание линейной оболочки системы собственных функций оператора \(L\)
образует подпространство \(\mathfrak N\) некоторой коразмерности \(\ge 1\).
Подпространство \(\mathfrak N^{\perp}\) инвариантно относительно \(L^*\),
конечномерно и, по построению, на этом подпространстве оператор \(L^*\) не
имеет собственных значений. Это противоречие, поэтому \(\omega=1\). Но тогда
система \(\{y_n\}_{n=1}^{\infty}\) есть базис Рисса, так как она минимальна.

В случае сближающихся собственных значений доказательство меняется следующим
образом. Известно (см.~\cite{Sh}), что в этом случае система собственных
функций \(\{y_n^0\}_{n=1}^{\infty}\) есть базис Рисса из подпространств,
причем в подпространства нужно объединять не более двух функций.
Сравнивая функции Грина операторов \(L_0\) и \(L\) так же, как в
теореме~\ref{tm:2.6}, получим (см. подробнее~\cite{Sh}), что для любого
\(\varepsilon>0\) существует число \(N\) такое, что
\[
    \sum\limits_{k=N}^{\infty}\|P_k-P_k^0\|^2<\varepsilon,
\]
где \(P_k\) и \(P_k^0\) --- проекторы на двумерные подпространства, отвечающие
сближающимся собственным значениям. Теперь доказательство можно завершить,
используя такие же аргументы, как в случае несближающихся собственных
значений. Теорема доказана.
\end{proof}
При доказательстве последней теоремы мы одновременно установили взаимно
однозначное соответствие между первыми собственными значениями
$\{\la_n\}_{n=1}^{N-1}$ и $\{\la^0_n\}_{n=1}^{N-1}$ операторов $L$ и
$L_0$. Это полностью восполняет недостающее звено в доказательстве
теоремы~\ref{tm:2.5}.

\section*{{\bf \S 3.\ \ Вторые члены в асимптотиках для собственных
значений и собственных функций}}\refstepcounter{section}

Цель этого параграфа --- выделить вторые члены в асимптотике собственных
значений и собственных функций оператора \(L\).

Мы ограничимся рассмотрением краевых условий Дирихле
\[
y(0)=y(\pi)=0,
\]
хотя применяемый метод несложно модифицировать на случай краевых условий
типа Штурма
\begin{equation*}
\left\{
\begin{array}{ll}
\cos\alpha y(0)+\sin\alpha y^{[1]}(0)=0,\\
\cos\beta y(\pi)+\sin\beta y^{[1]}(\pi)=0,
\end{array}\right.
\end{equation*}
где $\alpha$, $\beta\in[0,\pi)$. Для нераспадающихся краевых условий
общего вида формулы получаются слишком громоздкими.

Оценка остатков в асимптотических формулах для собственных значений и 
собственных функций будет проводится в зависимости от класса потенциала 
$q(x)$, т.е. в зависимости от пространства, которому принадлежит функция 
$u(x)$. Для изучения случая $u(x)\in W_2^\theta[0,\pi]$ мы специально 
ограничиваемся значениями $0\le\theta<1/2$.
Конечно, похожие результаты можно получить и при $\theta/ge1/2$, но мы 
видим, что они не оптимальны. Для получения точных результатов при 
$\theta/ge1/2$ нужны новые средства и авторы планируют написать об этом в 
другой работе.

\textbf{3.1. Обозначения.}
Итак, в этом праграфе, через $L$ обозначаем оператор, порождаемый
дифференциальным выражением \eqref{de1} и краевыми условиями Дирихле на
отрезке $[0,\pi]$. Через $\Lal$, $0<\al\le1$, обозначаем функции класса
Липшица порядка $\alpha$ на отрезке $[0,\pi]$, т.~е. те функции, чей
модуль непрерывности
$$
\omega(\delta;f):=\sup\limits_{|x_1-x_2|<\delta}|f(x_1)-f(x_2)|,\quad
\text{где}\ x_1,\,x_2\in[0,\pi]
$$
допускает оценку $\omega(\delta;f)\le C\delta^\al,$
с постоянной $C$, не зависящей от $\delta$.
Через $\Lal^p$, $0<\al\le1$, $p\ge1$  обозначаем функции, чей
интегральный модуль непрерывности
$$
\omega_p(\delta;f):=\sup\limits_{0<t\le\delta}
\left[\intl_0^\pi|f(x+t)-f(x)|^pdx\right]^{1/p}
$$
допускает оценку $\omega_p(\delta;f)\le C\delta^\al$,
с постоянной $C$, не зависящей от $\delta$ ( в определении
предполагается, что функция $f(x)$ продолжена периодически за пределы
отрезка $[0,\pi]$).  Через $V$ обозначим класс функций ограниченной
вариации на $[0,\pi]$.  Потенциал $q(x)$ не предполагается вещественным,
т.~е. оператор $L$ не обязательно самосопряжен.

Как прежде, через $\{\la_n\}_{n=1}^\infty$ обозначаем собственные
значения
оператора $L$. Сохраняем также все обозначения, введенные в п.2.1.
Из них
наиболее часто будем использовать функции $\upsilon(c,x,\la)$ и
$\Upsilon(c,\la)$. Положим
\begin{equation*}
\begin{array}{ccc}
\mu_n:=-\frac1\pi \upsilon(0,\pi,n^2)=
-\frac1{\pi}\intl_0^\pi u(t)\sin(2nt)dt
+\frac1{2\pi n}\intl_0^\pi u^2(t)\cos(2nt)dt-\\
-\frac2{\pi}\intl_0^\pi\intl_0^t u(t)u(s)\cos(2nt)\sin(2ns)dsdt
-\frac1{2\pi n}\intl_0^\pi u^2(t)dt=-\frac12b_{2n}+\frac1{4n}A_{2n}-
\w_{2n}-\frac1{2\pi n}U(\pi),
\end{array}
\end{equation*}
где
$$
\begin{array}{cc}
b_n:=\frac2{\pi}\intl_0^\pi u(t)\sin(nt)dt,\quad
A_n:=\frac2{\pi}\intl_0^\pi u^2(t)\cos(nt)dt,\\
\w_n:=\frac2{\pi}\intl_0^\pi\intl_0^t u(t)u(s)\cos(2nt)\sin(2ns)dsdt,\quad
U(\pi)=\intl_0^\pi u^2(t)dt.
\end{array}
$$
Обозначим также
$$
a_n:=\frac2{\pi}\intl_0^\pi u(t)\cos(nt)dt,\quad r_n=w_{2n}-\frac1{2\pi
n}U(\pi).
$$

\textbf{3.2. Асимптотика собственных значений.}
Наша ближайшая цель --- показать, что числа $\mu_n$ играют роль
второго члена в асимптотических формулах для последовательности
$\{\la_n\}$.  Точнее, в этом пункте мы проведем оценку разности
$\sqrt{\la_n}-\mu_n$. Конечно, не все слагаемые в выражении для $\mu_n$
имеют одинаковую силу. Позже мы выясним, что главную роль (т.е. роль
второго члена асимптотики) играет только слагаемое $-b_n/2$.
\begin{Theorem}\label{tm:3.1}
Справедливы равенства
\begin{equation}\label{lan}
\la^{1/2}_n=n+\mu_n+\rho_n, \quad\text{где}\
|\rho_n|\le M\Upsilon^2(\la_n),
\end{equation}
а постоянная $M$ зависит только от функции $u(x)$. В частности, если
$u(x)\in W_2^\theta[0,\pi]$ при некотором $0\le\theta<1/2$, то
$\{\rho_n\}_{n=1}^\infty\in l_1^{2\theta}$.
\end{Theorem}
\begin{proof}
Пусть функция $\Psi(x,\la)$ -- решение
уравнения~\eqref{main2} с начальными условиями $\Psi(0,\la)=0$,
$\Psi^{[1]}(0,\la)=1$. Тогда уравнение $\Psi(\pi,\la)=0$
определяет
собственные значения оператора $L$. Выражая решения через функции
Прюфера,
получаем $\Psi(x,\la)=r(1,x,\la)\sin\theta(0,x,\la)$.
Следовательно,
уравнение для собственных значений имеет вид
$\sin\theta(0,\pi,\la)=0$.
Воспользовавшись леммой~\eqref{lem:2.1}, перепишем это уравнение в
виде
\begin{equation}\label{3.1}
\la_n^{1/2}+\frac1\pi\upsilon(0,\pi,\la_n)+\rho(\la_n)=n,
\end{equation}
где $|\rho(\la_n)|<M\Upsilon^2(\la_n)$. Поскольку
$\mu_n=\frac1\pi\upsilon(0,\pi,n^2)$, то доказательство теоремы
сводится к доказательству оценки $|\upsilon(0,\pi,\la_n)-
\upsilon(0,\pi,n^2)|<M\Upsilon^2(\la_n)$.

Обозначим $\nu_n:=\la_n^{1/2}-n$ и заметим, что из~\eqref{3.1}
следует оценка
$|\nu_n|< M\Upsilon(\la_n)$. Имеем
$$
\begin{array}{cc}
\upsilon(0,\pi,\la_n)-
\upsilon(0,\pi,n^2)=\intl_0^\pi\left(\sin2\la_n^{1/2}t-
\sin2nt\right)dt+\left(\la^{-1/2}_n-n^{-1}\right)
\left(U(\pi)+B(0,\pi,\la_n)\right)+\\+n^{-1}
\left(B(0,\pi,\la_n)-B(0,\pi,n^2)\right)+2(\w(0,\pi,\la_n)-
\w(0,\pi,n^2)).
\end{array}
$$
Обозначим слагаемые в правой части последнего равенства в порядке их
очередности через $I_1$, $I_2$, $I_3$, $I_4$. Для оценки интеграла $I_1$
воспользуемся формулой
$$
\sin(2\la_n^{1/2}t)-\sin(2nt)=(2\nu_nt)\cos(2\la^{1/2}_nt)+
O\left(\nu_n^2\right),
$$
которая при $0\le t\le\pi$ следует из разложения функции $\sin x$
в точке
$x=2\la_n^{1/2}t$ в ряд Тейлора. Следовательно,
$$
\begin{array}{ccc}
I_1=\left|2\nu_n\intl_0^\pi tu(t)\cos\left(2\la_n^{1/2}t\right)dt+
O\left(\nu_n^2\right)\right|=\\=\left|2\nu_n\pi\intl_0^\pi u(t)
\cos\left(2\la_n^{1/2}t\right)dt-2\nu_n\intl_0^\pi\intl_0^t
u(s)\cos\left(2\la_n^{1/2}s\right)dsdt\right|
+O\left(\nu_n^2\right)\le\\\le
M|\nu_n|\left|\Upsilon(\la_n)\right|<
M^2\left|\Upsilon(\la_n)\right|^2.
\end{array}
$$
Так как $\la_n^{-1/2}-n^{-1}=O(n^{-2})$ и $n^{-
2}<M\Upsilon^2(\la_n)$, то
$|I_2|<M\Upsilon^2(\la_n)$. Аналогично, поскольку
$\cos\left(2\la_n^{1/2}t\right)-\cos(2nt)=O(\nu_n)$ и
$|\nu_nn^{-1}|<M\Upsilon^2(\la_n)$, то
$|I_3|<M\Upsilon^2(\la_n)$. Наконец, для оценки $I_4$ воспользуемся
равенством
$$
\begin{array}{cc}
\cos\left(2\la_n^{1/2}t\right)\sin\left(2\la_n^{1/2}s\right)-
\cos\left(2nt\right)\sin(2ns)=\cos\left(2\la_n^{1/2}t\right)
\left(\sin\left(2\la^{1/2}_ns\right)-\sin\left(2ns\right)\right)+\\
+\left(\cos\left(2\la^{1/2}_nt\right)-\cos\left(2nt\right)\right)\sin(2ns)=
2\nu_ns\cos\left(2\la_n^{1/2}t\right)\cos\left(2\la_n^{1/2}s\right)-
2\nu_nt\sin\left(2\la_n^{1/2}t\right)\sin\left(2\la^{1/2}_ns\right)+
O\left(\nu_n^2\right).
\end{array}
$$
Из этого равенства следует, что достаточно получить оценку двойного
интеграла
$$
\left|\intl_0^\pi
u(t)\cos\left(2\la_n^{1/2}t\right)
\intl_0^tsu(s)\cos\left(2\la_n^{1/2}s\right)dsdt
\right|<M\Upsilon(\la_n),
$$
и такого же интеграла, где косинусы заменены на синусы. Но эта оценка
следует из определения функции $\Upsilon(\la)$ после изменения порядка
интегрирования.

Итак, асимптотика~\eqref{lan} установлена. Для завершения доказательства
теоремы, остается заметить, что в случае, когда $u(x)\in
W_2^\theta[0,\pi]$, в силу леммы~\ref{lem:2.3}
$\{\Upsilon^2(\la_n\}_{n=1}^\infty\in l_1^{2\theta}$, поскольку
$\{\sqrt{\la_n}\}_{n=1}^\infty$ --- несгущающаяся последовательность.
Теорема доказана.
\end{proof}

Полезно иметь информацию о характере убывания последовательности
$\{\Upsilon^2(\la_n\}_{n=1}^\infty$
не только для функций $u(x)\in W_2^\theta$, но и для функций
других классов. Здесь приведем результат для функций ограниченной
вариации и функций, удовлетворяющих условию Липшица в среднем.

\begin{Proposition}\label{prp:3.2}
Если $u(x)\in V$, то
\begin{equation}
|\rho_n|\le Mn^{-2}.
\label{3.a}
\end{equation}
Если $\omega_1(\delta)$ --- интегральный модуль непрерывности функции
$u(x)$ в пространстве $L_1$, то
\begin{equation}
|\rho_n|\le M\omega_1^2(n^{-1}).
\label{3.b}
\end{equation}
В частности, если $u(x)\in\Lal^1[0,\pi]$, то
\begin{equation}
|\rho_n|\le Mn^{-2\alpha}.
\label{3.с}
\end{equation}
Здесь $\{\rho_n\}_{n=1}^\infty$ --- последовательность остатков в
асимптотической формуле~\eqref{lan}, а $M=const$.
\end{Proposition}
\begin{proof}
Согласно теореме~\ref{tm:3.1} достаточно оценить числа
$\Upsilon^2(\la_n)$. Для функции $u(x)\in V[0,\pi]$ интегрированием по
частям получаем
\begin{equation}
\left|\intl_0^x u(t)\cos(zt)dt\right|+\left|\intl_0^x u(t)\sin(zt)dt
\right|=O\left(\dfrac1{|z|}\right)
\label{3.d}
\end{equation}
при $z\to\infty$ внутри полосы $|\Im z|<\al$. Тогда из определения
функции $\Upsilon(\la)$ следует
$|\Upsilon(\la_n)|<M|\la_n|^{-1}<M_1n^{-1}$, что влечет~\eqref{3.a}.

Известно~\cite[гл. 2.4]{Z}, что при $z\in\R^+$ левая часть~\eqref{3.d}
мажорируется величиной $M\omega_1(z^{-1})$. Кроме того
(см.~\cite[гл. 2.3]{Z}), $\delta<M\omega_1(\delta)$ при $\delta\to0$.
В случае вещественной функции $u(x)$ собственные значения тоже
вещественны, откуда получаем
$|\Upsilon(\la_n)|<M\omega_1(\la_n)<M_1\omega_1(n^{-1})$, что
влечет~\eqref{3.b}.

Для доказательства~\eqref{3.b} в случае комплексной $u(x)$ необходимо
модифицировать прием, используемый в~\cite{Z}. Для оценки первого
слагаемого в правой части~\eqref{3.d} воспользуемся формулой
$$
\intl_0^x u(t)\cos(zt)dt=-\intl_{-\pi/\nu}^{x-\pi/\nu}
u\left(t'+\frac{\pi}\nu\right)\cos\left(zt'-\frac{\pi
i\sigma}\nu\right)dt',\quad z=\nu+i\sigma,\ t'=t-\frac{\pi}\nu.
$$
Так как $\cos(zt)-\cos(zt-\pi
i\sigma/\nu)=O(\nu^{-1})=O(|z|^{-1})$ при $z\in\Pi_\al$, то
\begin{equation}\label{3.e}
\left|\intl_0^x u(t)\cos(zt)dt\right|
\le M\omega_1\left(|z|^{-1}\right),\quad z\in\Pi_\al.
\end{equation}
Второй интеграл в правой части~\eqref{3.d} оценивается аналогично.
Тем самым, оценка~\eqref{3.b} получена и в комплексном случае. Предложение
доказано.
\end{proof}

\textbf{3.3. Упрощенная формула для второго члена асимптотики.}
Числа $\mu_n$, определяющие второй член асимптотики собственных значений
в формуле~\eqref{lan}, представлены четырьмя слагаемыми. В
действительности, только первое слагаемое играет доминирующую роль, а
сумму остальных трех слагаемых можно оценить "почти" также, как
остаток в теореме~\ref{tm:3.1}. Тем самым, формулу для собственных
значений можно записать в виде
\begin{equation}\label{3.2.1}
\sqrt{\la_n}=n-\frac12b_{2n}+s_n,\quad
s_n=A_{2n}-r_n+\rho_n,
\end{equation}
а оценка чисел $\rho_n$ уже проведена в теореме~\ref{tm:3.1}. Наша
ближайшая цель --- провести оценку $A_n$ и $r_n$.
\begin{Proposition}\label{prp:3.3}
Если $u(x)\in V[0,\pi]$, то справедлива асимптотика~\eqref{3.2.1}, где
$$
s_n=O(n^{-2}).
$$
\end{Proposition}
\begin{proof}
Если $u(x)$ имеет ограниченную вариацию, то $u^2(x)$ обладает тем
же свойством. Поэтому $A_n=O(n^{-2})$. Числа $\w_{2n}$ определяются
двойными интегралами. Интегрируя внутренний интеграл по частям, получаем
\begin{equation}\label{3.2.2}
\begin{array}{cc}
r_n=\frac1{2\pi n}\intl_0^\pi u^2(t)dt-
\frac1{\pi n}\intl_0^\pi u^2(t)\cos^2(2nt)dt+\\
+\frac{u(0)}{\pi n}\intl_0^\pi u(t)\cos(2nt)dt+
\frac1{\pi n}\intl_0^\pi u(t)\cos(2nt)\intl_0^t\cos(2ns)du(s)dt.
\end{array}
\end{equation}
Используя формулу $\cos^2(2nt)=(1-\cos(4nt))/2$, получаем, что первые
три слагаемые в правой части равенства~\eqref{3.2.2} дают вклад
$O(n^{-2})$. Последнее слагаемое после интегрирования по частям
принимает вид
$$
\frac1{2\pi n^2}\intl_0^\pi\sin(2nt)d(u(t)\omega_n(t)),\quad
\omega_n(t)=\intl_0^t\cos(ns)du(s).
$$
Очевидно, полная вариация функции $\omega_n(t)$ не зависит от $n$ и не
превосходит полной вариации функции $u(t)$. Но тогда полная вариация
функции $u(t)\omega_n(t)$ не зависит от $n$, а потому последнее слагаемое
есть $O(n^{-2})$. Предложение доказано.
\end{proof}
\begin{Note}
Для функций $u(x)\in V[0,\pi]$ оператор $L$ был определен ранее
другими методами и детально изучен. В частности,
предложение~\ref{prp:3.3} было доказано в~\cite{Zhi}, а в недавней
работе~\cite{ВС} для потенциалов такого класса выписана более точная
асимптотическая формула.
\end{Note}
\begin{Proposition}\label{prp:3.4}
Для непрерывной функции $u(x)$ справедлива асимптотика~\eqref{3.2.1}, где
$$
s_n\le M\omega^2(n^{-1}),
$$
а $\omega(\delta)$ --- модуль непрерывности функции $u(x)$. В частности,
если $u(x)\in\Lal$, $0<\al\le1$, то
$$
s_n=O(n^{-2\al}).
$$
\end{Proposition}
\begin{proof}
Заметим, что модуль непрерывности функции $u^2(x)$ оценивается велечиной
$M\omega(\delta)$. Кроме того, всегда $\delta^{-1}<M\omega(\delta)$
и $\omega_1(\delta)<\pi\omega(\delta)$ (см.~\cite[гл. 2.3]{Z}). Тогда
из оценки~\eqref{3.e} получаем
$A_n<Mn^{-1}\omega(n^{-1})<M^2\omega^2(n^{-1})$, а из
предложения~\ref{prp:3.3} $|\rho_n|<M\omega^2(n^{-1})$. Поэтому остается
оценить числа $r_n$.

Заметим, что формула~\eqref{3.2.2} остается справедливой не только для
$u(x)\in V[0,\pi]$, но и для непрерывной функции $u(x)$, если интеграл
понимать в смысле Римана--Стилтьеса. Это следует из известного факта
(см.~\cite[гл. 8.6]{Nat}): {\it существование одного из интегралов в
формуле интегрирования по частям влечет существование второго}. Очевидно,
сумма первых трех интегралов в правой части~\eqref{3.2.2} оценивается
величиной $Mn^{-1}\omega_1(n^{-1})<M^2\omega^2(n^{-1})$. Поэтому
предложение будет доказано, как только удастся доказать неравенство
\begin{equation}\label{3.2.3}
\left|\intl_0^\pi u(t)\cos(2nt)\intl_0^t\cos(2ns)du(s)dt\right|<
M\omega(n^{-1}).
\end{equation}
Вновь обозначим $\omega_n(t):=\intl_0^t\cos(2ns)du(s)$ и докажем, что
\begin{equation}\label{3.2.4}
\left|\omega_n\left(t+\frac\pi{2n}\right)-\omega_n(t)\right|\le
M\omega(n^{-1}),\quad t\in[0,\pi].
\end{equation}
Разложив $\cos(2ns)$ в
точке $t$ в ряд Тейлора, получим \begin{equation}\label{3.2.5}
\left|\intl_t^{t+\frac\pi{2n}}\cos(2ns)du(s)\right|\le
\suml_{k=0}^\infty\frac1{k!}(2n)^k\left|
\intl_t^{t+\frac\pi{2n}}(s-t)^kdu(s)\right|.
\end{equation}
Для $k=0$ имеем
$$
\left|\intl_t^{t+\frac\pi{2n}}du(s)\right|=
\left|u\left(t+\frac\pi{2n}\right)-u(t)\right|
\le M\omega(n^{-1}).
$$
При $k\ge1$ получаем
$$
\begin{array}{cc}
\left|\intl_t^{t+\frac\pi{2n}}(s-t)^kdu(s)\right|=
\left|u\left(t+\frac\pi{2n}\right)
\left(\frac\pi{2n}\right)^k-\intl_t^{t+\frac\pi{2n}}u(s)d(s-t)^k\right|=\\
=\left|\intl_t^{t+\frac\pi{2n}}\left(u\left(t+\frac\pi{2n}\right)-u(s)
\right)d(s-t)^k\right|\le
M\omega\left(\frac\pi{2n}\right)\left(\frac\pi{2n}\right)^k<
M\omega(n^{-1})\left(\frac\pi{2n}\right)^k.
\end{array}
$$
Эта оценка вместе с~\eqref{3.2.5} влечет~\eqref{3.2.4}.

Теперь применим стандартный прием. Воспользовавшись равенством
$\cos 2n(t+\pi/(2n))=-\cos2nt$, легко получить
$$
2\intl_0^\pi u(t)\omega_n(t)\cos(2nt)dt=
\intl_0^\pi\left[u\left(t+\frac\pi{2n}\right)\omega_n
\left(t+\frac\pi{2n}\right)-
u(t)\omega_n(t)\right]\cos(2nt)dt+O\left(n^{-1}\right).
$$
Функции под интегралом в правой части этого равенства оцениваются
величиной $M\omega(n^{-1})$. Это влечет неравенство~\eqref{3.2.3}.
Предложение доказано.
\end{proof}
Для функций $u(x)\notin C[0,\pi]$ интегрирование по частям во внутреннем
интеграле в представлении числа $\w_n$ проводить нельзя. Поэтому
нужны другие приемы для оценки чисел $r_n$. Наши оценки будут основаны
на следующей лемме.
\begin{Lemm}\label{lem:3.5}
Справедливы равенства
\begin{equation}\label{3.2.6}
r_n=-\frac1{4n}\suml_{k=1,\,k\ne2n}^\infty
\frac{k^2a_k^2}{(2n-k)(2n+k)}+\gamma_n,
\end{equation}
где
\begin{equation}\label{3.2.7}
\gamma_n=\frac3{16n}a_{2n}^2+\frac{a_0a_{4n}}{8n}+\frac{a_{2n}}{\pi}
\intl_0^\pi(\pi-s)u(s)\sin(2ns)ds.
\end{equation}
\end{Lemm}
\begin{proof}
В двойной интеграл, выражающий число $\w_n$, подставим ряд
\begin{equation*}
u(t)=\suml_{k=0}^\infty a_k\cos(kt)
\end{equation*}
и поменяем пределы интегрирования. В дальнейших преобразованиях
интегрируем по частям и используем формулы произведений синусов. В
результате получаем
$$
\begin{array}{cccc}
\w_n=\intl_0^\pi u(s)\sin(2ns)\intl_s^\pi\suml_{k=0}^\infty\left(
a_k\cos(kt) \cos(2nt)\right)dtds=\\=\frac1{2\pi n}\intl_0^\pi
u(s)\suml_{k=0}^\infty
a_k\left[-\cos(kt)\sin^2(2ns)+k\sin(2ns)\intl_s^\pi\sin(kt)\sin(2nt)
\right]dtds=\\=
\frac1{4\pi n}\intl_0^\pi u(s)\suml_{k=0}^\infty a_k\left[\cos(kt)
(-1+\cos(4ns))\right]+2na_{2n}(\pi-s)\sin(2ns)ds+\\
+\frac1{4\pi n}\intl_0^\pi u(s)\left[\suml_{k=1}^\infty\frac{ka_k}{2n+k}
\sin(2ns)\sin(2n+k)s-\suml_{k=1,\,k\ne2n}^\infty\frac{ka_k}{2n-k}
\sin(2ns)\sin(2n-k)s\right]ds.
\end{array}
$$
Вновь используя формулы для произведений косинусов и синусов, получаем
\begin{equation}\label{3.2.8-}
\begin{array}{cc}
2\w_n+\frac1{2\pi n}U(\pi)=\frac1{8n}\suml_{k=0}^\infty a_k(a_{4n+k}+
a_{4n-k})+\frac{a_{2n}}{\pi}\intl_0^\pi(\pi-s)u(s)\sin(2ns)ds+\\
+\frac{a_{2n}}{16n}(a_{2n}-a_{6n})+\frac1{8n}{\suml_{k=1}^\infty}' k\left(
\frac{a_k^2}{2n+k}-\frac{a_k^2}{2n-k}\right)+\frac1{8n}{\suml_{k=1}^\infty}'
k\left(\frac{a_ka_{4n-k}}{2n-k}-\frac{a_ka_{4n+k}}{2n+k}\right),
\end{array}
\end{equation}
где штрих в знаке суммы означает, что в сумме исключается индекс $k=2n$.
Воспользовавшись равенством
$$
-\frac1{8n}(a_{4n+k}+a_{4n-k})+\frac14\left(\frac{a_{4n+k}}{2n+k}+
\frac{a_{4n-k}}{2n-k}\right)=\frac{k}{8n}\left(\frac{a_{4n-k}}{2n-k}-
\frac{a_{4n+k}}{2n+k}\right),
$$
преобразуем последнюю сумму в правой части~\eqref{3.2.8} к виду
$$
-\frac1{8n}{\suml_{k=1}^\infty}' a_k(a_{4n+k}+a_{4n-k})+
\frac14{\suml_{k=1}^\infty}'\left(\frac{a_ka_{4n+k}}{2n+k}+
\frac{a_ka_{4n-k}}{2n-k}\right).
$$
Здесь первая сумма с точностью до слагаемого
$\frac{a_0a_{4n}}{4n}+\frac{a_{2n}(a_{2n}+a_{6n})}{8n}$ сокращается
с первой суммой в правой части~\eqref{3.2.8-}, а вторая сумма равна
$$
\frac14\left({\suml_{k=1}^{4n}}'\frac{a_ka_{4n-k}}{2n-k}\right)+
\frac14\left({\suml_{k=1}^\infty}'\frac{a_ka_{4n+k}}{2n+k}-
\suml_{k=4n+1}^{\infty}\frac{a_ka_{4n-k}}{k-2n}\right).
$$
После сокращения одинаковых слагаемых в первой и второй сумме останутся
только  числа $-\frac{a_0a_{4n}}{8n}$ и $-\frac{a_{2n}a_{6n}}{16n}$,
соответственно. Учитывая в формуле~\eqref{3.2.8-} проведенные вычисления,
приходим к формулам~\eqref{3.2.6} и~\eqref{3.2.7}. Лемма доказана.
\end{proof}

Теперь оценим числа $s_n$ в асимптотике~\eqref{3.2.1} в зависимости от
принадлежности функции $u(x)$ различным пространствам.
\begin{Proposition}\label{prp:3.6}
Если $u(x)\in\Lal^1[0,\pi]$, $0<\al\le1$, то числа $s_n$ в
формуле~\eqref{3.2.1} допускают оценку
$$
s_n<M(n^{-2\al}\ln n+n^{-1}\omega_2(n^{-1})).
$$
В частности, для $u(x)\in\Lal^2[0,\pi]$ справедлива оценка
\begin{equation}\label{3.2.9}
s_n=O(n^{-2\al}\ln n).
\end{equation}
\end{Proposition}
\begin{proof}
Из неравенства Коши--Буняковского следуют оценки
$$
\omega_1(u^2,\delta)<M\omega_2(u,\delta)\quad
\omega_1(u,\delta)<\sqrt{\pi}\omega_2(u,\delta).
$$
Тогда из
неравенства~\eqref{3.e} имеем $A_n<Mn^{-1}\omega_2(n^{-1})$. В частности,
для $u(x)\in\Lal^2\subset\Lal^1$ выполнены неравенства $A_n<Mn^{-1-\al}$.
Поэтому в силу леммы~\ref{lem:3.5} достаточно доказать для функций
$u(x)\in\Lal^1$ оценку
$$
\frac1n{\suml_{k=1}^\infty}'\frac{k^2|a_k|^2}{|2n-k|(2n+k)}=
O(n^{-2\al}\ln n).
$$
Для $u(x)\in\Lal^1$ имеем $a_k=O(k^{-\al})$, поэтому выписанный ряд
мажорируется суммой
$$
\frac1{2n}{\suml_{k=1}^{3n}}'\frac{k^{1-2\al)}}{|2n-k|}+
\frac3n\suml_{k=3n+1}^\infty|a_k|^2.
$$
Очевидно, первая сумма
есть $O(n^{-2\al}\ln n)$. При $\al\le1/2$ второй ряд
есть $O(n^{-1})=O(n^{-2\al})$, а при $\al>1/2$ после замены $|a_k|$ на
$O(k^{-\al})$ получаем такую же оценку. Предложение доказано.
\end{proof}
\begin{Proposition}\label{prp:3.7}
Пусть $u(x)\in L_2[0,\pi]$. Тогда последовательность $\{s_n\}_{n=1}^\infty$
в асимптотической формуле~\eqref{3.2.1} обладает свойством
$$
\{s_n\}_{n=1}^\infty\in l_p,\quad\{s_n\ln^{-p}(n+1)\}_{n=1}^\infty\in l_1
$$
при любом $p>1$.
\end{Proposition}
\begin{proof}
В силу теоремы~\ref{tm:3.1} $s_n=A_n+r_n+\rho_n$, где $A_n=o(n^{-1})$,
$\{\rho_n\}_{n=1}^\infty\in l_1$. Поэтому в силу леммы~\ref{lem:3.5}
достаточно доказать сформулированное утверждение для чисел
\begin{equation}\label{3.2.8}
\xi_n=\frac1n{\suml_{k=1}^\infty}'\frac{ka_k^2}{(2n-k)(2n+k)}=
{\suml_{k=1}^\infty}'a_k^2\left(\frac1{2n+k}+\frac1{2n-k}-\frac1n\right).
\end{equation}
Фиксируем любое число $p>1$ и рассмотрим произвольную последовательность
$\{c_n\}_{n=1}^\infty\in l_q$, $1/p+1/q=1$. В силу неравенства Гельдера
получаем
$$
\begin{array}{cc}
\suml_{n=1}^\infty|s_nc_n|\le\suml_{n=1}^\infty|c_n|{\suml_{k=1}^\infty}'
|a_k|^2\left(\frac1{|2n-k|}+\frac1{|2n+k|}+\frac1n\right)<\\
<3{\suml_{k=1}^\infty}'|a_k|^2\left(\suml_{n=1}^\infty|c_n|^q\right)^{1/q}
\left(\suml_{n=1}^\infty\frac1{|2n-k|^p}\right)^{1/p}.
\end{array}
$$
Так как ${\sum_n}'|2n-k|^{-p}<2/(p-1)$ (не зависит от $k$), то ряд
$\sum|s_nc_n|$ сходится для любой последовательности
$\{c_n\}_{n=1}^\infty\in l_q$. Следовательно,
$\{s_n\}_{n=1}^\infty\in l_p$. Далее, ряд
$$
{\suml_{n=1}^\infty}'\left(\frac1{|2n-k|}+\frac1{2n+k}+\frac1n\right)
\frac1{\ln^p(n+1)}
$$
сходится и его сумма мажорируется константой, не зависящей от $k$.
Следовательно, сходится ряд $\sum|s_n|\ln^{-p}(n+1)$.
Предложение доказано.
\end{proof}
Установим теперь аналог предложения~\ref{prp:3.7} для функций
$u(x)\in W_2^\theta$.
\begin{Proposition}\label{prp:3.8}
Если $u(x)\in W_2^\theta$, $0<\theta<1/2$, то последовательность
$\{s_n\}_{n=1}^\infty$ в асимптотической формуле~\eqref{3.2.1} обладает
свойствами
$$
\{s_n\ln^{-1}(n+1)\}_{n=1}^\infty\in l_1^{2\theta},\quad
\{s_n\}_{n=1}^\infty\in l_p^{2\theta}\quad\text{при любом}\  p>1.
$$
\end{Proposition}
\begin{proof}
В силу теоремы~\ref{tm:3.1} достаточно доказать это утверждение, когда
вместо $\{s_n\}$ фигурирует последовательность $\{A_n\}$ и $\{r_n\}$.

Рассмотрим сначала последовательность $\{A_n\}$. Пусть
$\|\cdot\|_{\al,p}$ --- норма в пространстве $W_p^\al$. Воспользуемся
мультипликативной оценкой (см.\cite{BS}):
$$
\|fg\|_{\al,p}\le M\|f\|_{\theta,2}\|g\|_{\theta,2},
$$
которая справедлива при
$$
0\le\al\le\theta,\quad1\le p\le\frac1{1-2\theta+\al},\quad
0\le\theta<\frac12.
$$
Положим $f=g=u(x)$. В случае $\theta\in[0,1/4]$ выберем $\al=0$,
$p=1/(1-2\theta)$. Согласно теореме Хаусдорфа--Юнга (ее можно применять,
т.к. $p\le2$, см.~\cite{Ba}) последовательность
$$
\{nA_n\}=\left\{\frac1{2\pi}\intl_0^\pi u^2(t)\cos(2nt)dt\right\}
$$
лежит в пространстве $l_q$, $1/p+1/q=1$. В силу неравенства Гельдера
$$
\suml_{n=1}^\infty\frac{n^{2\theta}|A_n|}{\ln^\beta(n+1)}\le
\left(\suml_{n=1}^\infty|nA_n|^q\right)^{1/q}\left(\suml_{n=1}^\infty
\left(\frac{n^{2\theta-1}}{\ln^\beta(n+1)}\right)^p\right)^{1/p}.
$$
Так как $p=1/(1-2\theta)$, то ряд в правой части сходится при любом
$\beta$, если $\beta/(1-2\theta)>1$, в частности при $\beta=1$.
Следовательно, $\{A_n\ln^{-1}(n+1)\}_{n=1}^\infty\in l_1^{2\theta}$.

В случае $\theta\in[1/4,1/2)$ воспользуемся той же мультипликативной
оценкой при $\al=2\theta-1/2$, $p=2$. Тогда $\al<1/2$,
$u^2(x)\in W_2^\al$, следовательно, $\{nA_n\}\in l_2^\al$. Из
неравенства Коши-Буняковского получаем
$$
\suml_{n=1}^\infty\frac{n^{2\theta}|A_n|}{\ln^\beta(n+1)}\le
\left(\suml_{n=1}^\infty n^{2(1+\al)}|A_n|^2\right)^{1/2}\left(
\suml_{n=1}^\infty\frac{n^{2(2\theta-1-\al)}}{\ln^{2\beta}n}
\right)^{1/2}.
$$
Все ряды сходятся при $\beta>1/2$, так как $2(2\theta-1-\al)=-1$.
Поэтому $\{A_n\ln^{-1}(n+1)\}_{n=1}^\infty\in l_1^{2\theta}$.

Теперь рассмотрим последовательность $\{r_n\}$. Оценка
последовательности $\{\gamma_n\}$, определенной в~\eqref{3.2.7}
тривиальна. Поэтому достаточно рассмотреть последовательность $\xi_n$,
определенную в~\eqref{3.2.8}. Воспользуемся оценками
$$
\frac{k^2}{n|2n-k|(2n+k)}<\left\{\begin{array}{cc}
kn^{-1}|2n-k|^{-1},\ \ &\text{если}\ 1\le n\le k,\\
k^2n^{-3},\ \ &\text{если}\ k<n<\infty.\end{array}\right.
$$
Тогда
$$
\suml_{k=1}^\infty|\xi_n|\frac{n^{2\theta}}{\ln(n+1)}\le
\suml_{k=1}^\infty|a_k|^2\left(\suml_{n=1}^k\frac{kn^{2\theta-1}}
{|2n-k|\ln(n+1)}+\suml_{n=k+1}^\infty k^2n^{2\theta-3}\right).
$$
Легко проверить, что обе суммы в скобках оцениваются величиной
$Mk^{2\theta}$, где постоянная $M$ не зависит от $k$. Но условие
$u(x)\in W_2^\theta$ влечет сходимость ряда $\sum|a_k|^2k^{2\theta}$,
следовательно, ряд в левой части также сходится. Таким образом,
$\{\xi_n\ln^{-1}(n+1)\}_{n=1}^\infty\in l_1^{2\theta}$.

Доказательство свойства $\{s_n\}_{n=1}^\infty\in l_p^{2\theta}$ при
$p>1$ проводится с помощью приема, примененного в
предложении~\ref{prp:3.7}. При этом оценки сходимости получающихся рядов
проще, нежели в первом рассмотренном случае и потому здесь опускаются.
Предложение доказано.
\end{proof}

\textbf{3.4. Случай функции из класса
Колмогорова--Силеверстова--Плеснера.}
Естественен вопрос: можно ли в случае $u(x)\in L_2[0,\pi]$ усилить
утверждение, сформулированное в предложении~\ref{prp:3.7} и доказать, что
$\{s_n\}\in l_1$? Позже мы покажем, что ответ на этот вопрос отрицателен,
а тогда естественно найти класс функций (конечно, включающий пространства
$W_2^\theta$ при $\theta>0$), для которых свойство
$\{s_n\}_{n=1}^\infty\in l_1$ выполнено. Наша цель --- показать, что это
свойство выполнено для функций класса Колмогорова--Силеверстова--Плеснера
(см.~\cite[гл. 8]{Ba}). Этот класс состоит из функций $u(x)$, для которых
\begin{equation}
\suml_{k=1}^\infty|a_k|^2\ln k<\infty,\quad a_k=\frac2{\pi}\intl_0^\pi
u(t)\cos(kt)dt.
\label{3.3.1}
\end{equation}
Для удобства сформулируем некоторые определения и результаты из теории
тригонометрических рядов, которые будем использовать.

Известно~\cite[гл. 8]{Ba}, что для произвольной суммируемой функции $f(x)$
почти при всех $x\in[0,\pi]$ существует предел
$$
\wf(x):=-\frac1{\pi}\lim_{\eps\to0}\intl_\eps^\pi\frac{f(x+t)-f(x-t)}
{2\tg(t/2)}dt.
$$
Функция $\wf(x)$ называется сопряженной к $f(x)$. Нижеследующие результаты
можно найти в монографиях~\cite{Ba},~\cite{Z}.\hfill\break
{\bf Теорема A.} (М. Рисс). \quad{\itshape Если
$f(x)\in L_p$  при $p>1$, то $\wf(x)\in L_p$.}\hfill\break
{\bf Теорема B.} (Ф.Рисс и М. Рисс). \quad{\itshape Если функция $f(x)$
абсолютно непрерывна, а сопряженная к ее производной
$\widetilde{f'}(x)$ суммируема, то ряд Фурье функции $f(x)$ абсолютно
сходится.}\hfill\break
{\bf Теорема C} (Лузин--Данжуа)\quad {\it Если тригонометрический ряд
$$
\frac{a_0}2+\suml_{n=1}^\infty a_n\sin(nx)+b_n\cos(nx)
$$
сходится абсолютно на множестве ненулевой меры, то ряд $\sum(|a_n|+|b_n|)$
сходится.}\hfill\break
{\bf Теорема D.} (Зигмунд)\quad {\it Если функция $f(x)$ измерима и
$|f(x)|\ln^+|f(x)|$ суммируема, то сопряженная к $f(x)$ также суммируема.}
\hfill\break
{\bf Теорема E.} (Плеснер)\quad {\it Условие~\eqref{3.3.1} эквивалентно
условию
\begin{equation}\label{3.3.2}
\intl_0^\pi\intl_0^\pi\frac{|u(x+t)-u(x-t)|^2}tdtdx<\infty
\end{equation}
(здесь предполагается, что $u(x)$ периодически продолжена с отрезка
$[0,\pi]$ на $\R$).}

Следующая лемма, возможно, известна, но мы не смогли найти ее в литературе
и потому приведем ее с доказательством.
\begin{Lemm}\label{lem:3.9}
Если $u(x)$ удовлетворяет условию~\eqref{3.3.2}, то сопряженная к функции
$u^2(x)$ суммируема.
\end{Lemm}
\begin{proof}
Поскольку разность $\tg(t/2)-1/(2t)$ непрерывна, то достаточно доказать, что
функция
$$
\widehat{u^2}(x):=\intl_0^\pi\frac{u^2(x+t)-u^2(x-t)}tdt
$$
суммируема. Очевидно, условие~\eqref{3.3.2} эквивалентно условию
$$
[u(x+t)-u(x-t)]t^{-1/2}\in L_2(Q),\quad Q=[0,\pi]\times[0,\pi].
$$
Это
свойство влечет
$$
[u(x+t)+u(x-t)-2u(x)]t^{-1/2}\in L_2(Q).
$$
Тогда
$$[u(x+t)-u(x-t)][u(x+t)+u(x-t)-2u(x)]t^{-1}\in L_1(Q).
$$
В силу теоремы Фубини почти для всех $x\in [0,\pi]$ существует интеграл
$$
\intl_0^\pi\frac1t(u^2(x+t)-u^2(x-t))dt-2u(x)
\intl_0^\pi\frac1t(u(x+t)-u(x-t))dt=\widehat{u^2}(x)-2u(x)\widehat{u}(x).
$$
Так как $u(x)\in L_2$, то в силу теоремы A получаем $\widetilde{u}(x)\in
L_2$, а потому $\widehat{u}(x)\in L_2$. Следовательно,
$\widehat{u^2}(x)\in L_1$, а тогда $\widetilde{u^2}(x)\in L_1$. Лемма
доказана.
\end{proof}
\begin{Proposition}\label{prp:3.10}
Если $u(x)$
подчинена условию~\eqref{3.3.1}, то последовательность $\{s_n\}$ в
асимптотической формуле~\eqref{3.2.1} принадлежит $l_1$.
\end{Proposition}
\begin{proof}
Достаточно показать, что $\{A_n\}$, $\{\xi_n\}\in l_1$, где $\xi_n$
определены в~\eqref{3.2.8}. Интегрируя по частям, получаем
$$
A_n=\frac1{2\pi n}\intl_0^\pi u^2(t)\cos(2nt)dt=\frac1{\pi}\intl_0^\pi
f(t)\sin(2nt)dt,\quad f(x):=\intl_0^x u^2(t)dt-x\intl_0^\pi
u^2(t)dt.
$$
Функция $f(x)$ абсолютно непрерывна, а согласно лемме~\ref{lem:3.9}
сопряженная к $f'(x)$ суммируема. В силу теоремы B ряд Фурье
функции $f(x)$ сходится абсолютно. В силу теоремы C тогда
$\{A_n\}\in l_1$.

Докажем, что $\{\xi_n\}_{n=1}^\infty\in l_1$. Имеем
$$
\suml_{n=1}^\infty|\xi_n|\le\suml_{k=1}^\infty|a_k|^2
{\suml_{n=1}^\infty}'\left|\frac1{2n+k}+\frac1{2n-k}-\frac1n\right|.
$$
Очевидно, внутренняя сумма мажорируется величиной $M\ln k$. Тогда
условие~\eqref{3.3.1} влечет сходимость ряда. Предложение доказано.
\end{proof}

\textbf{3.5. Обсуждение точности результатов. Итоговая теорема.}
Возникает вопрос: можно ли доказать, что $\{s_n\}_{n=1}^\infty\in l_1$
без требования условия~\eqref{3.3.1}? Мы уже упоминали, что
ответ на этот вопрос отрицателен. Более того, мы покажем, что
существует непрерывная вещественная функция $u(x)$, для
которой $\{s_n\}_{n=1}^\infty\notin l_1$. На первый взгляд это кажется
странным. Но противоречий нет: пространство $C[0,\pi]$ геометрически
"плохо"  вложено в $L_2[0,\pi]$, в частности, опеатор
вложения не компактен. Если же рассмотреть пространство
функций, в котором норма определена рядом~\eqref{3.3.1}, то
это пространство компактно вложено в $L_2[0,\pi]$. Поэтому не
все непрерывные функции удовлетворяют
условию~\eqref{3.3.1}.

\begin{Proposition}\label{prp:3.11}
Существует вещественная непрерывная функция $u(x)$, для которой
последовательность \(\{s_n\}_{n=1}^\infty\) из~\eqref{lan} не
принадлежит пространству \(l_1\).
\end{Proposition}
\begin{proof}
Положим
$u(x)=\sum a_n\cos(nx)$, где
$$
a_n=\left\{\begin{array}{ll}
\frac1{\sqrt{\ln(n/2)}}=p^{-2}&,\quad\text{если } n=2\cdot2^{p^4},\
p=1,2,\dots,\\ 0&,\quad\text{иначе.}
\end{array}\right.
$$
Очевидно, $u(x)$ --- вещественная функция, а ее непрерывность следует
из равномерной сходимости ряда
$$
\left|\suml_{n=1}^\infty
a_n\cos(nx)\right|\le\suml_{n=1}^\infty|a_n|=
\suml_{p=1}^\infty\frac1{p^2}<\infty.
$$
Очевидно также, что условие~\eqref{3.3.1} для этой функции не
выполнено.

Докажем, что $\{n^{-1}A_n\}_{n=1}^\infty\in l_1$, а
$\{r_n\}_{n=1}^\infty\notin l_1$, откуда будет следовать
$\{s_n\}_{n=1}^\infty\notin l_1$. Заметим, что функция
$\widetilde{u^2}(x)$ (сопряженная к $u^2(x)$) суммируема. Это следует из
непрерывности $u^2(x)$ и сформулированной в предыдущем пункте теоремы D.
Но тогда, повторяя доказательство предложения~\ref{prp:3.10}, получаем
$\{n^{-1}A_n\}_{n=1}^\infty\in l_1$.

Пусть числа $\xi_n$ определены равенствами~\eqref{3.2.8}. Докажем, что
$\{\xi_n\}_{n=1}^\infty\notin l_1$. Тогда из леммы~\ref{lem:3.5} будет
следовать $\{r_n\}_{n=1}^\infty\notin l_1$.

Положим $n_p=2^{p^4}$ и рассмотрим множества индексов
$\Gamma_p=\{n_p+j\}_{j=1}^{n_p}$, $p=1,2,\dots$. Эти множества не
пересекаются при различных $p$, а потому
\begin{equation}\label{3.3.11}
\suml_{n=1}^\infty|\xi_n|\ge
\suml_{p=1}^\infty\suml_{n\in\Gamma_p}|\xi_n|.
\end{equation}
Теперь из определения чисел $\xi_n$ при $n\in\Gamma_p$ получаем
$$
\xi_n=\frac{(2n_p)^2}{(2n-n_p)(2n+n_p)n\ln n_p}+\left(\suml_{s=1}^{p-1}-
\suml_{s=p+1}^\infty\right)\frac{2n_s^2}{|2n-2n_s|(2n+2n_s)n\ln n_s}>
$$
$$
>\frac1{(n-n_p)\ln n_p}-\suml_{s=p+1}^\infty\frac2{ns^4}.
$$
Здесь мы воспользовались явным видом для коэффициентов $a_k$ и
неравенствами $n+n_s\le2n_s$, если $s\ge p$ и $|2n-2n_s|>n_s$, если
$s\ge p+1$ (в обоих неравенствах $n\in\Gamma_p$). Теперь заметим, что ряд
$$
\suml_{p=1}^\infty\suml_{n\in\Gamma_p}\suml_{s=p+1}^\infty\frac2{ns^4}<
\suml_{p=1}^\infty\frac2{3p^3}\intl_{n_p}^{2n_p}\frac1xdx
$$
сходится, а ряд
$$
\suml_{p=1}^\infty\suml_{n\in\Gamma_p}\frac1{(n-n_p)\ln n_p}=
\suml_{p=1}^\infty\frac1{\ln n_p}\suml_{j=1}^{n_p}\frac1j
$$
расходится. Но тогда расходится ряд~\eqref{3.3.11}. Предложение доказано.
\end{proof}

В предложении~\ref{prp:3.7} мы доказали, что
$\{s_n\ln^{-p}(n+1)\}_{n=1}^\infty\in l_1$ при любом $p>1$, если $u(x)\in
L_2$. Построенный пример показывает, что множитель $\ln^{-p}(n+1)$ убрать
нельзя. Возникает естественный вопрос: можно ли положить $p=1$ (как это
было сделано в предложении~\ref{prp:3.8})? Ответ на этот вопрос нам найти
не удалось, но полезно заметить, что существуют функции $u(x)\in L_2$,
для которых $\{A_nn^{-1}\ln^{-1}(n+1)\}_{n=1}^\infty\notin l_1$. В
качестве примера рассмотрим функцию $$ f(x)=\suml_{n=1}^\infty\frac{\cos
nx}{\ln\ln(n+1)}.  $$ Последовательность $a_n=\frac1{\ln\ln(n+1)}$
монотонна и выпукла, а потому указанный ряд сходится при всех $x>0$ к
неотрицательной суммируемой функции $f(x)$ (см.~\cite[гл 5.1]{Z}). Но
тогда $u(x)=\sqrt{f(x)}\in L_2$ и $A_n=\frac1{\ln\ln(n+1)}$.
Следовательно, $\{A_nn^{-1}\ln^{-1}(n+1)\}_{n=1}^\infty\notin l_1$.

Полезно резюмировать полученные в этом параграфе основные результаты.
\begin{Theorem}\label{tm:3.12}
Собственные значения $\la_n$ оператора $L$ с краевыми условиями Дирихле
подчинены асимптотике~\eqref{3.2.1}, где последовательность $\{s_n\}$
в зависимости от функции $u(x)$ обладает следующими свойствами:
\begin{itemize}
\item[i)] $s_n=O(n^{-2})$, если $u(x)\in V$;
\item[ii)] $s_n=O(n^{-2\al})$, если $u(x)\in\Lal$;
\item[iii)] $s_n=O(n^{-2\al}\ln n)$, если $u(x)\in\Lal^2$;
\item[iv)] $s_n\in l_p$ и $\{s_n\ln^{-p}(n+1)\}_{n=1}^\infty\in l_1$
при любом $p>1$, если $u(x)\in L_2$;
\item[v)] $s_n\in l_p^{2\theta}$ при любом $p>1$ и
$\{s_n\ln^{-1}(n+1)\}_{n=1}^\infty\in l_1^{2\theta}$, если $u(x)\in
W_2^\theta$, $0<\theta<1/2$;
\item[vi)] $s_n\in l_1$, если $u(x)$ подчинена условию
Колмогорова--Силеверстова--Плесснера~\eqref{3.3.1}.
\end{itemize}
\end{Theorem}

\textbf{3.6. Асимптотика собственных функций.}
Здесь мы не будем подробно останавливаться на деталях и приведем итоговый
результат.
\begin{Theorem}~\label{tm:3.13}
Собственные функции оператора $L$ имеют асимптотику
\begin{equation}\label{3.6.1}
\begin{array}{cc}
y_n(x)=\sin(nx)+\sin nx\left(-\intl_0^x u(t)\cos(2nt)dt+\intl_0^\pi
u(t)\left(1-\frac{t}{\pi}\right)\cos(2nt)dt\right)+\\+\cos
nx\left(-\frac{x}\pi\intl_0^\pi u(t)\sin(2nt)dt+
\intl_0^x u^2(t)\sin(2nt)dt\right)+\psi_n(x),\quad n=1,2,\dots,
\end{array}
\end{equation}
а последовательность функций $\psi_n(x)$ в зависимости от функции $u(x)$
обладает следующими свойствами:\begin{itemize}
\item[i)]
$|\psi_n(x)|\le Mn^{-2}$, если $u(x)\in V$;
\item[ii)]
$|\psi_n(x)|\le Mn^{-2\al}$, если $u(x)\in \Lal$;
\item[iii)]
$|\psi_n(x)|\le Mn^{-2\al}\ln n$, если $u(x)\in \Lal^2$;
\item[iv)]
$\suml_{n=1}^\infty|\psi_n(x)|^p+\suml_{n=1}^\infty|\psi_n(x)|n^{1-p}<
M(p)$ при любом $p>1$, если $u(x)\in L_2$;
\item[v)]
$\suml_{n=1}^\infty|\psi_n(x)|^{p}n^{2\theta}+
\suml_{n=1}^\infty|\psi_n(x)|n^{1-p}<\infty$ при любом $p>1$,
если $u(x)\in W_2^\theta$, $0<\theta<1/2$;
\item[vi)]
$\suml_{n=1}^\infty|\psi_n(x)|<M$, если $u(x)$ подчинена
условию~\eqref{3.3.1}.
\end{itemize}
Во всех оценках постоянная $M$ не зависит от $x$.
\end{Theorem}
\begin{proof}
Подставим асимптотические формулы~\eqref{asth} и~\eqref{asr} для функций
$\theta(x,\la)$ и $r(x,\la)$ в формулу~\eqref{pr}. Получим, что решение
$y(x,\la)$ уравнения~\eqref{main2}, удовлетворяющее условию $y(0,\la)=0$,
имеет асимптотику
$$
y(x,\la)=\sin(\la^{1/2}x)\left(1-b(0,x,\la)-\frac12\la^{-1/2}A(0,x,\la)
\right)+\cos(\la^{1/2}x)\upsilon(0,x,\la)+\rho(x,\la),
$$
где $|\rho(x,\la)|<M\Upsilon^2(\la)$. Учитывая равенство
$\la_n^{1/2}=n-\frac12b_{2n}+s_n$, после несложных преобразований
получим представление~\eqref{3.6.1}. При этом в остаток $\psi_n(x)$
войдут несколько функций, из которых наиболее сложно оценивается
функция
\begin{equation}\label{3.22}
2\w(0,x,n^2)+\frac1{2n}U(x)=\intl_0^x\intl_0^tu(t)u(s)
\cos(2nt)\sin(2ns)dsdt+\frac1{2n}\intl_0^xu^2(t)dt.
\end{equation}
Ее оценка в случае принадлежности $u(x)$ пространствам $V$, $\Lal$,
$\Lal^p$ проводится также, как при $x=\pi$. В случае $u(x)\in
W_2^\theta$, $0\le\theta<1/2$, нужно воспользоваться тем, что оператор
умножения на характеристическую функцию $\chi_{[0,x]}$ непрерывен в
$W_2^\theta$ (см., например,~\cite[теорема 17.11]{BIN}). Заметим, что 
выражение~\eqref{3.22} совпадает с выражением 
$2\w(0,\pi,n^2)+\frac1{2n}U(\pi)$ сосчитанным для функции 
$\chi_{[0,x]}u(x)$. Но при $x=\pi$ все оценки были получены ранее.  
Функции $u(x)$, подчиненные условию~\eqref{3.3.1} входят в $W_2^\theta$ 
при любом $\theta>0$, поэтому для них проведенные рассуждения остаются в 
силе. Этим завершается доказательство теоремы.  
\end{proof}

\section*{{\bf \S 4.\ \ Операторы Штурма\,--\,Лиувилля с потенциалами высокой
сингулярности}}\refstepcounter{section}

В \S\,1.5 было показано, что попытки определить оператор Штурма--Лиувилля,
порожденный дифференциальным выражением~\eqref{de1}
с потенциалом, не принадлежащим классу распределений \(W_2^{-1}\),
сопряжены с трудностями. Во всяком случае, такой оператор не является
корреткно определенным с точки зрения предельного перехода. Основной вывод
состоит в следующем --- для построения оператора Штурма--Лиувилля с потенциалом
$q(x)$ --- распределением высокого порядка ($\int q(x)dx\notin L_2$)
недостаточно той информации, которая заключена в определении обобщенной
функции $q(x)$. Здесь мы покажем, как можно обойти эти трудности, если привлечь
дополнительную информацию о потенциале.
Нас будут интересовать потенциала вида
\begin{equation}
q(x)=x^\alpha
\label{1.1}
\end{equation}
на отрезке $[-1,1]$. Отметим, что $q(x)$ регулярен при
$\Re\,\alpha>-1$ и является распределением первого порядка при
$\Re\,\alpha>-3/2$ (если только $q(x)\ne \frac1{|x|}$, см. \S\,1.6, пример 1).
Наша цель --- корректно определить оператор Штурма-Лиувилля с данным
потенциалом при других значенях $\alpha$.
Мы предложим два метода определения этого оператора приводящие к
одинаковому результату.
Первый из них основан на аналитическом продолжении операторного
семейства, определенного при $\Re\,\alpha>-3/2$, в область
$-2<\Re\,\alpha\le-3/2$. Второй базируется на последовательной
регуляризации дифференциального выражения \eqref{de1}. Отметим, что
используемые здесь методы имеют аналогию с методами, которые были
использованы Бернштейном для определения обобщенных функций $x^\al$ (см.
\cite{GSh}). Рассмотренные здесь потенциалы являются модельными и оба
метода допускают обобщения на более широкий класс потенциалов.

{\bfseries 4.1 Построение оператора методом аналитического
продолжения.}
Здесь мы построим оператор Штурма-Лиувилля $L(\al)$ на отрезке $[-1,1]$
с потенциалом $q(x)=-|x|^\al$ (другие потенциалы типа (\ref{1.1}) будут
кратко рассмотрены ниже). Зафиксируем краевые условия Дирихле
$$
y(-1)=y(1)=0.
$$
Эти краевые условия мы выбираем для определенности, во избежаение
технических сложностей, не связанных с сутью задачи.

Мы начнем с поиска фундаментальной системы решений (ФСР) однородного
уравнения
\begin{equation}
-y''-|x|^\alpha y=0
\label{2.0.I}
\end{equation}
на промежутке $[0,1]$. Обозначим
\begin{equation}
\Pi_\alpha:=\{\alpha\ |\ \Re\,\alpha>-2,\quad \alpha\ne -2+\frac1n,\
n=1,2,3,\dots\}.
\label{2.1.I}
\end{equation}
Хорошо известно (см., например \cite{Фе}), что
если $\alpha\in \Pi_\alpha$, то ФСР имеет вид
$$
U_1(x)=m^{-1/m}\Gamma\Big(1-\dfrac1m\Big)\sqrt{x}J_{-1/m}\Big(\dfrac2m
x^{m/2}\Big),\quad
U_2(x)=m^{1/m}\Gamma\Big(1+\dfrac1m\Big)\sqrt{x}J_{1/m}\Big(\dfrac2m
x^{m/2}\Big),
$$
где $m=2+\alpha$, $\Re\, m>0$, $m\ne \frac1n, \ n=1,2,3,\dots$,
$\Gamma(\cdot)$ -- гамма функция, а $J_\nu(\cdot)$ -- функции Бесселя.
Аналогично можно выбрать ФСР на промежутке $[-1,0]$. При этом в силу
четности $q(x)$ можно взять
$$
U_1(-x)=U_1(x),\quad
U_2(-x)=-U_2(x).
$$
Выпишем асимптотические разложения функций $U_j(x)$ в нуле, $x>0$:
\begin{equation}
\begin{array}{ccc}
U_1(x)=1-\frac{x^m}{1!m(m-1)}+\frac{x^{2m}}{2!m^2(m-1)(2m-1)}-\dots,
\\ \phantom{.}\\
U_2(x)=x\big[
1-\frac{x^m}{1!m(m+1)}+\frac{x^{2m}}{2!m^2(m+1)(2m+1)}-\dots\big] .
\end{array}
\label{3.I}
\end{equation}
Обе функции $U_j(x)$, $j=1,2$ лежат в $L_2$. Для определения ФСР
уравнения (\ref{2.0.I}) на $[-1,1]$ необходимо 'склеить' ФСР на $[-1,0]$ и
$[0,1]$.  Далее мы будем считать функцию $U_1(x)$ четной, а функцию $U_2(x)$
нечетной на $[-1,1]$. Это фиксирует способ 'склейки' и, как будет
доказано позднее, полностью определяет искомый оператор.

Итак, приступим к построению оператора. Используя метод вариации
постоянных, можно записать общее решение неоднородного уравнения
\begin{equation}
-y''-|x|^\alpha y=f(x),\qquad f(x)\in L_2[-1,1],
\label{4.0.I}
\end{equation}
в виде
$$
y(x)=\left( C_1+\intl_{-1}^x U_2(t)f(t)dt\right) U_1(x)-
     \left( C_2+\intl_{-1}^x U_1(t)f(t)dt\right) U_2(x),
$$
или
\begin{equation}
y(x)=AU_1(x)+BU_2(x)+\intl_0^xU_2(t)f(t)dtU_1(x)-\intl_0^xU_1(t)f(t)dtU_2(x).
\label{5.I}
\end{equation}
Тогда
$$
y(x)=AU_1(x)+BU_2(x)+z(x),\quad\text{где}\ z(x)\in W_2^2[-1,1]\ \text{и}\
z(0)=z'(0)=0.
$$
Поскольку это общий вид решения уравнения (\ref{4.0.I}), то можно
определить оператор $L(\al)$ равенствами
\begin{equation*}
\begin{array}{cc}
L(\al)y=l(z),\\
\Dom(L(\al))=
\left\{y(x)=AU_1(x)+BU_2(x)+z(x)\bigg|
\begin{array}{rr}
A,B\in\C;\, z(x)\in W_2^2[-1,1];\\z(0)=z'(0)=0;\, y(\pm 1)=0.
\end{array}\right\}.
\end{array}
\end{equation*}
Теперь нам потребуется теорема о существовании и единственности
решения задачи Коши
\begin{equation}
L(\al)y=\lambda y+f
\label{6.I}
\end{equation}
а также о голоморфной зависимости этого решения от параметров $\alpha$ и
$\lambda$.

\begin{Proposition}
Пусть $\lambda\in \C$, $|\lambda|<\varepsilon$, где $\eps$ достаточно
малое число. Зафиксируем начальные условия
\begin{equation}
y(\xi)=a,\quad y'(\xi)=b,
\label{7.I}
\end{equation}
где $a,b\in\C$; $\xi\in[-1,1]$, $\xi\ne 0$. Тогда при любом \(f\in L_2\)
существует единственное решение уравнения (\ref{6.I}) из
$\mathfrak D(L(\al))$, т.~е.
$$
y(x)=AU_1(x)+BU_2(x)+z(x),\quad\text{где}\ z(x)\in W_2^2[-1,1]\ \text{и}\
z(0)=z'(0)=0,
$$
причем $y(x)$ удовлетворяет условиям (\ref{7.I}). Кроме того, функция
$y(x,\lambda,\alpha)$ голоморфна по $\alpha$ в $\Pi_\alpha$ для
фиксированных $x\in[-1,1]$ и $\lambda$, $|\lambda|<\varepsilon$ и
голоморфна по $\lambda$ при $|\lambda|<\varepsilon$ для фиксированных
$x\in[-1,1]$ и $\alpha\in \Pi_\alpha$.
\end{Proposition}
\begin{proof}
Запишем решение уравнения (\ref{4.0.I}) с начальными условиями
$y(\xi)=y'(\xi)=0$
$$
\widetilde{y}(x)=\intl_\xi^xU_2(t)f(t)dtU_1(x)-\intl_\xi^xU_1(t)f(t)dtU_2(x).
$$
Так как вронскиан пары $\{U_1(x),U_2(x)\}$ равен $1$ (это следует из
асимптотических формул (\ref{3.I})), то можно подобрать константы
$A$ и $B$, так что
$y_0(x)=AU_1(x)+BU_2(x)+\widetilde{y}(x)$ удовлетворяет условиям
(\ref{7.I}). Достаточно взять
\begin{equation}
A=b\,U_1(\xi)-a\,U_2'(\xi),\quad
B=a\,U_1'(\xi)-b\,U_2(\xi).
\label{8.I}
\end{equation}
Рассмотрим отображения из $L_2$ в $\Dom(L(\al))$
$$
y_0=F(f),\quad
\widetilde{y}=Ef.
$$
Отображение $E$ линейно и оба этих отображения ограничены:
$$
\begin{array}{ccc}
\|y_0\|_{L_2}\le\|U_1\|_{L_2}\intl_{-1}^1|U_2(t)||f(t)|dt+
\|U_2\|_{L_2}\intl_{-1}^1|U_1(t)||f(t)|dt
+|A|\|U_1\|_{L_2}+|B|\|U_2\|_{L_2}\le
\\ \phantom{.}\\
\le 2\|U_1\|_{L_2}\|U_2\|_{L_2}\|f\|_{L_2}+|A|\|U_1\|_{L_2}+
|B|\|U_2\|_{L_2}\le C_1\|f\|_{L_2}+C_2,
\end{array}
$$
поскольку $a$, $b$ и $\xi$ фиксированы;
$$
\|\widetilde{y}\|_{L_2}\le C_1\|f\|_{L_2}.
$$
Искомое решение уравнения (\ref{6.I}) найдем методом последовательных
приближений $y_0=F(f)$, $y_1=F(\lambda y_0+f)$, $\dots$,
$y_n=F(\lambda y_{n-1}+f)$, $\dots$.
Тогда
\begin{equation}
y_n(x)=y_0(x)+\sum\limits_{k=1}^n\left( y_k(x)-y_{k-1}(x)\right).
\label{9.I}
\end{equation}
Функции $v_k(x)=y_k(x)-y_{k-1}(x)$ удовлетворяют уравнению
$$
L(\al)v_k=\lambda v_{k-1},
$$
причем $v_k(\xi)=v_k'(\xi)=0$, и мы получаем оценки
$$
\|v_k\|_{L_2}\le C_1|\lambda|\|v_{k-1}\|_{L_2}\le\dots
\le\left( C_1|\lambda|\right)^{k-1}\|v_1\|_{L_2}.
$$
Таким образом, при малых $|\lambda|$ ряд (\ref{9.I}) сходится в $L_2$ и
определяет искомое решение уравнения (\ref{6.I}) $y(x,\lambda,\alpha)$.

Докажем единственность построенного решения. Требуется показать, что
функция $y(x)$, такая что
$$
y(x)=AU_1(x)+BU_2(x)+z(x),\quad\text{где}\ z(x)\in W^2_2[-1,1],\
z(0)=z'(0)=0\ \text{и}
$$
$$
-z''-|x|^\alpha z=\lambda y,
$$
удовлетворяющая начальным условиям $y(\xi)=y'(\xi)=0$ тождественно равна
нулю. Пусть, например, $\xi<0$. Тогда в силу классической теоремы
единственности, $y(x)\equiv 0$ на $[-1,0]$, т.е. $z(x)=-AU_1(x)-BU_2(x)$.
Так как $z(0)=0$, то из асимптотических формул
(\ref{3.I}) следует, что $A=0$. Аналогично, $z'(0)=0$ влечет
равенство $B=0$. Итак, $y(x)=z(x)$ на $[-1,1]$, где
\begin{equation}
z(0)=z'(0)=0,\qquad z(x)\in W^2_2[-1,1],
\label{10.I}
\end{equation}
\begin{equation}
-z''-|x|^\alpha z=\lambda z.
\label{11.I}
\end{equation}
Осталось доказать, что $z(x)\equiv 0$ на $[0,1]$. Заметим, что из условий
(\ref{10.I}) следует оценка $|z(x)|\le Cx^{3/2}$. Рассмотрим функцию
$w(x)=z(x)x^{-3/2}$. Пусть $\max\limits_{0\le x
\le\delta}|w(x)|=M$, где $\delta$ достаточно мало. Перепишем
уравнение (\ref{11.I}) в интегральном виде
$$
z(x)=-\intl_0^x\intl_0^t(\lambda+s^\alpha)z(s)dsdt.
$$
Отсюда легко следует оценка
$$
|z(x)|\le \frac43M(|\lambda|+1)x^{7/2+\al}
\le\frac43M(\eps+1)x^{7/2+\al}<Mx^{3/2}
$$
при надлежащем выборе $\delta$,
что влечет $M=0$ и $z(x)\equiv 0$ на $[0,\delta]$. Из классической
теоремы единственности получим, что $z(x)\equiv 0$ на всем отрезке $[0,1]$.

Теперь докажем голоморфность $y(x,\lambda,\alpha)$ по $\alpha$. Заметим,
что
\begin{equation}
|\widetilde{y}(x,\lambda,\alpha)|=|Ef(x)|\le
2\sup\limits_{[-1,1]}\{|U_1(x,\alpha)||U_2(x,\alpha)|\}\|f\|_{L_2}
\le C_3\|f\|_{L_2}.
\label{12.I}
\end{equation}
Уменьшая, если надо $|\lambda|$, получим равномерную сходимость
ряда (\ref{9.I}). В силу голоморфности $U_j(x,\alpha)$, $j=1,2$, функции
$\widetilde{y}(x,\alpha)$ и $y_0(x,\alpha)$ также голоморфны. Из
(\ref{5.I}) и (\ref{8.I}) следует голоморфность $y_n(x,\alpha)$ а значит
и $y(x,\lambda,\alpha)$. Аналогично доказывается голоморфность
$y(x,\lambda,\alpha)$ по $\lambda$.
\end{proof}

Построение оператора $L(\al)$ неоднозначно. Однако наше определение
$L(\al)$ оправдывается тем, что построенное операторное семейство является
голоморфным по $\al$. Напомним (см.~\cite{Kato}), что \textit{ограниченное
семейство операторов $T(\alpha)$ в гильбертовом пространстве
называется голоморфным, если для любых
векторов $f$ и $g$ функция $\left( T(\alpha)f,g\right)$ голоморфна;
неограниченное семейство операторов $T(\alpha)$ называется голоморфным, если
существуют ограниченные голоморфные семейства $U(\alpha)$ и $V(\alpha)$
такие что
\[
    U(\alpha)T(\alpha)=V(\alpha).
\]}
\begin{Theorem}
Область определения
$\Dom(L(\al))$ плотна в $L_2$. Резольвентное множество
$\rho(L(\al))$ не пусто и $\forall\lambda\in\rho(L(\alpha_0))$ резольвента
$R_\lambda=(L(\al_0)-\lambda)^{-1}$ компактна. Семейство операторов
$R_\lambda(\alpha)$ ограниченно голоморфно
в некоторой окрестности точки $\alpha_0$ по параметру $\alpha$ (при
фиксированном $\lambda$). Сами операторы $L(\alpha)$ образуют неограниченное
самосопряженное голоморфное семейство при $\alpha\in\Pi_\alpha$,
где область $\Pi_\al$ определена в (\ref{2.1.I}).
Спектр $L(\alpha)$ дискретен, а собственные значения
$\lambda_k(\alpha)$ аналитичны по $\alpha$ и более того, не имеют
особенностей. Все операторы $L(\alpha)$ полуограничены снизу.
\end{Theorem}
\begin{proof}
Плотность линеала $\Dom(L(\al))$ в $L_2[-1,1]$ очевидна. Докажем, что
$\rho(L(\alpha_0))\ne\varnothing$. Рассмотрим $\varphi(x)$ -- решение
уравнения (\ref{6.I}) для $f(x)\equiv 0$ с начальными условиями
$\varphi(-1)=0$,
$\varphi'(-1)=1$. В силу голоморфной зависимости от $\lambda$, функция
$\varphi(1,\lambda,\alpha_0)$ при $|\lambda|<\varepsilon$ обращается в
ноль на дискретном множестве. Пусть $\varphi(1,\lambda_0,\alpha_0)\ne 0$.
Для произвольной функции $f(x)\in L_2$ найдем решение уравнения (\ref{6.I})
$\psi(x)$ с
теми же начальными условиями $\psi(-1)=0$, $\psi'(-1)=1$. Тогда
$y(x)=\psi(x)-\frac{\psi(1)}{\varphi(1)}\varphi(x)$ удовлетворяет
(\ref{6.I}) и $y(x)\in\Dom(L(\al))$. Значит
$\lambda_0\in\rho(L(\alpha_0))$. Компактность оператора $R_{\lambda_0}$
следует из того, что операто $L(\al)$ является двумерным возмущением
оператора $L_\oplus$ --- прямой суммы операторов
$-\frac{d^2}{dx^2}-|x|^{-2+1/n}$ на $[-1,0]$ и $[0,1]$ с условиями Дирихле
на концах.
В силу голоморфности
$\varphi(1,\lambda_0,\alpha)$ по $\alpha$ можно выбрать достаточно малую
окрестность ${\mathcal O}(\alpha_0)$ точки $\alpha_0$, в которой эта функция не
обращается в ноль. Тогда семейство ограниченных операторов
$R_{\lambda_0}(\alpha)$ определено в окрестности ${\mathcal O}(\alpha_0)$ и
является голоморфным в этой окрестности. Из равенства $$
L(\alpha)R_{\lambda_0}(\alpha)=I+\lambda_0R_{\lambda_0}(\alpha)
$$
следует голоморфность семейства $L(\alpha)$ в ${\mathcal O}(\alpha_0)$. В
силу произвольности $\alpha_0$ получаем, что семейство $L(\alpha)$ голоморфно в
$\Pi_\alpha$. При $\Re\,\alpha>-1$ семейство $L(\alpha)$ является
самосопряженным голоморфным семейством. В силу теоремы единственности, это
свойство сохранится и при $\alpha\in\Pi_\alpha$. Так как резольвента
рассматриваемых операторов компактна, то спектры рассматриваемых операторов
состоят из собственных значений \(\lambda_k(\alpha)\), причем эти собственные
значения голоморфно зависят от \(\alpha\in\Pi_{\alpha}\) (см.~\cite{Kato}).
Из этого факта и из полуограниченности $L(\alpha)$ при $\Re\,\alpha>-1$
следует полуограниченность при $\alpha\in\Pi_\alpha$. Наконец компактность
резольвенты при всех $\lambda$, а не только при $|\lambda|<\varepsilon$
следует из резольвентного тождества
$$
R_\lambda-R_\mu=(\lambda-\mu)R_\lambda R_\mu.
$$
Теорема полностью доказана.
\end{proof}
\textbf{4.2. Определение операторов в исключительных точках.}

Нам осталось определить оператор в точках $\alpha=-2+\frac1n$,
$n=1,2,3\dots$.

\begin{Proposition}
Пусть $\alpha\to\alpha_n$, где $\alpha_n=-2+\frac1n$ для некоторого
$n\in\N$.  Тогда операторы $L(\alpha)$ сходятся в сильном резольвентном
смысле к оператору $L_\oplus$ -- прямой сумме операторов
$-\frac{d^2}{dx^2}-|x|^{-2+1/n}$ на $[-1,0]$ и $[0,1]$ с условиями Дирихле
на концах.
\end{Proposition}
\begin{proof}
Пусть
$$
a(m)=U_1(1)=U_1(-1),\quad
B(m)=U_2(1)=U_2(-1).
$$
\begin{equation}
\begin{array}{c}
Rf=y(x)=\left\{\frac{b}{2a}\intl_{-1}^1U_1(t)f(t)dt-\frac12\intl_{-1}^1
U_2(t)f(t)dt+\intl_{-1}^xU_2(t)f(t)dt\right\}U_1(x)-\\
-\left\{\frac{a}{2b}\intl_{-1}^1U_2(t)f(t)dt-\frac12\intl_{-1}^1
U_2(t)f(t)dt+\intl_{-1}^xU_1(t)f(t)dt\right\}U_2(x),
\end{array}
\label{13.I}
\end{equation}
если $a\ne 0$ и $b\ne 0$. В силу того, что $J_n(2n)\ne 0\ \forall n\in\N$
(см. \cite{Ват}) можно выбрать такую окрестность точки $\alpha_n$, что
$a(m)$ и $b(m)$ не обращаются в ноль когда $\alpha=-2+m$ лежит в этой
окрестности.

Определим сначала \(y(x)=Rf(x)\) при $x>0$. Если $f(x)$ четна, то
\begin{equation}
y(x)=\left\{\frac{b}{a}\intl_0^1U_1(t)f(t)dt-\intl_x^1U_2(t)f(t)dt\right\}
U_1(x)-\intl_0^1U_1(t)f(t)dtU_2(x).
\label{14.I}
\end{equation}
Здесь мы воспользовались четностью $U_1(x)$ и нечетностью $U_2(x)$.
В силу равенства
$$
J_\nu(z)=J_{-\nu}(z)+Y_\nu(z)\sin\pi\nu,
$$
получим
\begin{equation}
\begin{array}{ccc}
\hskip-15pt
U_1(x)=\frac{\Gamma(1-1/m)}{\Gamma(1+1/m)}m^{-2/m}U_2(x)-
m^{-1/m}\Gamma\left(1-\frac1m\right)\sin\frac\pi{m}\sqrt{x}
Y_{1/m}\left(\frac2mc^{m/2}\right)=\hskip-10pt\\ \phantom{.}\\
=\frac{\Gamma(1-1/m)}{\Gamma(1+1/m)}m^{-2/m}U_2(x)+U_3(x),\qquad x>0.
\end{array}
\label{15.I}
\end{equation}

Подставляя (\ref{15.I}) в (\ref{14.I}), имеем
\begin{equation}
\begin{array}{c}
y(x)=\intl_0^xU_2(t)f(t)dtU_3(x)-\left\{ \frac{U_3(1)}{U_2(1)}\intl_0^1
U_2(t)f(t)dt- \right. \\ \left.
-\intl_0^1U_3(t)f(t)dt+\intl_0^xU_3(t)f(t)dt \right\}
U_2(x)+o(1)\quad
\text{при}\ m\to\frac1n.
\end{array}
\label{16.I}
\end{equation}

Таким образом, в пределе получим действие на функцию $f$ резольвенты $R_+$
оператора $-\frac{d^2}{dx^2}-|x|^{-2+1/n}$ на $[0,1]$ с условиями
Дирихле. Если же $f(x)$ нечетна, то
\begin{equation}
y(x)=\intl_0^xU_2(t)f(t)dtU_1(x)-\left\{\frac{a}{b}\intl_0^1U_2(t)f(t)dt-
\intl_x^1U_1(t)f(t)dt\right\}U_2(x).
\label{17.I}
\end{equation}
Подставляя (\ref{15.I}) в (\ref{17.I}) имеем
$$
\begin{array}{c}
y(x)=\intl_0^xU_2(t)f(t)dtU_3(x)-\left\{\frac{U_3(1)}{U_2(1)}\intl_0^1
U_2(t)f(t)dt
-\intl_x^1U_3(t)f(t)dt\right\}U_2(x)+o(1)
\end{array}
$$
при $m\to\frac1n$,
что совпадает с (\ref{16.I}). Итак, для произвольной $f(x)\in L_2$
получаем равенство
$$
Rf=R_+f+o(1),\quad m\to\frac1n,\quad x>0.
$$

Точно так же можно найти \((Rf)(x)\) при $x<0$. Разница лишь в том, что $R_+$
надо заменить на $R_-$ -- резольвенту оператора
$-\frac{d^2}{dx^2}-|x|^{-2+1/n}$ на $[-1,0]$ с условиями Дирихле.
\end{proof}

Таким образом, мы определили оператор $L(\alpha)=
-\frac{d^2}{dx^2}-|x|^\alpha$, $x\in[-1,1]$ с условиями Дирихле на концах
при всех $\alpha$, $\Re\,\alpha>-2$.

\textbf{4.3. Другие голоморфные семейства.}
Опишем операторы $L(\alpha)$ для других голоморфных семейств потенциалов типа
(\ref{1.1}).
\hfill\break
{\bf 1.}\quad $q(x)=-c|x|^\alpha$, где $c>0$.
Все рассуждения и полученные ранее утверждения полностью сохраняются в этом
случае. В качестве ФСР надо взять функции
$c^{\frac1{2m}-\frac14}U_1(\sqrt{c}x)$ и
$c^{-\frac1{2m}-\frac14}U_2(\sqrt{c}x)$. \hfill\break
{\bf 2.}\quad $q(x)=c|x|^\alpha$, где $c>0$.
В этом случае можно провести аналогичные рассуждения. Заметим, что при
\(c=1\) ФСР уравнения $$
-y''+x^\alpha y=0
$$
на отрезке $[0,1]$ имеет вид
\begin{equation}
\begin{array}{cc}
V_1(x)=i^{1/m}m^{-1/m}\Gamma\left(1-\dfrac1m\right)\sqrt{-x}
J_{-1/m}\left(\dfrac2m i(-x)^{m/2}\right),\\
V_2(x)=-i^{-1/m}m^{1/m}\Gamma\left(1+\dfrac1m\right)\sqrt{-x}
J_{1/m}\left(\dfrac2m i(-x)^{m/2}\right).
\end{array}
\label{19.I}
\end{equation}
Эти функции имеют следующие асимптотические разложения в нуле
$$
\begin{array}{ccc}
V_1(x)=1+\frac{(-x)^m}{1!m(m-1)}+\frac{(-x)^2m}{2!m^2(m-1)(2m-1)}+\dots,
\\\phantom{.}\\
V_2(x)=x\left[1+\frac{(-x)^m}{1!m(m+1)}+\frac{(-x)^2m}{2!m^2(m+1)(2m+1)}+
\dots\right].
\end{array}
$$
В силу четности потенциала, эти же функции будут образовывать ФСР
на отрезке $[-1,0]$. Остается склеить ФСР на отрезке $[-1,1]$
из функций $c^{\frac1{2m}-\frac14}V_1(\sqrt{c}x)$
и $c^{-\frac1{2m}-\frac14}V_2(\sqrt{c}x)$. Первую функцию мы опять
продолжим четным, а вторую -- нечетным образом.

{\bf 3.\quad $q(x)=c\sign x|x|^\alpha$.}
В этом случае ФСР уравнения
\begin{equation}
-y''-\sign x|x|^\alpha y=0
\label{18.I}
\end{equation}
на $[0,1]$ будет такой же, как для уравнения (\ref{2.0.I}). Однако на
$[-1,0]$ в качестве ФСР надо взять функции $V_1(x)$ и $V_2(x)$.
В качестве ФСР уравнения (\ref{18.I}) на $[-1,1]$ возьмем функции
$$
U_1(x)=\left\{\begin{array}{c}c^{\frac1{2m}-\frac14}U_1(x),\quad x>0
\\
c^{\frac1{2m}-\frac14}V_1(x),\quad x<0
\end{array}\right.,
\quad
U_2(x)=\left\{\begin{array}{c}c^{-\frac1{2m}-\frac14}U_2(x),\quad x>0
\\
c^{-\frac1{2m}-\frac14}V_2(x),\quad x<0
\end{array}\right..
$$
Получаемое семейство операторов $L(\alpha)$, $\alpha\in\Pi_\alpha$ также
будет голоморфно и все утверждения теоремы 4.2 верны.
Отличие данного
семейства от предыдущих состоит в том, что в точках $\alpha_n=-2+\frac1n$,
$n=1,3,5,\dots$ предельный оператор будет отличен от прямой суммы.
Это легко понять, если заметить, что соотношение (\ref{15.I}) будет иметь
место для всех $x\in[-1,1]$ а не только для $x>0$, т.к. в (\ref{19.I})
коэффициенты $i^{1/m}$ и $-i^{-1/m}$ совпадут. Подставляя (\ref{15.I})
непосредственно в (\ref{13.I}) получим в пределе оператор, отличный от
прямой суммы. Нетрудно показать, что при четных $n$ предельный оператор
по прежнему будет прямой суммой.

\textbf{4.4. Определение оператора методом последовательной
регуляризации.}
Для определенности рассмотрим операторы \(L(\alpha)=-y''+|x|^\al y\).
Напомним, что основным шагом в построении оператора при
$\al>-3/2$, т.е. для случая $q(x)=u'(x)$, где $u(x)\in L_2$, было введение
квазипроизводной
$$
y^{[1]}=y'-\frac{\sign x|x|^{\al+1}}{\al+1}y
$$
и переход от резольвентного уравнения
\begin{equation}
-y''+|x|^\al y=\la y+f
\label{00.I}
\end{equation}
к системе дифференциальных уравнений первого порядка
\begin{equation}
\left(\begin{array}{cc} y_1\\y_2\end{array}\right)'=
\begin{pmatrix}\frac{|x|^{\al+1}}{\al+1}\sign x&1\\\phantom{.}\\
-\la+\frac{|x|^{2\al+2}}{(\al+1)^2}&\quad-\frac{|x|^{\al+1}}{\al+1}\sign x
\end{pmatrix}\left(\begin{array}{cc} y_1\\y_2\end{array}\right)+
\left(\begin{array}{cc} 0\\f\end{array}\right)
\label{1.I}
\end{equation}
с помощью замены $y_1=y$, $y_2=y^{[1]}$ (см. \S 1).
При $\al>-3/2$ элементы
матрицы системы лежат в $L_1$, что позволяет доказать теорему
существования и единственности и далее определить оператор $L$.

Изложенный в \S 1 метод можно развить таким образом, что он будет
давать результат и при других значениях параметра $\al$. Возьмем
$$
\begin{array}{ccc}
y_1(x)=y(x),\\\phantom{.}\\
y_2(x)=y'(x)-\left(\frac{|x|^{\al+1}}{\al+1}\sign x-
\frac{|x|^{2\al+3}\sign x}{(\al+1)^2(2\al+3)}\right)y(x).
\end{array}
$$
Непосредственная проверка показывает, что при такой замене матрица системы
уравнений будет иметь вид
$$
\begin{pmatrix}
\frac{|x|^{\al+1}}{\al+1}\sign x-\frac{|x|^{2\al+3}\sign
x}{(\al+1)^2(2\al+3)}&1\\ \phantom{.}\\
-\la+2\frac{|x|^{3\al+4}}{(\al+1)^3(2\al+3)}-\frac{|x|^{4\al+6}}{(\al+1)^4
(2\al+3)^2}&\phantom{.....}-\frac{|x|^{\al+1}}{\al+1}\sign
x+\frac{|x|^{2\al+3}\sign x}{(\al+1)^2(2\al+3)}
\end{pmatrix}.
$$
Элементы этой матрицы лежат в $L_1$ при $\al>-\frac53$,
$\al\ne-\frac32$, $-1$. Для того чтобы продвинуться дальше и обработать
значения параметра $\al>-\frac74$, $\al\ne-\frac53$, $-\frac32$, $-1$,
необходимо сделать замену
$$
\begin{array}{ccc}
y_1(x)=y(x),\\\phantom{.}\\
y_2(x)=y'(x)-\left(\frac{|x|^{\al+1}}{\al+1}\sign x-
\frac{|x|^{2\al+3}\sign x}{(\al+1)^2(2\al+3)}-\frac{2|x|^{3\al+5}\sign x}
{(\al+1)^3(2\al+3)(3\al+5)}\right)y(x).
\end{array}
$$
Процесс последовательной регуляризации можно неограниченно продолжать.
На каждом шаге необходимо добавлять очередное слагаемое в формулу для
$y_2$. Заметим, что точки $\al=-2+\frac1n$, $n=1,2,\dots$ являются
исключительными -- метод не позволяет определить оператор для таких
$\al$.

Опишем изложенный метод в общем виде. Запишем уравнение Риккати
$$
-w'+w^2+|x|^\al=0
$$
и будем искать его решение в виде ряда
\begin{equation}
\begin{array}{ccc}
w(x)=a_1|x|^{\al+1}\sign x+a_2|x|^{2\al+3}\sign x+a_3|x|^{3\al+5}\sign x+
\dots,\\\phantom{.}\\
\al\ne-2+\frac1n,\ n=1,2,\dots,\ \al>-2.
\end{array}
\label{2.I}
\end{equation}
Тогда мы получим уравнения на коэффициенты
$$
a_1=\frac1{\al+1},\ a_2=-\frac{a_1^2}{2\al+3},\ \dots,\
a_n=\frac1{n\al+n+1}\suml_{k=1}^{n-1}a_ka_{n-k}.
$$
Так как $\al>-2$, то степени у членов ряда (\ref{2.I}) $\al+1$,
$2\al+3$,...
растут, что доказывает его равномерную сходимость на $[-1,1]$ в силу
того, что $|a_n|\leq C\frac{(n/2)!}{(n-1)!}$. Кроме того, найдется номер
$N$ такой, что степень $N\al+2N-1>0$. Обозначим остаток ряда, начиная с
этого члена через $r_N(x)$, а начальную часть ряда $w(x)-r_N(x)=:w_N(x)$.
Если теперь в резольвентном уравнении (\ref{00.I}) осуществить переход к
системе с помощью замены
$$
y_1(x)=y(x),\ y_2(x)=y'(x)-w_N(x)y,
$$
то матрица системы (\ref{1.I}) будет иметь вид
\begin{equation}
\begin{pmatrix}
w_N&\quad1\\-\la+r_N'+2wr_N-r_N^2&\quad-w_N
\end{pmatrix}.
\label{4.I}
\end{equation}
При этом все ее элементы лежат в $L_1$.

Перейдем теперь к конструкции оператора \(L(\alpha)\).
В соответствии с описанной регуляризацией, назовем первой квазипроизводной
функцию
$$
y^{[1]}(x):=y'(x)-w_N(x)y(x),
$$
где функция $w_N(x)$ определена выше.

Регуляризуем дифференциальное выражение (\ref{de1}), а именно,
заменим его выражением
\begin{equation}
\widetilde{l}(y)=-(y^{[1]}(x))'-w_N(x)y^{[1]}(x)+
(r'_N(x)+2w_N(x)r_N(x)-r^2_N(x))y(x).
\label{nex.I}
\end{equation}
Везде в дальнейшем мы будем использовать выражение (\ref{nex.I}) вместо
(\ref{de1}).

Максимальный оператор $L_M$, связанный с дифференциальным
выражением (\ref{nex.I}) определим равенствами
$$
\begin{array}{ccc}
L_My=\widetilde{l}(y),\\\phantom{.}\\
\mathfrak D(L_M)=\{y(x)|\, y(x),\, y^{[1]}(x)\in W_1^1[-1,1];\,
\widetilde{l}(y)\in L_2[-1,1]\}.
\end{array}
$$
Минимальный оператор $L_m$ определим как сужение оператора $L_M$
на область
$$
\mathfrak D(L_m)=\{y|\, y(\pm1)=y^{[1]}(\pm1)=0\}.
$$
\begin{Proposition}
Оператор $L_m$ симметречен и замкнут, его индексы дефекта равны
$\{2,2\}$, а любое его самосопряженное расширение $L$ является сужением
оператора $L_M$ на область
$$
\mathfrak D(L)=\{y|\,y\in\mathfrak D(L_M), U_1(y)=U_2(y)=0\},
$$
где линейные формы имеют представление
\begin{equation}
U_j(y)=a_{j1}y(0)+a_{j2}y^{[1]}(0)+b_{j1}y(1)+b_{j2}y^{[1]}(1)=0,\quad
j=1,2,
\label{r1.I}
\end{equation}
для коэффициентов которых выполнены равенства~\eqref{eq:T.1.5}.
\end{Proposition}
\begin{proof}
Поскольку элементы матрицы (\ref{4.I}) системы (\ref{1.I})
лежат в пространстве $L_1$, для системы верна теорема существования и
единственности решения задачи Коши (см. \cite[\S 16]{Naj}). Далее
доказательство получается дословным повторением доказательств теорем 1.1 и
1.2 из \S\,1.
\end{proof}
\begin{Note}
Поскольку функция $w_N(x)$ является гладкой в окрестности точек $\pm1$, в
формулах для расширений (\ref{r1.I}) квазипроизводная
может быть заменена обычной производной. При этом коэффициенты форм $U_j$
изменятся, но соотношения \eqref{eq:T.1.5} сохранятся.
\end{Note}

\begin{Note}
Выбрав краевые условия в предложении 4.4 условиями Дирихле $y(-1)=y(1)=0$,
получим оператор, совпадающий с оператором $L(\al)$, построенным
методом аналитического продолжения для потенциала
$|x|^\al$ (см. \S 4.3 пункт {\bf 2}).
\end{Note}
Действительно, из аналитичности матриц \eqref{4.I} по
параметру $\al$ и классической теоремы из курса обыкновенных
дифференциальных уравнений следует аналитическая зависимость решений
однородного уравнения
$$
-y''+|x|^{\al}y=0
$$
по $\al$. Это означает, что фундаментальная система решений однородного
уравнения совпадает с функциями $V_1(x)$, $V_2(x)$.

В заключение остановимся на вопросе о возможности приближения построенных
операторов операторами Штурма--Лиувилля с гладкими потенциалами. В \S\,1.2 мы
показали, что имеет место равномерная резольвентная сходимость операторов
\(L_{\varepsilon}\) с гладкими потенциалами \(q_{\varepsilon}(x)\) к оператору
\(L\) с потенциалом \(q\in W_2^{-1}\), если \(q_{\varepsilon}\to q(x)\) при
\(\varepsilon\to 0\) в пространстве \(W_2^{-1}\).

В нашем случае потенциалы $|x|^\al$ не лежат в пространстве $W_2^{-1}$ при
$\al<-\frac32$, а значит теорема не гарантирует возможность приближения
построенных операторов. Более того, в \S\,1.5 был приведен пример,
показывающий, что различные последовательности операторов с потенциалами,
сходящимися в пространстве $W_p^{-1}$ при $p<2$ могут сходиться к разным
пределам или не сходиться вовсе. Тем не менее можно описать класс
последовательностей гладких потенциалов, таких, что операторы Штурма--Лиувилля
с этими потенциалами будут приближать построенные операторы в смысле сильной
резольвентной сходимости.

Рассмотрим матрицы
$$
A_\eps=\begin{pmatrix} w_{N,\eps}&1\\-\la+v_{N,\eps}&-w_{N,\eps}
\end{pmatrix},
$$
где функции $w_{N,\eps}$ и $v_{N,\eps}$ гладкие и выполнены условия
$$
w_{N,\eps}\to w_N,\quad v_{N,\eps}\to r_N'+2wr_N-r_n^2
$$
в пространстве $L_1$ при $\eps\to0$. От системы
\begin{equation}
\begin{pmatrix} y_{1,\eps}\\y_{2,\eps}\end{pmatrix}'=A_\eps
\begin{pmatrix} y_{1,\eps}\\y_{2,\eps}\end{pmatrix}-\begin{pmatrix}
0\\f\end{pmatrix}
\label{l.I}
\end{equation}
можно совершить обратный переход к уравнению второго порядка
$$
-y''_\eps+(v_{N,\eps}+w_{N,\eps}'+w_{N,\eps}^2)y_\eps=\la y_\eps+f,
$$
где $y_\eps=y_{1,\eps}$.
В силу теоремы о непрерывной зависимости решения системы дифференциальных
уравнений от параметра, решения системы (\ref{l.I}) будут сходиться и мы
получим сильную резольвентную сходимость операторов
$L_\eps=-\dfrac{d^2}{dx^2}+v_{N,\eps}+w_{N,\eps}'+w_{N,\eps}^2$ к
построенному оператору $L$.

\medskip

{\bf Савчук А.~М.:}\hfill\break
Московский государственный университет,\hfill\break
Механико--математический факультет,\hfill\break
119899 Москва, Россия\hfill\break
email: $\rm artem\_savchuk\symbol{64}mail.ru$

\medskip
{\bf Шкаликов А.~А.:}\hfill\break
Московский государственный университет,\hfill\break
Механико--математический факультет,\hfill\break
119899 Москва, Россия\hfill\break
email: $\rm ashkalikov\symbol{64}yahoo.com$

\end{document}